\documentclass{amsart}
\usepackage[utf8]{inputenc}

\usepackage{kotex,amsmath,amsthm,verbatim,amssymb,amsfonts,amscd,graphicx,calc}
\usepackage{xcolor}
\usepackage{graphics}

\usepackage{euscript, enumerate,hhline,array}
\usepackage{url,hyperref}

\newcommand{\p}{\partial}

\newcommand{\la}{\langle}
\newcommand{\ra}{\rangle}

\newcommand{\e}{\varepsilon}

\newcommand{\eps}{\varepsilon}

\newcommand{\be}{\begin{equation}}
\newcommand{\ba}{\begin{aligned}}
\newcommand{\bee}{\begin{equation*}}
\newcommand{\ee}{\end{equation}}
\newcommand{\ea}{\end{aligned}}
\newcommand{\eee}{\end{equation*}}
\newcommand{\bea}{\begin{equation} \begin{aligned} }
\newcommand{\eea}{\end{aligned}\end{equation} }

\newtheorem{manualtheoreminner}{Theorem}
\newenvironment*{manualtheorem}[1]{%
  \renewcommand\themanualtheoreminner{#1}%
  \manualtheoreminner
}{\endmanualtheoreminner}

\newtheorem{manualcorollaryinner}{Corollary}
\newenvironment*{manualcorollary}[1]{%
  \renewcommand\themanualcorollaryinner{#1}%
  \manualcorollaryinner
}{\endmanualcorollaryinner}

\theoremstyle{plain}
\newtheorem{theorem}{Theorem}[section]

\newtheorem{corollary}[theorem]{Corollary}

\newtheorem{lemma}[theorem]{Lemma}
\newtheorem{proposition}[theorem]{Proposition}

\theoremstyle{remark}

\newtheorem{remark}[theorem]{Remark}

\theoremstyle{definition}
\newtheorem{definition}[theorem]{Definition}

\numberwithin{equation}{section}

\title{Translating surfaces under flows by sub-affine-critical powers of Gauss curvature}

\author{Beomjun Choi}

\address{BC: Department of Mathematics, POSTECH, 77 Cheongam-ro, Nam-gu, Pohang, Gyeongbuk 37673, Republic of Korea}
\email{bchoi@postech.ac.kr}

\author{Kyeongsu Choi}
\address{KC: School of Mathematics, Korea Institute for Advanced Study, 85 Hoegiro, Dongdaemun-gu, Seoul 02455, Republic of Korea}
\email{choiks@kias.re.kr}

\author{Soojung Kim}
\address{SK: Department of Mathematics, Soongsil University,  Seoul 06978,    Republic of Korea}
\email{soojungkim@ssu.ac.kr}

\date{}

\begin{document}

\begin{abstract} We classify the surfaces translating under the flows by sub-affine-critical powers of the Gauss curvature. This, in particular, lists all translating solitons possibly model Type II singularities for convex closed solutions in all positive powers.   The surfaces are entire graphs, and therefore our result corresponds to the Liouville theorem for the degenerate Monge--Amp\`ere  equations $\det D^2 u=(1+|Du|^2)^{2-\frac{1}{2\alpha}}$ on $\mathbb{R}^2$ in the range $0< \alpha <1/4$. The result also reveals that the moduli spaces of solutions are homeomorphic to either Euclidean spaces or cylinders.     
\end{abstract}
 
\maketitle 
\tableofcontents

\section{Introduction}

A one-parameter family of smooth embeddings $X:M^2 \times [0,T)\to \mathbb{R}^3$ with complete convex images $\Sigma_t :=X(M^2,t)$ is a solution to the $\alpha$-Gauss curvature flow  ($\alpha$-GCF) if 
\begin{equation}
\tfrac{\partial}{\partial t} X=-K^\alpha \nu 
\end{equation}
holds, where $K$ and $\nu$ denote the Gauss curvature and  the unit outward-pointing\footnote{The convex hull of the image surface $\Sigma_t$ is the inside region, when $\Sigma_t$ is non-compact.} normal vector at $X$, respectively.

\bigskip

When $\alpha=\frac14$, the $\frac{1}{4}$-GCF\footnote{The $\alpha$-Gauss curvature flow with $\alpha=\frac{1}{n+2}$ is the affine normal flow in  $\mathbb{R}^{n+1}$.} is called the affine normal flow, because a solution remains as the flow is invariant under any volume-preserving affine transformation of the ambient space $\mathbb{R}^3$. This affine invariant property yields a beautiful result in singularity analysis that every closed solution to the $\frac{1}{4}$-GCF converges to an ellipsoid at the singularity after rescaling; see \cite{andrews1996contraction}. In the super-affine-critical case $\alpha>\frac{1}{4}$, there is a strong rigidity in singularity formation that every closed solution to the $\alpha$-GCF  converges to a
round sphere at the singularity after rescaling;  see \cite{andrews1999gauss, andrews2011surfaces, andrews2016flow,  brendle2017asymptotic, choi-daskalopoulos16, guan2017entropy}.

\bigskip

In stark contrast to the above, the theory for sub-affine-critical powers $\alpha \in (0,\frac{1}{4})$ is rudimentary. For example, multiple closed shrinkers have been discovered for small $\alpha \ll 1$ in \cite{andrews2000motion}, but we do not know yet if it is the full list of closed shrinkers\footnote{The shrinkers under the $\alpha$-curve shortening flow were already classified by Andrews \cite{andrews2003classification}.}. Moreover, the experts conjectured that closed solutions to the $\alpha$-GCF with $\alpha \in (0,\frac{1}{4})$ generically develop type II singularities, namely closed solutions would not converge to shrinkers in general ; see \cite{andrews2002non}. Therefore, one needs to characterize translators to study generic singularities. Indeed, we often obtain the result that blow-up limits at  Type II singularities of flows satisfying certain reasonable conditions must be translators; see \cite{brendle2019uniqueness, choi2020type}.
\bigskip

The main result of this paper is the characterization of all translating surfaces to the $\alpha$-Gauss curvature flow with $\alpha\in (0,\frac{1}{4})$. We say that a complete\footnote{We say that a surface is complete and convex if it is the boundary of a convex set with non-empty interior, so called a convex body.  } non-compact convex smooth surfaces $\Sigma\subset \mathbb{R}^3$ is a smooth translator moving with the velocity $\omega \in \mathbb{R}^3$ if
\begin{equation}\label{eq:translator-curvature}
K^\alpha=\langle -\nu, \omega\rangle 
\end{equation}
holds for some $0\not = \omega\in \mathbb{R}^3$. In this paper, let us fix $\omega=\bf{e_3}$ as other translators can be obtained by rescalings and rigid motions. Then, by the authors' recent result in \cite{CCK_regularity}, every translator $\Sigma$ is the graph of an entire convex function $u:\mathbb{R}^2\to \mathbb{R}$ with a certain polynomial growth rate.
\bigskip

An $\alpha$-GCF translator in $\mathbb{R}^{n+1}$, with $\omega=\bf{e_{n+1}}$ and $\alpha>0$, has the \textit{profile} function $u$ over the open set $\Omega\subset \mathbb{R}^n$ defined by $\Sigma \subset \overline{\Omega}\times \mathbb{R}$ and $(x,u(x))\in \Sigma$ for all $x\in \Omega$. Then, the convex function $u$ is a solution to the Monge--Amp\`ere equation
\begin{equation}\label{eq:u_in_R^n}
\det D^2 u=(1+|Du|^2)^{\frac{n+2}{2}-\frac{1}{2\alpha}}
\end{equation}
satisfying  $\displaystyle \lim_{x\to x_0} |Du(x)|=+\infty$ for every $x_0 \in \partial \Omega$ if $ \Omega \neq \mathbb{R}^n$.

\bigskip

Urbas \cite{Urbas1988, U98GCFsoliton} showed that if $\alpha>\frac{1}{2}$, the domain $\Omega=\text{Int}(\pi(\Sigma))\subset\mathbb{R}^n $ is bounded. Conversely, given a convex bounded open domain $\Omega \subset \mathbb{R}^n$ and $\alpha>\frac{1}{2}$, there exists a unique translator $\Sigma$ under the $\alpha$-GCF which is asymptotic to $\partial\Omega \times \mathbb{R}$. Moreover, he proved that $\Sigma$ is strictly convex and smooth in $\Omega\times \mathbb{R}$. However, if $\Sigma$ touches $\partial \Omega\times \mathbb{R}$, it might be neither strictly convex nor smooth on $\partial \Omega\times\mathbb{R}$ as shown in \cite{CDL}. In general, the translators attract nearby non-compact flows; see \cite{choi2018convergence, choi2019convergence}. Moreover, the translator is the only non-compact ancient solution in each $\Omega\times \mathbb{R}$; see \cite{choi2020uniqueness}.   

\bigskip

 Urbas' results \cite{Urbas1988, U98GCFsoliton} can be considered as the well-posedness of the optimal transport  $Du$ with respect to the cost function $c(x,y)=-x\cdot y$ between the bounded domain $\Omega$ equipped with the standard measure $dx$ and the unbounded domain $\mathbb{R}^n$ equipped with a finite measure $(1+|y|^2)^{\frac{1}{2\alpha}-\frac{n+2}{2}}dy$; c.f. the celebrated result of Ma-Trudinger-Wang \cite{ma2005regularity} dealt with optimal transports between two bounded domains with general cost functions and finite measures. 

\bigskip

The measure $(1+|y|^2)^{\frac{1}{2\alpha}-\frac{n+2}{2}}dy$ over $\mathbb{R}^n$ becomes an unbounded measure for $\alpha\leq \frac{1}{2}$ and thus $\Omega$ should be unbounded. Moreover, the optimal transport is not unique. H.Jian-X.J.Wang  \cite{JW} constructed infinitely many entire solutions to \eqref{eq:u_in_R^n} for $\alpha\leq \frac{1}{2}$, and also showed $\Omega$ should be $\mathbb{R}^n$ for $\alpha <\frac{1}{n+1}$. To be specific, they showed that given an ellipsoid $E\subset\mathbb{R}^n$, there exists a solution $u$, satisfying $u(0)=0$ and $Du(0)=0$, such that $E$ is homothetic to the John ellipsoid of the level set $\{x\in \mathbb{R}^n:u^*(x)=1\}$, where $u^*$ is the Legendre transformation of $u$; c.f. K.S.Chou-X.J.Wang \cite{chou1996entire}.
\bigskip

Note that the affine critical case $\alpha=\frac{1}{n+2}$ corresponds the simple equation $\det D^2 u=1$, and the entire convex solutions are quadratic polynomials \cite{nitsche1956elementary, jorgens1954losungen, calabi1958improper, cheng1986complete, caffarelli1995topics}; see also \cite{caffarelli2003extension, caffarelli2004liouville, jin2014liouville, jin2016solutions, li2019monge} for more Liouville theory of Monge--Amp\`ere equations. Therefore, in this affine critical case, the list of solutions constructed by H.Jian-X.J.Wang matches to the full list of entire solutions. However, our result will show that there are more solutions in the case $\alpha <\frac{1}{9}$ and $n=2$; see Theorem \ref{thm:classification-1}. 

\subsection{Main results}

 In this subsection, we describe the main result concerning the Liouville theorem for the translating surfaces. Meanwhile, we will recall our previous works regarding (i) the regularity and the estimate on the growth rate of solutions in \cite{CCK_regularity}, and (ii) the construction of solutions with prescribed asymptotic behavior \cite{CCK_existence}, which will indeed describe all possible solutions.

First, note that the translating surfaces which appear as singularity models are not a priori smooth, and hence the definition of translating surface should be understood in a generalized sense. The main result of our previous work \cite{CCK_regularity}  reduces the classification of weak translators, defined in a geometric Alexandrov sense as Urbas \cite[pp 272]{U98GCFsoliton}, into those surfaces which are smooth entire graphs having homogeneous growth rate:
 
 \begin{manualtheorem}{A}[Growth rate \cite{CCK_regularity}]\label{thm-smooth} Let $\Sigma\subset \mathbb{R}^3$ be a translator with $\alpha \in (0,1/4)$ in a geometric Alexandrov sense. Then, its profile $u: \mathbb{R}^2\to\mathbb{R}$ is an entire smooth strictly convex function such that  
\begin{equation}\label{eq-440}
C^{-1}\le  \liminf_{|x|\to \infty}  |x|^{-\frac{1}{1-2\alpha}}u(x) \le\limsup_{|x|\to \infty}  |x|^{-\frac{1}{1-2\alpha}}u(x)  \le C
\end{equation} 
holds for some uniform constant $C>1$ depending only on $\alpha$.
\end{manualtheorem}

Toward the classification, the next step is to understand the possible finer asymptotics of solutions. In view of the blow-up (or blow-down) technique in geometric analysis, this amounts to investigate the limit(s) of rescaled profiles  
\be \label{eq-translatorgraph} u_{\lambda}(x) = \lambda ^{-\frac{1}{1-2\alpha}}u(\lambda x), \ee
as $\lambda\to \infty$.   As $u(x)$ solves  $ \det D^2 u = (1+|Du|^2)^{2-\frac1{2\alpha}}$, any reasonable limits of $u_\lambda$ will be at least a weak solution to the blow-down equation \begin{equation}\label{eq-asymptranslatorgraph}
\det D^2 \bar u = |D\bar u|^{4-\frac1{\alpha}}.
\end{equation}

Note that the limit $\bar u$ is a homogeneous solution $\bar u (x) = |x|^{\frac{1}{1-2\alpha}}g(x/|x|)$ to the blow-down equation, and it is completely determined by its level curve at height $1$, $\{x\in\mathbb{R}^2\,:\, \bar  u(x)=1\}$. One of our crucial observation is that the level curve of the homogeneous solution has to be a shrinking solitons to the $\frac{\alpha}{1-\alpha}$-curve shortening flow; see Remark \ref{remark-shrinker}. Our first result is the proof of uniqueness for the blow-downs, indicating that each translator has a unique homogeneous solution that is asymptotic to it near infinity.

\begin{theorem}[Uniqueness of blow-downs]\label{thm:tangent_flow-shrinker} For $\alpha\in (0,\frac 14)$, let $\Sigma\subset\mathbb{R}^3$ be a given translator with the height function $u(x)$. Then the rescaled functions $u_\lambda (x):=\lambda^{-\frac{1}{1-2\alpha}} u(\lambda x)$ converge to a homogeneous function $\bar u(x)=r^{\frac{1}{1-2\alpha}}g(x/|x|)$  in $ C^{\infty}_{loc}(\mathbb{R}^2 \setminus \{0\}) \cap C^0 _{loc}(\mathbb{R}^2)$ as $\lambda\to +\infty$. Furthermore, the blow-down $\bar u$ solves \eqref{eq-asymptranslatorgraph} and the level curve $\{x\in \mathbb{R}^2:\bar u(x)=1\}$ is a convex embedded smooth shrinker under the $\frac{\alpha}{1-\alpha}$-curve shortening flow\footnote{The $\alpha$-curve shortening flow is the $\alpha$-GCF in $\mathbb{R}^2$, and we will denote this by $\alpha$-CSF.}.
\end{theorem}

This limit blow-down $\bar u$ reminds the tangent flow at infinity of an ancient mean curvature flow or the tangent cone at infinity of a complete minimal hypersurface. Indeed, the self-similarity of blow-down level curves $\{ \bar u=l\}$ is an analogue to the self-similarity in tangent flows and tangent cones.

\begin{figure}
\centering
\def\svgscale{.7}{
\begingroup%
  \makeatletter%
  \providecommand\color[2][]{%
    \errmessage{(Inkscape) Color is used for the text in Inkscape, but the package 'color.sty' is not loaded}%
    \renewcommand\color[2][]{}%
  }%
  \providecommand\transparent[1]{%
    \errmessage{(Inkscape) Transparency is used (non-zero) for the text in Inkscape, but the package 'transparent.sty' is not loaded}%
    \renewcommand\transparent[1]{}%
  }%
  \providecommand\rotatebox[2]{#2}%
  \newcommand*\fsize{\dimexpr\f@size pt\relax}%
  \newcommand*\lineheight[1]{\fontsize{\fsize}{#1\fsize}\selectfont}%
  \ifx\svgwidth\undefined%
    \setlength{\unitlength}{352.15235066bp}%
    \ifx\svgscale\undefined%
      \relax%
    \else%
      \setlength{\unitlength}{\unitlength * \real{\svgscale}}%
    \fi%
  \else%
    \setlength{\unitlength}{\svgwidth}%
  \fi%
  \global\let\svgwidth\undefined%
  \global\let\svgscale\undefined%
  \makeatother%
  \begin{picture}(1,0.7580683)%
    \lineheight{1}%
    \setlength\tabcolsep{0pt}%
    \put(0,0){\includegraphics[width=\unitlength,page=1]{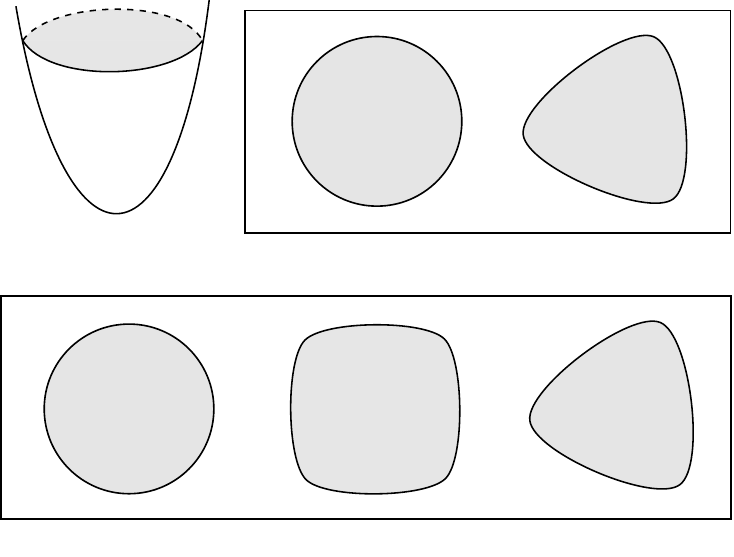}}%
    \put(0.4516907,0.39957521){\makebox(0,0)[lt]{\lineheight{1.25}\smash{\begin{tabular}[t]{l}$\frac{1}{16}\le \alpha <\frac 1 9$\end{tabular}}}}%
    \put(0.44997956,0.00571263){\makebox(0,0)[lt]{\lineheight{1.25}\smash{\begin{tabular}[t]{l}$\frac 1{25} \le \alpha <\frac 1{16}$\end{tabular}}}}%
    \put(0.65457146,0.5862513){\makebox(0,0)[lt]{\lineheight{1.25}\smash{\begin{tabular}[t]{l}$\text{or}$\end{tabular}}}}%
    \put(0.32557378,0.18891638){\makebox(0,0)[lt]{\lineheight{1.25}\smash{\begin{tabular}[t]{l}$\text{or}$\end{tabular}}}}%
    \put(0.65619214,0.18979606){\makebox(0,0)[lt]{\lineheight{1.25}\smash{\begin{tabular}[t]{l}$\text{or}$\end{tabular}}}}%
  \end{picture}%
\endgroup%

\caption{possible blow-downs in different ranges of $\alpha$}\label{fig1}}
\end{figure}

\begin{remark}\label{cor:class.blow-down} The asymptotic behavior in Theorem \ref{thm:tangent_flow-shrinker} is complemented by the classification of  the self-shrinking $\frac{\alpha}{1-\alpha}$-CSFs, shown by Andrews \cite{andrews2003classification}. According to the classification, 
given $\alpha\in [\frac{1}{(m+1)^2},\frac{1}{m^2})$ with $2\leq m \in \mathbb{N}$, there exist exactly $(m-1)$ smooth homogeneous solutions to \eqref{eq-asymptranslatorgraph} up to rotations. These are the rotationally symmetric homogeneous solution and the $k$-fold\footnote{See Definition \ref{def:k-fold}.} symmetric homogeneous solutions with integers $k\in [3,m]$ (if $m\geq 3$).  See Figure \ref{fig1}.
\end{remark}

\bigskip

Following Theorem \ref{thm:tangent_flow-shrinker}, the next step is to classify translators for each prescribed blow-down at infinity. If we linearize   \eqref{eq-asymptranslatorgraph} around a fixed blow-down, then the kernel of linearized operator (i.e. the Jacobi fields) accounts possible variations among translators. In Section \ref{subsection-ODE}, we investigate the Jacobi fields by using the spectrum of linearized operator around the level curve shrinker\footnote{This reminds the study of minimal cone using the spectrum of Jacobi operator around the minimal surface made by the intersection between the cone and the sphere \cite{simon1986minimal}\cite{Chan1997CompleteMH}\cite{MR809969}.}. Among the Jacobi fields, effective ones that actually generate translators are slowly decaying ones.

In the previous work \cite{CCK_existence}, we showed the existence of translators that correspond to the slowly decaying Jacobi fields of a fixed blow-down at infinity. Note this notion of slowly decaying Jacobi fields was previously employed in  \cite{simon1986minimal,Chan1997CompleteMH} to study the complete minimal surfaces asymptotic to a fixed minimizing cone at infinity.

\begin{manualtheorem}{B}[Existence \cite{CCK_existence}]\label{thm:classification-1} For $\alpha\in(0, 1/4)$, let  $m:=\lceil \alpha^{-1/2}\rceil -1 $.\footnote{Recalling the definition $\lceil x\rceil :=\min \{n\in \mathbb{Z}\,:\, x\le n \}$, $m$ is the unique integer such that $\frac{1}{(m+1)^2}\le \alpha < \frac{1}{m^2}$.  
} There exists a  $(2m+1)$-parameter family of smooth translators under the $\alpha$-Gauss curvature flow  whose blow-downs are the rotationally symmetric homogeneous solution to \eqref{eq-asymptranslatorgraph}. In addition, for each integer $k\in [3,m]$ (which necessarily implies $\alpha<1/9$), there exists a ${(2k-1)}$-parameter family of smooth translators whose blow-downs are a given $k$-fold symmetric homogeneous solution to \eqref{eq-asymptranslatorgraph}. {Moreover, the parametrization is continuous in the topology of locally uniform convergence of surfaces.}
\end{manualtheorem}

\begin{table}[ht] 
  \centering{\renewcommand{\arraystretch}{2.5}
\begin{tabular}{ |c|c|c|c|c| c} 
 \hline 
  & \includegraphics[scale=.7]{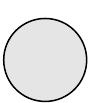}& \includegraphics[scale=.67]{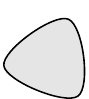} & \includegraphics[scale=.7]{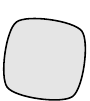} & \includegraphics[scale=.7]{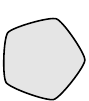}  & 
$\ba \cdots \\ \\\ea $ \\ 
\hhline{|=|=|=|=|=|=} 
 $\frac 19 \le \alpha <\frac 14$ & 5 & X & X  & X & $\cdots$  \\ 
\hline
 $\frac 1{16} \le \alpha <\frac 19$ & 7 & 5 & X & X & $\cdots$ \\ 
 \hline
  $\frac 1{25} \le \alpha <\frac 1{16}$ & 9 & 5 &7 & X& $\cdots$  \\ 
 \hline
  $\frac 1{36} \le \alpha <\frac 1{25}$ & 11 & 5 & 7 & 9 & $\cdots$ \\ 
 \hline
 $\vdots$ & $\vdots$ & $\vdots$ & $\vdots$ & $\vdots$ & $\ddots $ \\ 
\end{tabular}
}
\vspace{.4cm}

\caption{the number of parameters for fixed blow-down homogeneous solutions}
\label{table1}
\end{table}

\begin{remark} We explain in Section \ref{subsec:existenceconstruction} the behavior and the estimate of solutions constructed in \cite{CCK_existence} to an extent that is required in this paper. Let us also note that each family is invariant under actions of translating solution in $\mathbb{R}^3$. In other words, the family contains the translated copies of each solution.  In  the asymptotically rotationally symmetric case, the family is invariant under rotations about $x_3$-axis. This last fact is, however, not shown in the construction, but a corollary of classification theorem. See Table \ref{table1} for a summary on the number of parameters and shapes depending on the ranges of $\alpha$. 

\end{remark}

Now we can state our main classification result. 

\begin{theorem}[Classification] \label{thm:classification-2} For $\alpha\in (0,\frac{1}{4})$, a given translating surface under the $\alpha$-Gauss curvature flow must be one of the smooth translators given in Theorem \ref{thm:classification-1}.
\end{theorem}

For a sake of simplicity, we will say radially symmetric profile function as $\infty$-fold symmetric. Then, by the unique blow-down Theorem \ref{thm:tangent_flow-shrinker}, the profile $u$ of a translator $\Sigma$ is asymptotic to a $k$-fold homogeneous solution to \eqref{eq-asymptranslatorgraph} for some $k\in \{\infty,3,4,\cdots\}$. Hence, for each $\alpha \in (0,\frac{1}{4})$  we can consider the space 
\[\mathcal{X}^{k}=\{\text{translators whose profiles have $k$-fold symmetric blow-downs\},}\] where $\mathcal{X}^k$ is equipped with the topology of uniform convergence on compact sets. The last main theorem exhibits the topology of $\mathcal{X}^k$. 

\begin{theorem}[Topology of moduli space]\label{thm:moduli.space.topology}
If $\alpha \in [\frac{1}{9},\frac{1}{4})$, the space $\mathcal{X}^{\infty}$ is homeomorphic to $\mathbb{R}^5$ and every translator belongs to $\mathcal{X}^{\infty}$.  If $\alpha\in [\frac{1}{(m+1)^2},\frac{1}{m^2})$ for some  integer $ m\ge 3$, then $\mathcal{X}^{\infty}$ is homeomorphic to $\mathbb{R}^{2m+1}$, and $\mathcal{X}^k$ with $k\in \{3,\ldots,m\}$ is homeomorphic to $ \mathbb{R}^{2k-1}\times\mathbb{S}^1$. Moreover, every translator belongs to one of $\mathcal{X}^{\infty},\mathcal{X}^3,\cdots,\mathcal{X}^m$.
\end{theorem}

Indeed, the proof of Theorem \ref{thm:classification-2} shows a bijection between $\mathcal{X}^k$ and $ \mathbb{R}^{2k-1}\times\mathbb{S}^1$ (or between $\mathcal{X}^\infty$ and $\mathbb{R}^{2m+1}$, $m=\lceil \alpha^{-1/2}\rceil -1$,) and the continuity statement of Theorem \ref{thm:classification-1} shows the continuity of this bijection on one direction. Our contribution here is to establish the continuity of the other direction. 

\bigskip

{Finally, we note that the main Theorems \ref{thm:tangent_flow-shrinker}, \ref{thm:classification-2},  and \ref{thm:moduli.space.topology} can be shown {\it mutatis mutandis} for blow-down translators which are graphs of (Alexandrov) solutions to \eqref{eq-asymptranslatorgraph}. Corresponding results of Theorem \ref{thm-smooth} and Theorem \ref{thm:classification-1} were previously shown for blow-down translators: See \cite[Theorem 1.2 and Corollary 1.4]{CCK_regularity} and \cite[Remark 5.3]{CCK_existence}, respectively. The families of blow-down translators share the same properties of the families of translators (in terms of their number of parameters and asymptotic behaviors). A blow-down translator, however, is only smooth away from the point whose subdifferential contains the zero vector.}

\subsection{Outline}
{In Section \ref{sec-prelim}, we set up the main problem and introduce some preliminaries.  As in \cite{CCK_existence}, we work with the support function of the level curves: $S(l,\cdot)$ denotes the representation of level curve at level $x_3=l$ by the support function 
\[S(l,\theta)= \sup_{p\in \Sigma \cap \{x_3=l\}} \langle  (\cos\theta,\sin\theta,0),p \rangle .\] The classification theorem is a consequence of fine asymptotic analysis on the dynamics of the support function $S(l,\theta)$ as $l$ goes to infinity. As it solves an elliptic equation, the translator equation can be viewed as a second order nonlinear ODE in function space when the solution is parametrized in terms of its levels. A novelty of this work lies in Section \ref{subsection-ODE} where we vectorize the equation and transform it as a first order ODE system so as to apply dynamical system arguments at ease. We introduce the inner product and spectral theory of this space which is closely related with a  suitable function  space adapted to the linearization of the problem. Finally in Section \ref{subsec:existenceconstruction}, we remind the main existence theorems of \cite{CCK_existence} to an extent that is needed in the classification: see Theorem \ref{thm-existence} and Corollary \ref{cor-wexpression}.

\bigskip

The remaining Sections \ref{sec-spectral}--\ref{sec-uniq-bd2} are dedicated to the study of the classification and topology of the moduli space, and prove Theorem \ref{thm:tangent_flow-shrinker}, Theorem \ref{thm:classification-2}, and Theorem \ref{thm:moduli.space.topology}.  Here, we enumerate the logical steps towards the theorems. As a preliminary, by Theorem \ref{thm-smooth}, note every translator $\Sigma$ is the graph of an entire function $u(x)$ with a specific growth rate $|x|^{\frac{1}{1-2\alpha}}$.

\begin{enumerate}[$(1)$]
\item At large heights, the translator becomes similar to the graph of a homogeneous profile. More precisely, for any $\lambda_i \to \infty$, $u_{\lambda_i}(x)= \lambda_i^{-\frac{1}{1-2\alpha} } u(\lambda_i x)$ subsequentially converges to a $({1-2\alpha})^{-1}$-homogeneous function whose level curve is a shrinking soliton to the $\frac{\alpha}{1-\alpha}$-CSF. In terms of $S(l,\theta)$,  
\[S_{\lambda_i} (l,\theta) = \lambda_i ^{-(1-2\alpha)}S(\lambda_i l,\theta) \longrightarrow \frac{l^{1-2\alpha}}{1-2\alpha} h(\theta)  \] along a subsequence, where $h$ represents the support function of a shrinker as in  \eqref{eq-shrinker}. We call the limit convex hypersurface represented by $(1-2\alpha)^{-1}l^{1-2\alpha}h(\theta)$ a blow-down of $\Sigma$. At this stage, one may obtain {\it possibly} different shrinkers along different rescaling parameters $\lambda_i\to \infty$. 

\item The limits $h(\theta)$ along different $\lambda_i\to \infty$ in step (2) may only differ by rotations. i.e.  for a given $\Sigma$, there is a shrinker $h=h(\Sigma)$ such that for any $\lambda_i \to \infty$, there is a subsequence such that
\[S_{\lambda_i} (l,\theta) = \lambda_i ^{-(1-2\alpha)}S(\lambda_i l,\theta) \longrightarrow \frac{l^{1-2\alpha}}{1-2\alpha} h(\theta-\theta_0)\text{ for some }\theta_0 \in [0,2\pi].   \]
  
\item Assuming (2), the limit homogeneous profile is indeed unique (i.e. $\theta_0$ is unique for all $\lambda_i\to \infty$) and hence full convergence holds in step (2). Moreover, there is an improved convergence rate which is exponential in $s=\log l$: for some $\e =\e (\alpha)>0$,
\[S(l,\theta)- (1-2\alpha)^{-1} l^{1-2\alpha} h(\theta) = o(l^{1-2\alpha-\e}) \text{ as }l\to \infty.\]  
 
 \item Assuming (3), the given translator must be one of translators constructed  in Theorem \ref{thm-existence} in Section \ref{subsec:existenceconstruction}. 

\item The space of translating solitons that have the same limiting homogeneous profile forms a topological manifold which is homeomorphic to an Euclidean space whose dimension is equal to the number of parameters in Table \ref{table1}.  

\end{enumerate}

For smoother exposition, the paper will follow a different order. Assuming the consequence of (2) are given (we call it the unique limit shape property; see Definition \ref{def-41}),  we prove (3)-(4) in Section \ref{sec-spectral}.  We prove the unique limit shape (1)-(2) in Section \ref{sec-uniq-bd1} and verify the topology of moduli space (5) in Section \ref{sec-uniq-bd2}.

\bigskip

Let us explain rough ideas. In Section \ref{sec:set.up.spect}, we first introduce preliminary notations. Once the equation is vectorized, we analyze the behavior of a translator for large $s$ by studying the dynamical system for the projections of a given solution onto three spaces spanned by the eigenvectors of $\mathcal{L}$ in \eqref{eq-mathcalL}. In Section \ref{sec:norotation}, (3) is shown for the case of non-radial limiting shrinking curves; see Theorem \ref{thm-exponential1}. We prove it by showing the amount of extra rotation required to perfectly fit a given translator to the limit homogeneous profile from one scale to the next scale decays as a geometric sequence converging to zero.  This is a reminiscence of the celebrated technique of Allard-Almgren \cite{AA} in showing the uniqueness of tangent cones of minimal surfaces. As a corollary, we also obtain that solution converges to the limit exponentially fast in $s$.

In Section \ref{sec-exponential-round}, (3) is shown for the case of the radial shrinking curve ($h$ is constant); see Theorem \ref{thm-exponential2}. As the unique limit shape property vacuously implies the uniqueness of blow-down in this case, the theorem is to verify an exponential convergence rate.  The only possibility of non-exponential decay is when the projection onto the eigenvectors that tend to grow like homogeneous profiles, called neutral modes, dominates the solution. We show there is no such case by a contradiction and this follows from a similar argument of \cite{ChoiSun} in the (backward) parabolic case.

In Section \ref{sec-uniqueness}, we show (4), the classification under the assumption that solutions have unique limit shapes as in Definition \ref{def-41}; see Theorem \ref{thm-72}. The family of translators in Theorem \ref{thm-existence} is parametrized by $\mathbf{y}\in \mathbb{R}^K$ so that the asymptotic difference of two given translators is described by the last non-zero coordinate of difference of their parameters; see {\it (ii)} in Theorem \ref{thm-existence}. This asymptotic difference is represented by a Jacobi field. Conversely in Theorem \ref{thm-72}, we prove the difference of two translators is asymptotic to a kernel of linearized equation, a Jacobi field. The theorem follows by successive applications of the previous observation between the given translator and a translator in the constructed family: at each application, one finds another translator in the family which is closer to the given translator near infinity. Once we reach two translators converging to each other at sufficiently fast exponential rate, then the two must be identical by the comparison principle and then the classification follows.

 In Section \ref{sec-rigidityhomogeneous}, we shall prove Theorem \ref{thm-rigiditywithcondition}, a lemma that will be used in the proof of Theorem \ref{thm:tangent_flow-shrinker} (uniqueness of blow-down) in Section \ref{sec-uniq-bd1}. We present the theorem here as the proof uses similar ideas which are already introduced in the other subsections of Section \ref{sec-spectral}. This is one main reason  why the steps (3)-(4) precede   (1)-(2).

In Section \ref{sec-uniq-bd1}, we directly prove Theorem \ref{thm:tangent_flow-shrinker}. The steps (1)-(2) will be presented in the course of the proof. Then Theorem \ref{thm:classification-2} simply follows by collecting previous results. The proof of Theorem \ref{thm:tangent_flow-shrinker} mainly uses the function $J_{u^*}(\cdot)$ in Definition \ref{def:J_v} where $u^*$ is the Legendre dual of $u$. This function behaves like an entropy which detects a self-similar behavior of a solution at infinity. For instance, given a solution $v$  to $\det D^2 v = |x|^{\frac{1}{\alpha}-4}$, $J_{v}(\cdot)$ is a constant if $v$ is a homogeneous function, and conversely (and more importantly) $v$ is a homogeneous function if $J_v(\cdot)$ attains an interior maximum or minimum. In view of Theorem \ref{thm-smooth} and a compactness argument, the proof reduces to obtain the convergence of $J_{u^*}(x)$ as $|x|\to \infty$. This convergence follows by a contradiction argument. Assuming the $\limsup$ and the  $\liminf$ are different, we may pick a sequence of scaling factors along which the blow-down is not a homogeneous function, but at the same time satisfies the conditions of Theorem \ref{thm-rigiditywithcondition} which asserting that the blow-down has to be a homogeneous function. 

Finally in Section \ref{sec-uniq-bd2}, we prove the space of translators with a prescribed blow-down, if it is equipped with the topology of locally uniform convergence, is homeomorphic to Euclidean space. As mentioned after Theorem \ref{thm:moduli.space.topology}, it suffices to prove the continuity of the bijection from the space of translators to Euclidean space. This follows by showing a sequence of translators in Theorem \ref{thm-existence} with unbounded parameters cannot converge to another translator having the same blow-down (Lemma \ref{lem-para.bound}). In this proof, the quantitative estimate for the constructed translators from Corollary \ref{cor-wexpression} in Section \ref{subsec:existenceconstruction} is used. }

\section{Preliminary and Notation}\label{sec-prelim}

 \subsection{Change of coordinates and entropy} \label{section-2.2}

   Given a translator $\Sigma$, let $u(x)$ be smooth strictly convex entire height function given by Theorem \ref{thm-smooth}.  For each $l> \inf u$, we may consider the support function of the level convex closed curve $\{u(x)=l\}\subset \mathbb{R}^2$  {at height $l$}  as follows: for $(l,\theta)\in  \{l\in \mathbb{R}: l >\inf u\}\times \mathbb{S}^1$, 
  \bea \label{eq-supportfunction} 
  S(l,\theta) := \sup_{\{u(x)=l\}} \la x, e^{i\theta} \ra  .
  \eea 
The support function $S(l,\theta)$ is smooth on $ (\inf u,\infty)\times \mathbb{S}^1$, and it becomes a solution to the equation 
\begin{equation}\label{eq-translatorsupport}
 S_{\theta\theta}+S+(1+S_l^{-2})^{\frac1{2\alpha}-2}S_l^{-4}S_{ll} =0 .
\end{equation} 
If $u(x)$ solves the blow-down equation \eqref{eq-asymptranslatorgraph}, then $S$ is a solution to the equation
\bea \label{eq-asymptranslatorsupport}
 S_{\theta\theta}+S +S_{ll}S_l^{-\frac1\alpha} =0.
 \eea 
See \cite[Appendix A]{CCK_existence} for the derivation of two equations above.

Next, we recall a crucial quantity from \cite[Definition 5.4]{DSavin}. This function, which was originally defined for solutions to the dual equation $\det D^2 v(y)= |y|^{\frac1\alpha -4}$, behaves like a sort of entropy, which detects a self-similar behavior of a given solution. Later, we use it to show each translator has a unique limit shape (see Definition \ref{def-41}) in Section \ref{sec-uniq-bd1}. 

\begin{definition} \label{def-J*}
Given a strictly convex function $u$, we define a function $J^*_u$ on $\{x\,:\, Du(x) \neq 0\}$ by
\begin{equation} \label{eq-0566}
J^*_u: =   { u^{kk}}{ ( u^{ij}u_iu_j)^{4\alpha-1}}, 
\end{equation}
where the matrix $\{u^{ij}\}$ denotes the inverse matrix of $D^2u$ and $u^{kk}=u^{11}+u^{22}$.
  We will denote    $J^*_u$  simply by $J^*$  when no confusion is  possible on $u$.
\end{definition}

\begin{lemma}\label{prop-J^*-supp}
In terms of $S $,   we have the following expression for $J^* $:
 \begin{equation}\label{eq-J^*-supp}
  J^* =   \left[  S_l    (-S_{ll})^{-1}\right]^{4\alpha }  [    S_l^2+S_{\theta l}^2   - S_{ll} (S_{\theta\theta}+S)].
\end{equation}
%
\end{lemma}

 \begin{proof} In view of Corollary A.1 in \cite{CCK_existence}, where $Du$ and $D^2u$ were derived in terms of $S$, we have if $Du \neq 0$ and  $Du /|Du|=(\cos \theta,\sin\theta)$, then 
$D u = \begin{bmatrix}  (S_l)^{-1} &   0 \end{bmatrix} R_{\theta}^{\,T}$ and
\[[D^2  u]^{-1}  =    \frac{S_l}{-S_{ll}}R_{\theta}^{\, T}\begin{bmatrix}  S_{l}^2  & S_l S_{l\theta} \\ S_l S_{l\theta} & S_{l\theta}^2 -S_{ll}(S+S_{\theta\theta}) \end{bmatrix} R_{\theta} .\] Here  $Du$ is a row matrix and $R_\theta$ is the matrix representing the rotation by the angle $\theta$.  
Since in the matrix notation
\[J^*_u= \mathrm{tr}([D^2u]^{-1})(Du [D^2u]^{-1} Du^T)^{4\alpha-1},\] the expression \eqref{eq-J^*-supp} directly follows.
 \end{proof}

 \medskip

The next remark reveals an important connection between the  asymptotic behavior of   translating solitons for  the $\alpha$-GCF at infinity and the shrinking solitons for  the curve shortening flow whose speed is given by the $\frac{\alpha}{1-\alpha}$th power of its curvature, namely $\frac{\alpha}{1-\alpha}$-CSF. 

 \begin{remark}\label{remark-shrinker}
  Let $\bar u(x)$ be a solution to the blow-down equation \eqref{eq-asymptranslatorgraph} that is a homogeneous function. Then the degree of homogeneity has to be $1/\sigma$ with 
\be
 \sigma:=1-2\alpha.
 \ee
The corresponding support function $S(l,\theta)$ has to be $\sigma$-homogeneous in $l$ and we may write
 $ 
 S(l,\theta)=   {\sigma}^{-1}l^{\sigma} h(\theta)$ for some $h(\theta)$. 
 In light of \eqref{eq-asymptranslatorsupport},    $h :  \mathbb{S}^1\to \mathbb{R}$ becomes a solution to  \bea\label{eq-shrinker}  h_{\theta\theta}+h=  2\alpha \sigma h^{\frac{\alpha-1}{\alpha}}.\eea  Indeed, this $h(\theta)$ represents the support function of a shrinking soliton.
If  the  support function  $h$ of a  convex  curve $\Gamma$ satisfies \eqref{eq-shrinker}, then   
 \bea  \kappa^{\frac{\alpha}{1-\alpha}}=(2\alpha \sigma )^{\frac{\alpha}{\alpha-1}} \langle x, \nu  \rangle,\eea
 where $\kappa$  and     $\nu $ denote    the curvature and   the unit outward-pointing normal vector.  
Namely, 
  the level curves of the homogeneous solution $\bar u$     to \eqref{eq-asymptranslatorgraph} are shrinking solitons for the curve shortening flow whose speed is given by the $\frac{\alpha}{1-\alpha}$-CSF. 
  \end{remark}
 
\medskip
 
  In Section \ref{sec-uniq-bd1}, we show    that  the blow-down of the  translator is a $\sigma^{-1}$-homogeneous solution to the blow-down equation \eqref{eq-asymptranslatorgraph}. i.e. the level curves of  the  translator become asymptotic to a shrinking soliton (after suitable rescaling). Note that the entropy $J^*$ is a constant for homogeneous solution.
\medskip

\begin{remark}\label{rmk-homo-J}
If $S(l,\theta)=\sigma^{-1} l^\sigma h(\theta)$ for some $h$ satisfying \eqref{eq-shrinker}, then the entropy \eqref{eq-J^*-supp} becomes
\begin{equation}\label{eq-J^*}
J^*= (2\alpha )^{-4\alpha}    \left(h^2+h'^2+ 4\alpha^2 h^{2-\frac{1}{\alpha}}\right).
\end{equation}
Moreover, \eqref{eq-shrinker} implies
\begin{equation}
\partial_\theta J^* = 2(2\alpha )^{-4\alpha} h'   \left(h+h''- 2\alpha \sigma h^{1-\frac{1}{\alpha}}\right)=0.
\end{equation}
Thus, $J^*=J^*(l,\theta)$ is a constant if  $S$ is a $\sigma$-homogeneous solution to \eqref{eq-asymptranslatorsupport}.

\end{remark}

 \medskip

Once we know a given translator must be asymptotic to a homogeneous solution, the classification of homogeneous solutions becomes essential. Indeed, B. Andrews \cite{andrews2003classification}  classified all possible shrinking solitons for  the  CSF and thereby the homogeneous solutions are also classified.

\begin{definition}[$k$-fold symmetry]\label{def:k-fold} 
Let $\Gamma\subset \mathbb{R}^2$ be    a closed embedded curve  that contains the origin in its inside. 
We say $\Gamma$ is $k$-fold symmetric if   there is  the largest positive integer   $k$  such that $\Gamma$ is invariant under the rotation by angle $\frac{2\pi}{k}$.\end{definition} Similarly,    an $l$-homogeneous function $u(x)=|x|^l g(x/|x|)$ on $\mathbb{R}^2$ is  said to be   $k$-fold symmetric if its level curves are $k$-fold symmetric, i.e., $g(\theta+ \frac{2\pi}{k}) = g(\theta)$ for all $\theta\in \mathbb{S}^1$ and there is no larger integer with the same property. According to the definition, circles  are not $k$-fold symmetric for any integer $k$.  
    Unless  the curve  $\Gamma$ in Definition \ref{def:k-fold}
  is a circle, there is a  unique positive integer $k$.    From the classification of shrinkers in \cite{andrews2003classification}, the following corollary follows. 
 
 \begin{corollary}[Classification of homogeneous solutions. c.f. \cite{andrews2003classification}]\label{cor-anderews2003} 
 	Given $\alpha \in (0,1/4)$,    let $m$ be the largest integer satisfying $m^2 < \tfrac 1\alpha$. 
  Then modulo rotations in $x_1x_2$-plane, there are exactly $(m-1)$ homogeneous solutions to \eqref{eq-asymptranslatorgraph} which are smooth away from the origin: they are the rotationally symmetric solution, and  (if $m\geq3$) $k$-fold symmetric solutions for  $3\leq k\leq m$.	  
 \end{corollary}

 \medskip
 
 In order to   study translators for the     $\alpha$-GCF,   we   consider the asymptotic expansion  of the support function  near infinity. Let   $h$ be  a solution to  the equation \eqref{eq-shrinker}. Suppose $S(l,\theta)$ represents a translator and it is asymptotic to homogeneous function $  {\sigma}^{-1} l^{\sigma} h(\theta)$, which is a solution to \eqref{eq-asymptranslatorsupport}, near $l=\infty$. 
To study the higher order asymptotics, we express the solution $S(l,\theta)$  as follows:
 \bea \label{eq-w}
 S(l,\theta)=  {\sigma}^{-1}   e^{\sigma s} h(\theta) + w(s,\theta) \quad \text{ for }  s= \ln l
 \eea 
with some function $w=w(s,\theta)$. In this paper, we often consider the equations of translators \eqref{eq-translatorsupport} and blow-down translators \eqref{eq-asymptranslatorsupport} at the same time. In fact, one can write down the two equations in one form: 
\begin{equation}\label{Se} 
S+S_{\theta\theta}+(\eta   +S_l^{-2})^{\frac1{2\alpha}-2}S_l^{-4}S_{ll} =0,
\end{equation} 
where $\eta\in \{0,1\}$. The next lemma shows the equations of $w$ for both cases.
 
 \begin {lemma} [Lemma A.1 in \cite{CCK_existence}] \label{lem-w-eq}
 Let $S(l,\theta)$ be a solution to \eqref{Se} with $\eta=0$ or $1$  and let 
  $w$ be a function given by  \eqref{eq-w}. 
Then  $w$   solves 
 \bea \label{eq-w'} 
 w_{ss}+w_s+{h}^{\frac1\alpha}( w_{\theta\theta}+w)   =E(w)=-E_1(w) -\eta E_2(w), 
 \eea 
 where 
  \bea\label{eq-Ew} 
 E_1(w):=  \tfrac1\alpha h^{\frac1\alpha-1} {e^{- \sigma s}}  {w_s}( w_{\theta\theta}+w) + ( 2\alpha  e^{\sigma s}  h^{1-\frac{1}{\alpha}}  +w_{\theta\theta}+   w )N_1(w), \eea 
 \begin{equation}
   N_1(w):= (h+{ e^{-\sigma s}}{w_s})^{\frac1\alpha} -h^{\frac1\alpha} -\tfrac1{\alpha} h^{\frac1\alpha-1}e^{-\sigma s}{w_s}, \text{ and}
 \end{equation}
 \bea \label{eq-Ew2}
&E_2(w):=( h+e^{-\sigma s}w_s )^{\frac1\alpha} ( 2\alpha  e^{\sigma s}  h^{1-\frac{1}{\alpha}}  +w_{\theta\theta}+   w )N_2(w),
\eea
\begin{equation}
     N_2(w):= \left (1+ e^{-4\alpha s} (  h+e^{-\sigma s}w_s )^2\right)^{2-\frac{1}{2\alpha}}-1.
\end{equation}

 \end{lemma}
 If a given translator is asymptotic to the homogeneous solution, then the normalized profile function $e^{-\sigma s}w$ and its derivatives converge to $0$ as $s\to \infty$. Note that the terms in $e^{-\sigma s}E_1$  are at least quadratic in (the derivatives of) $w$ and those in $e^{-\sigma s}E_2$ are exponentially decaying in $s$. Meanwhile, one can view $e^{-\sigma s}(w_{ss} +w_s+h^{\frac1 \alpha} (w_{\theta\theta}+w))$  as a linear operator applied to $e^{-\sigma s} w$. For this reason, we will call $w_{ss} +w_s+h^{\frac1 \alpha} (w_{\theta\theta}+w)$ the Jacobi operator and a kernel of this operator as a Jacobi field associated with the asymptotic homogeneous profile. 
 
 \subsection{ODE formulation} \label{subsection-ODE} 
 
   In this subsection, we consider a transformation of the   second order  elliptic equation  \eqref{eq-translatorsupport}  into    a system of   first-order  ODEs on Hilbert spaces, which  enables us  to apply the Merle-Zaag  ODE Lemma  \ref{lem-MZODE} in  its analysis.    
   
   Let  us first fix     a smooth positive solution  $h$   to  the equation \eqref{eq-shrinker}.  The   elliptic equation \eqref{eq-w'}  can be written as   a system of ODEs for the height variable $s$ as follows: 
\bea\label{eq-vectorw}
\frac{\p}{\p s}
\begin{bmatrix} w \\ \p_s w \end{bmatrix}=\mathcal{L}\begin{bmatrix} w \\ \p_s w \end{bmatrix}+\begin{bmatrix}0  \\ E(w) \end{bmatrix}, 
\eea 
where the operator $\mathcal{ L} $ on the variable $\theta\in \mathbb{S}^1$ is given by 
\bea \label{eq-mathcalL}
\mathcal{ L} \begin{bmatrix}w_1 \\w_2\end{bmatrix}  := \begin{bmatrix} w_2 \\- {h}^{\frac1\alpha} ( (w_1)_{\theta\theta} +w_1 ) -w_2  \end{bmatrix}.  
\eea 
To discuss a functional space, let us   consider the operator $L$    defined  by 
 \bea \label{eq-linearopL}
 L\varphi=:h^{\frac1{\alpha}}( \varphi_{\theta\theta}+\varphi) ,\eea
and the Hilbert function space $L^2_h(\mathbb{S}^1)$ which denotes the functions in $L^2(\mathbb{S}^1)$  equipped with the inner product
 \bea \label{eq-innerh} 
 ( f,g)_h := \int_{\mathbb{S}^1} fg\, h^{-\frac1\alpha}d\theta , \eea
which is  adapted to  the operator $L$.
\begin{lemma} \label{lemma-lambdavarphi}
 The operator  $L$ is self-adjoint with respect to    $ ( \cdot , \cdot )_h$. There is an orthonormal basis of $L^2_h$, namely $\{\varphi_j\}_{j=0}^\infty$, consists of eigenfunctions of $L$ with decreasing order of eigenvalues: $L\varphi_j=\lambda_j\varphi_j$, and 
\[\lambda_0\ge \lambda_1 \ge \lambda_2 \ge \ldots . \]  
Each eigenvalue has finite multiplicity and $\lambda_j\to -\infty$ as $j\to \infty$. 
\begin{proof}  This follows by a standard spectral decomposition theorem. See more details in \cite[Section 2.2]{CCK_existence}
\end{proof}
\end{lemma} In light of  \cite[Lemma 5]{andrews2002non}, we may choose 
\begin{equation} \label{eq-c0c1c2}
\varphi_0:=m_0 h, \quad  \varphi_1:= m_1 \cos \theta ,\quad \varphi_2:=m_2 \sin\theta,
\end{equation} 
for some normalizing positive constants $m_i$, $i=0,1,2$. The corresponding eigenvalues are $\lambda_0=2\alpha(1-2\alpha)>0$, and $\lambda_1=\lambda_2=0$. In the remaining of this paper, let us fix $\{\varphi_j\}_{j\ge0}$ as suggested above. 

\medskip 

 Regarding the operator  $\mathcal{L}$ given in \eqref{eq-mathcalL}, we define a symmetric bilinear form $\la\cdot ,\cdot  \ra_{0}$  by 
 \bea \label{eq-0norm}
 \la  \mathbf{f}, \mathbf{g}\ra_{0} := \int_{\mathbb{S}^1} \left\{   (f_1)_\theta ( g_1)_\theta -f_1g_1 + {h}^{-\frac1\alpha}f_2g_2 \right\} d\theta .\eea 
 for $\mathbf{f} =   (f_1,f_2), \mathbf{g} = (g_1,g_2) \in    H^1(\mathbb{S}^1)\times L^2(\mathbb{S}^1)$. 
One directly checks that $\mathcal{L}$  is self-adjoint with respect to $\la\cdot,\cdot\ra_0$ and $\mathcal{L}$ has the following eigenfunctions:  
 \bea \label{eq-basis}
 \left\{  (\varphi_j, \beta^\pm_j \varphi_j) \right\}_{j=0}^\infty 
 \eea  with eigenvalues $\beta^\pm_j$, respectively. 
Here,  $\beta^+_j $ and $\beta^-_j $   are two distinct real roots of       $ \beta^2+ \beta + \lambda_j =0 $,
\bea\label{eq-def-beta_j-pm}
 \beta^\pm_j := -  \frac 12 \pm  \frac{\sqrt {1- 4 \lambda_j  }}{2}.
\eea 
Note that   $\beta_0^+= -2\alpha$, $\beta_0^-= -\sigma$, $\beta_1^+=0=\beta_2^+$ and $\beta_1^-=-1=\beta_2^-$.  
Another way of interpretation of eigenfunctions above follows by observing $ w= e^{\beta^+_j s}\varphi_j $ and $ e^{\beta^-_j s}\varphi_j$, for $j\ge0$, are Jacobi fields as discussed at the end of Section \ref{section-2.2}. i.e.   
\be w_{ss} +w_s+h^{\frac1 \alpha} (w_{\theta\theta}+w)=0, \quad \text{for } w= e^{\beta^\pm_j s}\varphi_j .\ee

\medskip 

As we will see the symmetric bilinear form $\langle \cdot,\cdot \rangle_0$ is not positive definite, our goal is to modify it so to become an inner product. A direct computation shows that for $ \mathbf{f} =(f_1,f_2) \in H^1( \mathbb{S}^1)\times L^2( \mathbb{S}^1)$ and $j\ge0$,  
\bea\label{eq-orthogo-varphi}
 \la      \mathbf{f} ,   (\varphi_j, \beta^\pm_j \varphi_j)\ra_{0} =     (  f_1,  -L \varphi_j)_h  + \beta_j^\pm  (f_2 ,  \varphi_j)_h = -\lambda_j ( f_1,  \varphi_j)_h +  \beta_j^\pm  ( f_2, \varphi_j )_h .\eea 
Thus $ \left\{ (\varphi_j, \beta^\pm_j \varphi_j) \right\}_{j=0}^\infty $ is pairwise  orthogonal with respect to $\la\cdot ,\cdot\ra_0$ 
since  $\beta_j^{\pm} $ satisfy   $ \beta^2+ \beta + \lambda_j =0 $.   
The symmetric bilinear form $\la\cdot, \cdot\ra_0$, however, is not positive definite as the following  three eigenfunctions 
 \bea
   (h,-2\alpha h),\quad   (\cos\theta ,0),\quad    ( \sin \theta, 0) \eea 
 among \eqref{eq-basis} are of  negative  or  zero   with respect to $\la\cdot ,\cdot\ra_0$.      
 Indeed,  we have that 
\bea\label{eq-sym-bf-not}
 &\la  (\varphi_0, \beta^+_0 \varphi_0 ), (\varphi_0, \beta^+_0 \varphi_0 ) \ra _0= - 2\alpha(1-4\alpha)<0,\quad \la  (\varphi_i, \beta^+_i \varphi_i ), (\varphi_i, \beta^+_i \varphi_i ) \ra _0=0  \,\, \hbox{ for $i=1,2$.}
  \eea
We modify so that $\{(\varphi_i,\beta^\pm_i \varphi_i)\}_{i=0}^\infty $ still serves as an orthogonal basis (and thus $\mathcal{L}$  is self-adjoint).

  \begin{definition}[Inner product] \label{def-innerprod}
 For a  solution $h$ to  \eqref{eq-shrinker}, we define the inner product $\la \cdot,\cdot\ra_h$  on $H^1( \mathbb{S}^1)\times L^2( \mathbb{S}^1)$  as follows: for $\mathbf{f}=(f_1,f_2)$ and $\mathbf{g}=(g_1,g_2)$  $\in H^1( \mathbb{S}^1) \times L^2( \mathbb{S}^1)$,  
  \bea \label{eq-def-innerprod}
 \la \mathbf{f},\mathbf{g} \ra_h :=& \la \mathbf f, \mathbf  g \ra _{0} - 2 \, \frac{ \la  \mathbf f ,    ({h},-2\alpha {h}) \ra_{0}\, \la  \mathbf g , ( {h},-2\alpha {h})  \ra_{0}        }{\la ( {h},-2\alpha {h})   , ({h},-2\alpha {h})\ra_{0}} \\
 &+ \left\{ (f_1, \cos\theta  )_h+ (f_2, \cos\theta )_h\right\}  \left\{ (g_1,\cos\theta )_h+ (g_2, \cos\theta )_h\right\} \\
 &+ \left\{ (f_1, \sin\theta)_h+ (f_2, \sin\theta)_h\right\}  \left\{ (g_1, \sin\theta)_h+ (g_2, \sin\theta)_h\right\} .
 \eea  
The  induced norm  is denoted   by $\Vert  \cdot     \Vert _h$ on   $ H^1( \mathbb{S}^1)\times L^2( \mathbb{S}^1)$.    
 \end{definition} 
 
 Here we summarize and prove the desired properties of $\la \cdot,\cdot \ra_h$. 
 \begin{proposition} \label{prop-normequivalence}
 For a  solution $h$ to  \eqref{eq-shrinker}, let $\la \cdot, \cdot \ra_h$ be given by  \eqref{eq-def-innerprod}. Then,   
 $(H^1(\mathbb{S}^1)\times L^2(\mathbb{S}^1),\la \cdot, \cdot \ra_h)$ is a  Hilbert space, and  the following statements hold true:
 \begin{enumerate}[(a)]
  \item The induced norm $\Vert \cdot  \Vert_h$ is equivalent to the standard norm on $H^1(\mathbb{S}^1)\times L_h^2(\mathbb{S}^1)$, i.e.,  for any $(f_1,f_2)\in H^1(\mathbb{S}^1)\times L^2(\mathbb{S}^1)$, 
 \be\label{eq-euiv-norm}
 C^{-1} \left(\Vert f_1\Vert _{H^1(\mathbb{S}^1)} +\Vert f_2 \Vert_{L^2_h(\mathbb{S}^1)} \right)\le \Vert (f_1,f_2)\Vert_h \le C \left(\Vert f_1\Vert _{H^1(\mathbb{S}^1)} +\Vert f_2 \Vert_{L^2_h(\mathbb{S}^1)} \right)
  \ee
  with some   constant $C=C(\alpha,h)>0$.
   \item $\{(\varphi_i ,\beta^\pm_i \varphi_i)\}_{i=0}^\infty$  	forms an orthogonal basis for $(H^1(\mathbb{S}^1)\times L^2(\mathbb{S}^1),\la \cdot, \cdot \ra_h)$.
 \item $(\varphi_i,\beta^\pm_i \varphi_i)$ are eigenfunctions to $\mathcal{L}$ with eigenvalues $\beta^\pm_i$, respectively.
  \item $\mathcal{L}$ on $\mathcal{D}\left(H^1(\mathbb{S}^1)\times L^2(\mathbb{S}^1)\right)$ is self-adjoint with respect to $\langle\cdot ,\cdot \rangle_h$.

 \end{enumerate}
  	 \end{proposition}
 
\begin{proof} {\it Step 1.} We first  prove that $\la \cdot, \cdot \ra_h$ is an inner product on $H^1(\mathbb{S}^1)\times L^2(\mathbb{S}^1)$ and the equivalence of norms \eqref{eq-euiv-norm} holds. It is immediate to see from the definition that $\langle\cdot ,\cdot  \rangle_h$ is symmetric bilinear and the second inequality in \eqref{eq-euiv-norm} holds for some $C=C(\alpha,h)$. It suffices to prove the first inequality in \eqref{eq-euiv-norm}, which in particular proves that $\langle\cdot ,\cdot \rangle_h $ is positive definite.

 Let  $\mathbf{f}$ be given function in $H^1(\mathbb{S}^1)\times L^2(\mathbb{S}^1)$. By viewing   $\mathbf{f}$ as functions in $L^2_h(\mathbb{S}^1)\times L^2_h(\mathbb{S}^1)$, we     decompose  $\mathbf{f}$  uniquely into
\bea
 \mathbf{f} &=  \sum_{i=0}^2 \left[ A^+_i (\varphi_i, \beta^+_i \varphi_i) + A^-_i (\varphi_i,\beta^-_i\varphi_i)\right ] + (v_1,v_2) .
\eea
Here, $v_i$ are orthogonal to $\{\varphi_0,\, \varphi_1, \,\varphi_2 \}$ in the space $L^2_h(\mathbb{S}^1)$, and $A^{\pm}_i$are coefficients which expand the remaining parts of $\mathbf{f}$.
Then  it follows that
\bea
\la \mathbf{f},\mathbf{f}\ra _0& = \la(v_1,v_2),(v_1,v_2) \ra _0 + \sum_{i=0} ^2 \left[|A^+_i| ^2  \la  (\varphi_i, \beta^+_i \varphi_i ), (\varphi_i, \beta^+_i \varphi_i ) \ra _0+ |A^-_i|^2  \la  (\varphi_i, \beta^-_i \varphi_i ), (\varphi_i, \beta^-_i \varphi_i ) \ra _0\right].  \eea 

Since  $(\varphi_1,\cos\theta)_h=\Vert \cos\theta \Vert _{L^2_h(\mathbb{S}^1)}  $, $(\varphi_1,\sin \theta)_h=\Vert \sin \theta \Vert _{L^2_h(\mathbb{S}^1)} $ and \eqref{eq-sym-bf-not},
we have 
\bea \label{eq-156}
\la \mathbf{f},\mathbf{f}\ra _h & = \la(v_1,v_2),(v_1,v_2) \ra _0 + \sum_{i=0} ^2 |A^-_i|^2  \la  (\varphi_i, \beta^-_i \varphi_i ), (\varphi_i, \beta^-_i \varphi_i ) \ra _0    \\ & \quad + |A^+_0|^2 \left[-  \la(\varphi_0, \beta^+_0 \varphi_0 ),(\varphi_0, \beta^+_0 \varphi_0 )\ra _0 \right] +| A^+_1|^2 \Vert \cos\theta \Vert _{L^2_h(\mathbb{S}^1)}^2 + |A^+_2|^2\Vert \sin\theta \Vert _{L^2_h(\mathbb{S}^1)}^2 .
\eea 
 This in particular implies that 
\be\label{eq-est-f-v_i-0}
\la \mathbf{f},\mathbf{f} \ra_h \ge \la (v_1,v_2),(v_1,v_2)\ra _0 +C^{-1} \sum _{i=0}^2 ( |A^+_i |^2 + |A^-_i|^2 ) 
\ee for some $C=C(\alpha,h)>0$.  Since $v_1$ is orthogonal to $\{h, \, \cos\theta, \, \sin\theta\}$, basis for non-negative eigenspace of $L$, we have a Poincare type estimate: 
 \bea\label{eq-est-v_i-0}
  \la (v_1,v_2),(v_1,v_2) \ra _0 &= \int_{\mathbb{S}^1}\left\{ |\partial_\theta v_1|^2 - |v_1|^2 + |v_2|^2 h^{-\frac 1\alpha} \right\}  d\theta  =  \int_{\mathbb{S}^1} \left(-v_1  Lv_1 + v_2^2\right)h^{-\frac1 \alpha } \, d\theta  \\
  &\ge- \lambda_3 \Vert v_1\Vert_{L^2_h(\mathbb{S}^1)}^2 +\Vert v_2\Vert_{L^2_h(\mathbb{S}^1)}^2,
   \eea 
where $\lambda_3<0$ is the eigenvalue from Lemma \ref{lemma-lambdavarphi}.   This implies 
   $\Vert v_1\Vert_{^2_h(\mathbb{S}^1)}^2 \le C \la (v_1,v_2),(v_1,v_2) \ra_0 $, and hence  
   \be
    \Vert (v_1,v_2) \Vert _{H^1(\mathbb{S}^1) \times L^2_h(\mathbb{S}^1)}^2 =\int_{\mathbb{S}^1}\left\{| \partial_\theta v_1 |^2 + |v_1|^2 + |v_2|^2 h^{-\frac 1\alpha}\right\} \, d\theta \le C\la (v_1,v_2),(v_1,v_2)\ra _0. 
    \ee
 Then it follows from \eqref{eq-est-f-v_i-0} that 
      \bea 
\la \mathbf{f},\mathbf{f} \ra_h &\ge  C^{-1}\left[  \Vert (v_1,v_2) \Vert _{H^1(\mathbb{S}^1) \times L^2_h(\mathbb{S}^1)}^2 + \sum _{i=0}^2 ( |A^+_i |^2 + |A^-_i|^2 )\right] 
\eea 
for some constant $C=C(\alpha,h)>0$.

\medskip

 {\it Step 2.}  Next, we will prove  (b). By polarizing the identity \eqref{eq-156},  one sees  that  $\{(\varphi_i ,\beta^\pm_i \varphi_i)\}_{i=0}^\infty$ is orthogonal with respect to $\la \cdot, \cdot\ra _h$ due to their orthogonality in $\la\cdot ,\cdot \ra _0$. 
To show it forms a basis,  suppose that  $  (f_1,f_2)  \in H^1( \mathbb{S}^1)\times L^2( \mathbb{S}^1)$ satisfies 
\be\label{eq-inn-pr}
 \la     (f_1,f_2) , (\varphi_i, \beta^\pm_i \varphi_i)   \ra_{h}=0\qquad \forall i=0,1,2,\ldots.
\ee 
 In view of \eqref{eq-orthogo-varphi},   it follows from  the  orthogonality of   $\{\varphi_i \}_{i=0}^\infty $ in $(\cdot, \cdot)_h$ that 
       \bea
 \la      (f_1,f_2) , (\varphi_0, \beta^\pm_0 \varphi_0)   \ra_{h}
 = &   -\lambda_0(  f_1,  \varphi_0)_h  + \beta_0^{\pm}  ( f_2,  \varphi_0)_h \\
 &- 2\frac{\la  (\varphi_0, \beta^\pm_0 \varphi_0)  ,  ({h},-2\alpha {h})   \ra_{0}  }{\la   ({h},-2\alpha {h})     , ({h},-2\alpha {h})  \ra_{0}}    \left\{  -\lambda_0  (  f_1, h)_h + \beta_0^{+}  ( f_2,  h )_h\right\}.
 \eea  
 This   implies  
\be
  -\lambda_0(  f_1,  \varphi_0)_h  + \beta_0^{+}  ( f_2,  \varphi_0)_h=0,\qquad 
-\lambda_0(  f_1,  \varphi_0)_h  + \beta_0^{-}  ( f_2,  \varphi_0)_h=0.
\ee
Solving the above system of the equations for $ (  f_1, \varphi_0)_h$ and $ (  f_2, \varphi_0)_h$, it holds  that 
$
(  f_1, \varphi_0)_h=(  f_2, \varphi_0)_h=0. 
$
Arguing similarly with the use of \eqref{eq-inn-pr}, we obtain that 
\be
(  f_1, \varphi_i)_h=(  f_2, \varphi_i)_h=0\quad\forall i=0,1,2,\ldots,
\ee 
which implies $  (  f_1,   f_2)=(0,0).$ Therefore,  $\{(\varphi_i ,\beta^\pm_i \varphi_i)\}_{i=0}^\infty$            is    an orthogonal basis of  $\left(H^1( \mathbb{S}^1)\times L^2( \mathbb{S}^1),  \la     \cdot  , \cdot  \ra_{h} \right) $.

      \medskip

 {\it Step 3.}  as $\beta_i^{\pm} $ are solutions to   $ \beta^2+ \beta + \lambda_i =0 $, from the definition, $
  \mathcal{L} (\varphi_i, \beta^\pm_i \varphi_i)= \beta^\pm_i (\varphi_i, \beta^\pm_i \varphi_i)$ for all $i=0,1,2,\ldots$, and this shows (c).  Lastly,  we shall prove  that $\mathcal{L}$  is self-adjoint with respect to $\la\cdot,\cdot\ra_h$.  Let  $\mathbf{f}$  and $\mathbf{g}$  be any functions in $H^1(\mathbb{S}^1)\times L^2(\mathbb{S}^1)$. From the definition of $\mathcal{L}$ \eqref{eq-mathcalL}, the fact that $L$ is self-adjoint with respect to $(\cdot,\cdot)_h$, and $\cos \theta$ is a kernel of $L$, we have
 \be\label{eq-cL-saj}
  ([\mathcal{L}  \mathbf{f}]_1,\cos\theta)_h+   ([\mathcal{L}  \mathbf{f}]_2,\cos\theta)_h = - (Lf_1,\cos\theta)_h=- (f_1, L\cos\theta)_h =0 ,\ee
and the same for $\sin \theta$. 
Finally  \eqref{eq-cL-saj} and the facts that $\mathcal{L}$ is self-adjoint with respect to $\la\cdot,\cdot\ra_0$, and $(h,-2\alpha h)$ is an eigenfunction to $\mathcal{L}$,   we have  $  \la\mathcal{L}  \mathbf{f},\mathbf{g}\ra_h = \la \mathbf{f},\mathcal{L}  \mathbf{g}\ra_h , $ which proves (d).  
\end{proof}

Lastly,  we introduce an    $L^2$-norm  and several  H\"older norms  which will be used  to study the asymptotic behavior at infinity. 

\begin{definition}[Norms on annuli and exterior domains]\label{def:norm.abbr}
Let $w=w(s,\theta)$ be a function  over    $I \times  \mathbb{S}^1$ for an interval $I\subset \mathbb{R}$.  Let   $k$ be  a non-negative integer,   let $0<\beta<1$ be a    constant and  let   $t$ and  $a$  be    constants such that $[t-a,t+a]\subset I$. 
Let us define 
\bea \Vert w\Vert_{L^2_{t,a}} := \Vert w \Vert_{L^2([t-a,t+a]\times \mathbb{S}^1)}, \eea 
\bea  
\Vert w \Vert _{C^{k}_{t,a}}:= \Vert w  \Vert_{C^{k}([t-a,t+a]\times \mathbb{S}^1)} ,\quad \hbox{and}\quad  \Vert w \Vert _{C^{k,\beta}_{t,a}}:= \Vert w  \Vert_{C^{k,\beta}([t-a,t+a]\times \mathbb{S}^1)}  . 
 \eea

For    a function    $w $      over       $[R,\infty)\times  \mathbb{S}^1$ with a  constant  $R>0$ and  for  a constant     $\gamma\in\mathbb{R}$,
we      define a $\gamma$-weighted $C^{k,\beta}$ norm on          $[R,\infty)\times  \mathbb{S}^1$ as follows: 
\bea \label{eq-1566}\Vert w \Vert_{C^{k,\beta,\gamma}_R} =  \sup_{t\ge R+1}} e^{-\gamma t}\Vert w \Vert_{C^{k,\beta}_{t,1}  .  \eea 
We denote  by $C^{k,\beta,\gamma}_R$  the   function   space   equipped with the  norm    $\Vert \cdot \Vert_{C^{k,\beta,\gamma}_R} $.
\end{definition}

\subsection{Family of translating surfaces} \label{subsec:existenceconstruction} We remind the main results of \cite{CCK_existence} concerning the construction of translating solitons and thereby provide details of Theorem \ref{thm:classification-1}. As discussed after Lemma \ref{lem-w-eq}, it is natural to view 
\be  \label{eq-Jacobi}
w_{ss}+w_s+{h}^{\frac1\alpha}(w_{\theta\theta}+w)   =0 
\ee
as the Jacobi equation (linearized equation.) The Jacobi fields, solutions to the Jacobi equation \eqref{eq-Jacobi}, serve as ansatzs for $w(s,\theta)$ near infinity. Indeed, $\varphi_j e^{\beta^+_j s}$ and $\varphi_j e^{\beta^-_j s}$ are Jacobi fields that have crucial role in the construction and the classification theorems.    

We say the rate    $ \beta^+_j $ corresponds to `slow decay' if $\beta^+_j< \sigma=1-2\alpha$.   
  The dimension of `slowly decaying' Jacobi fields  is denoted by 
  \bea \label{eq-K}
  K:= \# \left\{j\in \mathbb{N}\cup\{0\}:  \beta^+_j < \sigma \right\}.
\eea  
Here,    the term  `slowly decaying'    was previously used in \cite{Chan1997CompleteMH} for minimal hypersurfaces   asymptotic to minimal cones near infinity.   
Let $\mathcal{J}$ denote  the  space of  the   slowly decaying Jacobi fields: 
\be
 \mathcal{J}:=   \Big\{\,\sum_{j=0}^{K-1} y_j \,  l^{\beta^+_j} \varphi_j(\theta)  \Big\}_{ (y_0,\ldots,y_{K-1})\in \mathbb{R}^{K} } . 
 \ee
We will construct the translators which roughly correspond to the space of slowly decaying Jacobi fields $\mathcal{J}$. 
 
  \begin{remark}\label{remark-K}
 The dimension $K$ can be expressed in terms of $\alpha$ and $h$ as follows.
 \begin{enumerate}[(a)]
     \item If $h$ is radial, then $K = 2m+1 $ with $m=\lceil \alpha^{-1/2}\rceil -1$.  
     \item If  $h$ is $k$-fold symmetric, then $K  = 2k-1.$
 \end{enumerate}
 In fact, 
if $h$ is radial, i.e., $h\equiv (2\alpha\sigma)^{\alpha}$, then   $\lambda_j$ and $\beta^\pm_j$ are explicitly computed,  
and hence  we obtain that  
$ K = 2m+1 $ in the case.  Next, if  $h$ is $k$-fold symmetric,   utilizing the result  of   \cite[Theorem 1.2]{ChoiSun} implies that 
$ K  = 2k-1.$
 \end{remark}
\begin{manualtheorem}{C}[Existence, Theorem 2.9 in \cite{CCK_existence}]\label{thm-existence}
For given $\alpha \in (0,\frac{1}{4})$ and shrinker $h$, there exists $K$-parameter (see Remark \ref{remark-K}
 for $K$) family of translators, namely $\{\Sigma_{\mathbf{y}}\}_{\mathbf{y} \in \mathbb{R}^{K}}$,   whose level curves are asymptotic to the shrinker after rescalings.
 More precisely, the surfaces $\Sigma_{\mathbf{y}}$ and their support functions $S_{\mathbf{y}}$ defined in \eqref{eq-supportfunction} satisfy the followings:  for some  constant $\eps>0$ depending only on $\alpha$   and $h$, 
\begin{enumerate}[(i)]
\item  $\Sigma_{\mathbf{y}}$ is the graph of a smooth strictly convex entire function  and $S_{\mathbf{y}}$ has the asymptotics  
\bee  S_{\mathbf{y}}(l,\theta) = \sigma^{-1} { l^{\sigma} } h(\theta) + o (l^{\sigma-\eps })\quad \text{ as } l\to \infty.
\eee  
  \item For distinct vectors $\mathbf{\hat y}=(\hat y_0, \ldots,\hat y_{K-1} )$, $\mathbf{\bar y }=(\bar  y_0, \ldots,\bar y_{K-1} )$ in $\mathbb{R}^{K}$, let     $j=\max\{i\,: \, \hat y_i \neq \bar  y_i\}$ be the last parameter that does not coincide. 
Then there holds 
 \bee 
 S_{\mathbf{\hat y}}(l,\theta) -S_{\mathbf{\bar  y}}(l,\theta)= (\hat y_j-\bar y_j)l^{\beta^+_j} \varphi_{j}(\theta)+ o(l^{\beta^+_j-\eps}).
 \eee  

\item The first three parameters correspond to translations in $\mathbb{R}^3$:  for $\mathbf{y} \in\mathbb{R}^{K} $  and $a_j\in  \mathbb{R}$  ($j=0,1,2$),  
\bee  
\Sigma_{\mathbf{y}+(a_0,a_1,a_2,\mathbf{0}_{K-3} )}=  \Sigma_{\mathbf{y}} -a_0m_0 \mathbf{e}_3  +a_1m_1\mathbf{e}_1+a_2m_2\mathbf{e}_2,
 \eee 
 where $m_j$'s are the positive  normalizing constants given in \eqref{eq-c0c1c2}. 

\item The parametrization $\mathbf{y}\mapsto \Sigma_{\mathbf{y}}$ is continuous in the $C^0_{\text{loc}}$-convergence topology of surfaces.
\end{enumerate}
\end{manualtheorem}

\begin{remark} It is a direct consequence of $(ii)$ that  $\Sigma_{\mathbf{\hat y}}\neq \Sigma_{\mathbf{\bar y }}$ if $\mathbf{\hat y} \neq \mathbf{\bar y}$. \end{remark}

The construction is designed so that the proof shows more quantitative estimates than  $(i)$ and $(ii)$ as stated below. We need this result in Section \ref{sec-uniq-bd2} when we reveal the topology of moduli spaces. We need a few definitions.

 \begin{definition}[Effective coordinates] \label{def-b(a)}
We define a transformation    $\mathbf{b} : \mathbb{R}^K \rightarrow \mathbb{R}^K$  which maps $\mathbf{a}=(a_0,\ldots, a_{K-1})$ to $\mathbf{b}(\mathbf{a})=(b_{0},\ldots, b_{K-1})$ with  $b_i=\textrm{sgn}(a_i)|a_i|^{1/({\sigma-\beta^+_i})}$. 
\end{definition}
The $b_i$ are effective in the sense that $a_i l^{\beta^+_i-\sigma } \varphi_i$, a renormalization of Jacobi field according to the growth rate of solution, becomes comparable to unit scale at $l=|b_i|$ regardless of $i$ and the size of $a_i$. Since we work with $s=\ln l$, for each $\rho\in \mathbb{R}$ and $\mathbf{a}\in \mathbb{R}^{K}$, let us define the effective level as  
\[R_{\rho,\mathbf{a}}:= \rho + \ln(|\mathbf{b}(\mathbf{a})|+1).\]
We now state our construction with quantitative estimates.
\begin{manualcorollary}{D}[Corollary 5.1 in \cite{CCK_existence}]\label{cor-wexpression} Let $\{\Sigma_{\mathbf{y}}\}$ be the family constructed in Theorem \ref{thm-existence} and $w_{\mathbf{y}} $ be the representation of $\Sigma_{\mathbf{y}}$ as in \eqref{eq-w}. There are $\rho <\infty $ and $\varepsilon>0$ depending  on $\alpha $ and  $h $ such that 
the translator $\Sigma_{\mathbf{a}}$ satisfies the estimate  \be 
\Vert  w_{\mathbf{a}} \Vert_{C^{2,\beta,  \sigma -\eps  }_{R_{\rho,\mathbf{a}}} } \le       e^{\eps  R_{\rho, \mathbf{a}}} .
\ee
Moreover, there are families of functions $\{ g_{\mathbf{a}} \}_{\mathbf{a}\in\mathbb{R}^K}$ and $\{f_{\mathbf{a}}\}_{\mathbf{a} \in \mathbb{R}^K}$ {defined on on $s\ge R_{\rho,\mathbf{a}}$} such that 
 \begin{equation} \label{eq-wexpression-0} 
 w_{\mathbf{a}}= g_{\mathbf{0}}+\sum_{j=0}^{K-1}  a_j\big(  \varphi_j e^{\beta^+_j s} + g_{ \pi_j(\mathbf{a})} \big) +f_\mathbf{a}\qquad \hbox{on $s\ge R_{\rho,\mathbf{a}}$}.
 \end{equation}
Here,  $ \pi_j(\mathbf{a}) $  denotes $(  \mathbf{0}_{j}, a_j,\ldots, a_{K-1})$, the projection of $\mathbf{a}\in\mathbb{R}^K$ onto $\{\mathbf{0}_j\}\times \mathbb{R}^{K-j}$.
The families $\{g_{\mathbf{a}}\}$ and $\{f_{\mathbf{a}}\}$ have the following estimates: 
\begin{enumerate}[(i)] 
 
\item $\Vert f_\mathbf{a} \Vert _{C^{2,\beta,-\frac{1}{2} }_{R_{\rho,\mathbf{a}}}} \le {\frac{1}{2}  e^{(\sigma +{1}/{2} )R_{\rho,\mathbf{a}}}}$ ,

\item  $\| g_{\mathbf{0}} \|_{C^{2,\beta,\sigma -\eps}_{R_{\rho,\mathbf{a}}}}\ \le    e^{\eps  R_{\rho,\mathbf{a}}}$, and  $\| g_{ \pi_j(\mathbf{a})} \|_{C^{2,\beta,\beta^+_{j}-\eps}_{R_{\rho,\mathbf{a}}}}\ \le    e^{\eps  R_{\mathbf{a}}}$. 

\end{enumerate} 

\end{manualcorollary}

\section{Spectral dynamics around limit profile} \label{sec-spectral}

In the remaining sections, we prove the full classification  Theorem \ref{thm:classification-2}. In Section \ref{sec-uniq-bd1}, we will show that every translator has a blow-down which is unique modulo rotations. The goal in this section is to show the classification under this assumption. From now on, we say such solutions have unique limit shapes. Let us first write down a precise definition.

\begin{definition}[Unique limit shape] \label{def-41}
Let $\Sigma=\partial \{x_3>u(x)\}\subset \mathbb{R}^3$ be the graph of a striclty convex entire function $u:\mathbb{R}^2\to \mathbb{R}$. We denote by  $S$ the level support function  of $\Sigma$ as defined in \eqref{eq-supportfunction}. Then, we say $\Sigma$ has a \textit{unique limit shape} if there exists a function $h\in C^\infty(\mathbb{S}^1)$ such that for any divergent sequence of positive values $\{\lambda_i\}_{i\in\mathbb{N}}$, the sequence of rescaled functions $S_{\lambda_i}(l,\theta) = \lambda_i ^{-(1-2\alpha)}S(\lambda_il,\theta)$ has a subsequence converging to $\tfrac1{1-2\alpha} l^{1-2\alpha}h(\theta + \theta_0)$ in  the $C^{\infty}_{loc}(\{l>0\})$-topology  for some $\theta_0$.
\end{definition}
 
As in the previous sections, we denote  $\sigma :=1-2\alpha$, the frequently appearing constant. The first main result of this section proves the uniqueness of blow-downs and the fast convergence of solutions to the unique asymptotic.

\begin{theorem}[Unique blow-down and fast convergence]\label{thm-42}
 Let a smooth solution  $S $ to \eqref{Se}, with a fixed $\eta\in\{0,1\}$, have the unique limit shape represented by a  solution $h_0 \in C^\infty(\mathbb{S}^1)$ to \eqref{eq-shrinker}.
Then,  there exists  a solution $h\in C^\infty(\mathbb{S}^1)$ to \eqref{eq-shrinker} such that $h(\theta)=h_0(\theta+\theta_0)$ for some constant $\theta_0$ and 
\begin{equation}
    \lim_{\lambda\to +\infty}\lambda^{-\sigma} S(\lambda l,\theta)=\sigma^{-1}l^{\sigma}h(\theta)
\end{equation}
 holds in  the $C^{\infty}_{loc}(\{l>0\})$-topology. Moreover, there is some small constant $\e>0$ such that 
\bea \label{eq-1526}
S(l ,\theta)= \sigma^{-1} l^{\sigma} h(\theta)+o(l^{\sigma-\e})\quad \text{ as }l\to \infty.\eea 
\end{theorem} 
\begin{remark} 
 If the limit shape is radial, namely $h$ is a constant, then the uniqueness of the blow-down is vacuous.
\end{remark}

Now, by Theorem \ref{thm-42}, every translator has a unique blow-down. The next result shows that the translator has to be one of translators constructed in \cite{CCK_existence}, which is re-stated in Theorem \ref{thm-existence} in Section \ref{subsec:existenceconstruction}.

 \begin{theorem} [Classification]\label{thm-72} 
Suppose that a translator $\Sigma$, under the $\alpha$-GCF in $\mathbb{R}^3$ with $\alpha \in (0,\frac{1}{4})$, has the unique blow-down whose level set support function is $\sigma^{-1} l^{\sigma} h(\theta)$. Then \bea \Sigma \equiv \Sigma_{\mathbf{y}} \text{ for some } \mathbf{y} \in \mathbb{R}^{K},\eea
where $\{\Sigma_{\mathbf{y}}\}_{\mathbf{y}\in\mathbb{R}^{K}}$ denotes the $K$-parameter family of translators with the blow-down whose level set support function is $\sigma^{-1} l^{\sigma} h(\theta)$, which are constructed in Theorem \ref{thm-existence}.
\end{theorem}

\medskip

We here state a brief idea and an organization of subsections.  When $S(l,\theta)$ is close to an asymptotic profile $\sigma ^{-1} l^{\sigma}h $, we may use the ODE system approach of Section \ref{subsection-ODE} and consider the asymptotic behavior of a translator as the dynamics between the coefficients of eigenvectors in \eqref{eq-basis}.  As suggested by Merle-Zaag ODE Lemma (c.f. Lemma \ref{lem-MZODE}), the solution will be dominated either by the eigenspace of the eigenvalue $\sigma$ or by the eigenspace of the eivenvalues strictly less than $\sigma$.  In the later case, the solution satisfies a fast decay condition \bea \label{eq-decay11}   S(l,\theta )- \sigma^{-1} l^\sigma h(\theta )=o(l^{\sigma-\e} ) \text{ for some } \e>0\text{ as }l\to \infty ,\eea and this is the comparably easier case to prove the classification theorem. (See Section \ref{sec-uniqueness}.)

In Section \ref{sec:norotation} and \ref{sec-exponential-round}, we show that the former case can not happen. In Section \ref{sec:norotation} we prove it for non-radial limit shapes. In this case, the eigenspace of the eigenvalue $\sigma$ corresponds to the Jacobi fields generated by the  surface rotation around $x_3$-axis. In Section \ref{sec-exponential-round}, we prove it for the radial limit shape. In this case, such an eigenspace can exist only for certain powers  $\alpha^{-1}=k^{-2}$ where  $k\ge 3$ is an integer. Theorem \ref{thm-42} directly follows by combining two previous results and its proof is presented at the Section \ref{sec-exponential-round}.

Finally in Section \ref{sec-uniqueness}, we prove Theorem \ref{thm-72}. The fast decay \eqref{eq-decay11} is important to guarantee that the nonlinear term on the right hand side of \eqref{eq-w'} is negligible compared to the linear part. 

 Notice that in Section \ref{sec-rigidityhomogeneous} we show a rigidity theorem \ref{thm-rigiditywithcondition} by using Proposition \ref{prop-coeff-extract}, the main result in Section \ref{sec-uniqueness}. The rigidity theorem \ref{thm-rigiditywithcondition} will be employed in Sections \ref{sec-uniq-bd1} and \ref{sec-uniq-bd2} as a lemma.

  \subsection{Setting up spectral analysis}\label{sec:set.up.spect} We define some notations which will be used throughout this section. Given a smooth function $f(s,\theta)$, we define the vector-valued function $[f]$ by
  \begin{equation}
      [f](s,\theta):=(f(s,\theta),\partial_s f(s,\theta)).
  \end{equation}
Next, we recall the operator $L$ in \eqref{eq-linearopL} and its eigenfunctions $\varphi_j$ with eigenvalues $\lambda_j$ in Lemma \ref{lemma-lambdavarphi}. Then, by using $\beta^+_j$ in \eqref{eq-def-beta_j-pm}, we define
\begin{equation}
    \vec{\varphi}_j:= (\varphi_j,\beta_j^+ \varphi_j )\|(\varphi_j,\beta_j^+ \varphi_j )\|_h^{-1}.
\end{equation}
We notice that the operator $\mathcal{L}$ in \eqref{eq-mathcalL} satisfies
\begin{equation}
   \mathcal{L} \vec{\varphi}_j=\beta^+_j\vec{\varphi}_j.
\end{equation}
Remembering $\beta_0^+=-2\alpha$, given $\beta \geq \beta_0^+=-2\alpha$, we define the projections
\begin{align}
 & P_{>\beta}[f]:= \sum_{ \beta_j^+ > \beta } \langle  [f],\vec{\varphi}_j \rangle_h, &&  P_{=\beta}[f]:= \sum_{ \beta_j^+ = \beta } \langle  [f],\vec{\varphi}_j \rangle_h, &&   P_{<\beta}[f]:=[f]-P_{>\beta}[f]-P_{=\beta}[f].
\end{align}
Notice that $P_{<\beta}$ [resp. $P_{>\beta}$] is the projection to the space spanned by the eignfunctions whose eigenvalues are less [resp.  greater] than $\beta$. We also define $P_{\geq  \beta},P_{\leq \beta}$ in the same manner. In addition, we frequently use
\begin{align}
    & P_+=P_{>\sigma}, && P_\sigma=P_{=\sigma}, && P_-=P_{<\sigma}.
\end{align}
Also, given a smooth function $\mathbf{f}: \mathbb{S}^1\to \mathbb{R}^2$, we recall the standard $H^m\times H^n$-norm that
\begin{align}\label{def:H^m}
&   \|\mathbf{f}\|_{H^m\times H^n}^2=  \int_0^{2\pi}\sum_{i=0}^m |\partial_\theta^i f_1|^2+\sum_{i=0}^n|\partial_\theta^i f_2|^2\,  d\theta, \quad  \text{where}\; \mathbf{f}=(f_1,f_2).
\end{align}
We may abuse notation so that $H^0$ denotes $L^2$-norm, and  $\|\mathbf{f}\|_{H^m\times H^n}$ denotes $\|\mathbf{f}(s,\cdot)\|_{H^m\times H^n }$ if $\mathbf{f}$ is a function of $(s,\theta)\in [s_1,s_2]\times \mathbb{S}^1$. In what follows, we will mostly use the norms of type $H^{m+1} \times H^{m}$.

 \bigskip

  \begin{proposition}\label{prop:ODE_general}
 Suppose that smooth $\mathbf{f},\mathbf{g}: [s_1,s_2]\times \mathbb{S}^1\to \mathbb{R}^2$ satisfy $\partial_s \mathbf{f}=\mathcal{L}\mathbf{f}+\mathbf{g}$. Then, given $\beta \ge - 2\alpha $,
 \begin{align}  
  &\left|  \tfrac{d}{ds} \|P_{=\beta}\mathbf{f}\|_h -\beta \|P_{=\beta}\mathbf{f}\|_h  \right|  \leq   \| \mathbf{g} \|_{h}, \\
 & \tfrac{d}{ds}\|P_{\leq \beta}\mathbf{f}\|_h  \leq  \beta \|P_{\leq \beta}\mathbf{f}\|_h  + \| \mathbf{g} \|_{h},\\
  &  \tfrac{d}{ds}\|P_{\geq \beta}\mathbf{f}\|_h  \geq  \beta \|P_{\geq \beta}\mathbf{f}\|_h  -\| \mathbf{g} \|_{h},
 \end{align}
 hold for a.e. $s$ in $[s_1,s_2]$. 
\end{proposition}

\begin{proof}
We remind that $P_{\leq \beta}$ is the projection to the space spanned by eigenfunctions whose eigenvalues are less than or equal to $\beta$. Hence, we have the second inequality as follows.
\begin{equation}
    \tfrac{1}{2}\tfrac{d}{ds}\|P_{\leq \beta}\mathbf{f}\|_h^2=\langle P_{\leq \beta}\mathbf{f}, \partial_s \mathbf{f}\rangle_h=\langle P_{\leq \beta}\mathbf{f},\mathcal{L}\mathbf{f}+\mathbf{g}\rangle_h\leq \beta \|P_{\leq \beta}\mathbf{f}\|_h^2+\|P_{\leq \beta}\mathbf{f}\|_h \|\mathbf{g} \|_{h}.
\end{equation}
We can obtain the other inequalities in the same manner.
\end{proof}
 \bigskip

Let $h (\theta)$ and $S(l,\theta)$ be smooth solutions  to \eqref{eq-shrinker} and \eqref{Se} with $\eta \in \{0,1\}$, respectively. Then, the difference
\begin{equation}
 w(s,\theta):=S(e^s,\theta)-\sigma^{-1}e^{\sigma s}h(\theta) 
\end{equation}
satisfies 
\begin{equation}
\partial_s [w]=\mathcal{L}[w]-(0,E_1(w)+\eta E_2(w)),
\end{equation}
where $E_1,E_2$ are given in \eqref{eq-Ew} and \eqref{eq-Ew2}.
 \bigskip

\begin{proposition}\label{prop:Error12}
Let $h \in C^\infty(\mathbb{S}^1)$ be a $k$-fold symmetric solution to \eqref{eq-shrinker} for some  $\alpha \in (0,\frac{1}{4})$. Then, there exist some positive constants $\varepsilon,C$ only depending on $\alpha,k$ with the following significance. Given a solution $w$ to \eqref{eq-w'}, if $e^{-\sigma s}\|w\|_{C^4_{s,1}} \leq \varepsilon$ holds for $s \in [s_1,s_2]\subset [0,\infty)$ then  
\begin{align}
& \|(0,E_1(w)) \|_{h}  \leq  C e^{-\sigma s}  \|w  \|_{C^4_{s,1}}  \| [w] \|_h, &&|E_2(w)|\leq Ce^{-4\alpha s},
\end{align} 
hold for $s \in [s_1,s_2]$.
\end{proposition}
  
\begin{proof}
If  $e^{-\sigma s}\|w\|_{C^4_{s,1}}$ is sufficiently small, then we have
\begin{equation}
\big|  (h+ e^{- \sigma s} w_s )^{\frac1\alpha} -h^{\frac1\alpha} -\tfrac1{\alpha} h^{\frac1\alpha-1}e^{- \sigma s} w_s \big| \leq Ce^{-2\sigma \tau}w_s^2.
\end{equation}
This yields $|E_1(w)|\leq Ce^{-\sigma \tau}\|w\|_{C^4_{s,1}}w_s$, and thus  by Proposition \ref{prop-normequivalence} we obtain the first inequality. The second inequality is obvious by observing
 \begin{align}\label{eq:E2_expo_decay}
|E_2(w)|\leq C \big [  (1+ e^{-4\alpha s} (  h+e^{-\sigma s}w_s )^2)^{2-\frac{1}{2\alpha}}-1  \big ]\leq Ce^{-4\alpha s}.
 \end{align}
  \end{proof}

\begin{proposition}\label{prop:Error12_round}
Let $h\equiv (2\alpha \sigma)^\alpha$ for $\alpha \in (0,\frac{1}{4})$, which is the constant solution to \eqref{eq-shrinker}. Then, for each integer $m\geq 0$, there exist some positive constants $\varepsilon,C$ only depending on $\alpha,m$ with the following significance. Given a solution $w$ to \eqref{eq-w'}, if $e^{-\sigma s}\|w\|_{C^{m+4}_{s,1}} \leq \varepsilon$ holds for $s \geq s_0\geq 0$ then  
\begin{align}
& \|(0,\partial_\theta^m E_1(w)) \|_{h}  \leq  C e^{-\sigma s}  \|w  \|_{C^{m+4}_{s,1}}  \| [w] \|_{H^{m+1}\times H^m }, &&|\partial_\theta^m E_2(w)|\leq Ce^{-4\alpha s},
\end{align} 
hold for $s \geq s_0$.
\end{proposition}
  
\begin{proof}
Assuming that $e^{-\sigma s}\|w\|_{C^{m+4}_{s,1}}$ is sufficiently small, we have
\begin{equation}
\Big| \partial_\theta^m \big[  (h+ e^{- \sigma s} w_s )^{\frac1\alpha} -h^{\frac1\alpha} -\tfrac1{\alpha} h^{\frac1\alpha-1}e^{- \sigma s} w_s \big]\Big| \leq Ce^{-2\sigma \tau} \sum_{i+j\leq m} |\partial_\theta^i w_s||\partial_\theta^j w_s|.
\end{equation}
Hence,  Proposition \ref{prop-normequivalence} yields the first inequality. On the other hand, $e^{-\sigma s}\|w\|_{C^{m+4}_{s,1}} \leq \varepsilon$ implies
 \begin{align}
\big| \partial_\theta^i   (1+ e^{-4\alpha s}  (  h+e^{-\sigma s}w_s  )^2 )^{2-\frac{1}{2\alpha}}   \big |\leq Ce^{-4\alpha s},
 \end{align}
 for $i\geq 1$. Therefore, combining with \eqref{eq:E2_expo_decay} gives us the second inequality.
  \end{proof}

\bigskip

  \begin{proposition}\label{prop:Error_est_interval}
  Let $h \in C^\infty(\mathbb{S}^1)$ be a $k$-fold symmetric solution to \eqref{eq-shrinker} with $\alpha \in (0,\frac{1}{4})$. Then, there exist some positive constants $\varepsilon,C$ only depending on $\alpha,k$ with the following significance. Given a solution $w$ to \eqref{eq-w'}, if $e^{-\sigma s}\|w\|_{C^4_{s,1}} \leq \varepsilon$ holds for $s \in [s_1,s_2]\subset [0,\infty)$ then we have 
 \begin{align} 
\tfrac{d}{ds} e^{-\sigma s}\|P_{+}[w]\|_h & \geq  \lambda^+  e^{-\sigma s}\|P_{+}[w]\|_h  -C(e^{-2\sigma s}  \|   w  \|_{C^4_{s,1}}\|[w] \|_h +\eta e^{-(1+2\alpha) s }), \\
 \left|\tfrac{d}{ds}e^{-\sigma s}\|P_{\sigma}[w]\|_h\right| &\leq C(e^{-2\sigma s}  \|   w  \|_{C^4_{s,1}}\|[w] \|_h +\eta e^{-(1+2\alpha) s }),\\
  \tfrac{d}{ds}e^{-\sigma s}\|P_{-}[w]\|_h &\leq -\lambda^- e^{-\sigma s}\|P_{-}[w]\|_h + C(e^{-2\sigma s}  \|   w  \|_{C^4_{s,1}}\|[w] \|_h +\eta e^{-(1+2\alpha) s }), 
 \end{align}
almost every $s$ in $[s_1,s_2]$, where   $\sigma+\lambda^+$ is the least eigenvalue of $\mathcal{L}$ greater than $\sigma$, and  $\sigma-\lambda^-$ is the greatest eigenvalue of $\mathcal{L}$ less than $\sigma$. 
\end{proposition}

\begin{proof}
Since $4\alpha+\sigma=1+2\alpha$,  combining Proposition \ref{prop:ODE_general} and Proposition \ref{prop:Error12} yields the desired result.
\end{proof}

  \bigskip
  
  \begin{proposition}\label{prop:Error_est_round}
  Let $h\equiv (2\alpha \sigma)^\alpha$ for $\alpha \in (0,\frac{1}{4})$, which is the constant solution to \eqref{eq-shrinker}. Then, for each integer $m\geq 0$, there exist some positive constants $\varepsilon,C$ only depending on $\alpha,m$ with the following significance. Given a solution $w$ to \eqref{eq-w'}, if $e^{-\sigma s}\|w\|_{C^{m+4}_{s,1}} \leq \varepsilon$ holds in $s \in [ s_0,+\infty)$ for some $s_0\geq 0$, then we have 
 \begin{align} 
\tfrac{d}{ds} e^{-\sigma s}\|P_{+}[\partial_\theta^m w]\|_h & \geq  \lambda^+  e^{-\sigma s}\|P_{+}[\partial_\theta^m w]\|_h  -C(e^{-2\sigma s}  \|   w  \|_{C^{m+4}_{s,1}}\|[w] \|_{H^{m+1}\times H^m } +\eta e^{-(1+2\alpha) s }), \\
 \left|\tfrac{d}{ds}e^{-\sigma s}\|P_{\sigma}[\partial_\theta^m w]\|_h\right| &\leq C(e^{-2\sigma s}  \|    w  \|_{C^{m+4}_{s,1}}\|[w] \|_{H^{m+1}\times H^m } +\eta e^{-(1+2\alpha) s }),\\
  \tfrac{d}{ds}e^{-\sigma s}\|P_{-}[\partial_\theta^m w]\|_h &\leq -\lambda^- e^{-\sigma s}\|P_{-}[\partial_\theta^m w]\|_h + C(e^{-2\sigma s}  \|   w  \|_{C^{m+4}_{s,1}}\|[w] \|_{H^{m+1}\times H^m } +\eta e^{-(1+2\alpha) s }), 
 \end{align}
almost every $s$ in $[s_0,+\infty)$, where   $\sigma+\lambda^+$ is the least eigenvalue of $\mathcal{L}$ greater than $\sigma$, and  $\sigma-\lambda^-$ is the greatest eigenvalue of $\mathcal{L}$ less than $\sigma$. 
\end{proposition}

\begin{proof}
We  combine Proposition \ref{prop:ODE_general} and Proposition \ref{prop:Error12_round} to obtain the desired result.
\end{proof}
\bigskip
  
  \begin{theorem}\label{thm:MZ_ODE}
 Let $h\equiv (2\alpha \sigma)^\alpha$ for $\alpha \in (0,\frac{1}{4})$, which is the constant solution to \eqref{eq-shrinker}. Given a solution $w$ to \eqref{eq-w'}, if $ e^{-\sigma s}\|w\|_{C^{4}_{s,1}} \to 0$ as $s\to +\infty$, then there is some $s_0\geq 0$ with the following significance: the solution satisfies  $\| P_{\sigma}[  w]\|_{h}>0$ for $s\ge s_0$  and
 \begin{equation}\label{eq:MZ_neutral_dom}
  e^{-\sigma s}\|[  w]-P_{\sigma}[  w]\|_{H^{5}\times H^4 }+\left| \tfrac{d}{ds}e^{-\sigma s}\| P_{\sigma}[  w]\|_{h}\right|=o(1)\|e^{-\sigma s}P_{\sigma}[ w]\|_{h} \text{ as } s\to \infty, 
 \end{equation}
unless $\|[w]\|_h\leq Ce^{-\delta s} $ holds for  $s\ge s_0$ with some $C,\delta>0$.
  \end{theorem}

\begin{proof}
    
We observe that $\lim_{s\to +\infty} e^{-\sigma s}\|w\|_{C^{4}_{s,1}} = 0$ implies $\lim_{s\to +\infty} e^{-\sigma s}\|w\|_{C^{10}_{s,1}} = 0$ by \cite[Lemma B.1 (interior regularity)]{CCK_existence}. Thus, by Proposition \ref{prop:Error_est_round}, the absolutely continuous non-negative functions
\begin{align}
 & x=\sum_{j=0}^{4} e^{-\sigma s}\|P_{+}[\partial_\theta^j w]\|_h,&& y=\sum_{j=0}^{4} e^{-\sigma s}\|P_{\sigma}[\partial_\theta^j w]\|_h, && z= e^{- s}+\sum_{j=0}^{4} e^{-\sigma s} \|P_{-}[\partial_\theta^j w]\|_h
\end{align}
 satisfy
\begin{align}
 &x'\geq \lambda x -o(x+y+z), &&|y'|\leq o(x+y+z) , && z'\leq -\lambda z +o(x+y+z)
\end{align} 
for  $\lambda :=\min\{\lambda_+,\lambda_-,1\}$. Notice that we also have  $x+y+z\geq  e^{-s}>0$. Therefore, by Merle-Zaag ODE Lemma \ref{lem-MZODE}, either   $x+z=o(y)$ or $x+y=o(z)$ holds.

\medskip

We first consider the case $x+z=o(y)$. By the assumption $h_\theta=0$, we have $\langle \textbf{f}_\theta ,\textbf{g}\rangle_h=\langle \textbf{f} ,\textbf{g}_\theta\rangle_h$. Also, if $\mathcal{L}\vec{\varphi}=\beta \vec{\varphi}$, then $\mathcal{L}\vec{\varphi}_\theta=\beta \vec{\varphi}_\theta$. Thus, 
\begin{equation}
 y \leq \sum_{i=0}^4 \sum_{\beta^+_j=\sigma} e^{-\sigma s} |\langle [w],\partial_\theta^i \vec{\varphi}_j\rangle_h| \leq  C'e^{-\sigma s } \| P_{\sigma}[  w]\|_h   
\end{equation}
holds for some $C' =C'(\alpha)>1$ (cf. \cite[(3.23)]{ChoiSun}). This infers that
\begin{equation}
e^{-\sigma s}\|[w]-P_\sigma [w]\|_{H^{5}\times H^4}\leq o(1)e^{-\sigma s}\|P_\sigma [w]\|_h.
\end{equation}
In addition, Proposition \ref{prop:Error_est_interval} yields $ |\frac{d}{ds}e^{-\sigma s}\| P_{\sigma}[  w]\|_h|\leq o(1)e^{-\sigma s}\| P_{\sigma}[  w]\|_h+C  e^{-(1+2\alpha)s} \le  o(1)e^{-\sigma s}\| P_{\sigma}[  w]\|_h$, where  in the last inequality we used $C'e^{-\sigma s} \| P_{\sigma}[  w]\|_h \ge z\ge  e^{-s}>0 $ for sufficiently large $s$. This shows \eqref{eq:MZ_neutral_dom} and $\| P_{\sigma}[  w]\|_h>0$ for large $s\ge s_0$. 

In the second case  $x+y=o(z)$, Proposition \ref{prop:Error_est_interval} gives us $z'\leq -\lambda z+o(z)$ where $\lambda=\min\{\lambda_-,1\}$. Hence, we obtain $z\leq Ce^{-\delta s}$ for any $\delta <\lambda$. This shows $\| [w]\|_h \le C e^{-\delta s} $ for $s\ge s_0$. 

 \end{proof}
 
 \bigskip

 \subsection{Unique blow-down and rate of convergence I}\label{sec:norotation}

\begin{theorem}\label{thm-exponential1}
Let $S:[l_0,+\infty)\times \mathbb{S}^1\to \mathbb{R}$ be a smooth solution to \eqref{Se} in $\{l>l_0\}$ and assume $S(l,\theta)$ has the unique limit shape represented by a non-radial $k$-fold solution $h \in C^\infty(\mathbb{S}^1)$ to \eqref{eq-shrinker}. Then, there are some $\delta=\delta(\alpha, k)>0$ and  $\theta_\infty \in \mathbb{S}^1$ such that 
\begin{equation}
\lim_{\tau \to +\infty}e^{\delta\tau} \Vert  (1-2\alpha)S(e^s,\theta ) e^{-(1-2\alpha)s} -h(\theta-\theta_\infty)\Vert _{C^{4,\beta}_{\tau,1}}=0 \text{ where }\beta=2^{-1}.
\end{equation} 
\end{theorem}

\bigskip
The $C^{4,\beta}$-norm is introduced to use the Schauder's interior estimates in \cite[Lemma B.1 (interior regularity)]{CCK_existence}. Any $\beta\in(0,1)$ will be okay, but we fix $\beta=2^{-1}$ throughout this section so that the constants appearing in this section do not depend on $\beta$. Let us define 
\begin{align}
&h_{\theta'}(\theta):=h(\theta+\theta'), && S_{\theta'}(l,\theta):=  S(l,\theta+\theta'), && w_{\theta'}(s,\theta) :=   S_{\theta'}(e^s,\theta) -\sigma^{-1}e^{\sigma s}h(\theta),
\end{align}
and
\begin{equation}
\mathbf{h}(\theta) := (  \sigma ^{-1}  h'(\theta)  ,   h'(\theta)).
\end{equation}
Note that by Lemma 7.3 in  \cite{andrews2003classification},  $h'$ is an eigenfuction of the operator $L$ defined in \eqref{eq-linearopL} with the eigenvalue $2\sigma (\alpha-1)$. More importantly, the kernel of $L-2\sigma (\alpha-1)I$ is spanned by $h'$, where $I$ denotes the identity operator $If:=f$. Hence, $\mathbf{h}$ is the only eigenvector (upto scaling) with the eigenvalue $\sigma$ of the operator $\mathcal{L}$ defined in \eqref{eq-mathcalL}. Namely,
\begin{equation}\label{eq:int.Jacobi.field}
\ker (\mathcal{L}-\sigma \mathbf{I})=\text{span}\{   \mathbf{h}\},
\end{equation}
where $\mathbf{I}$ is the identity operator.
\bigskip

The theorem will follow from an iteration of following lemma. 
\begin{lemma}[Refinement by rotations] \label{lem-key}
	Let $h\in C^\infty(\mathbb{S}^1)$ be  a non-radial $k$-fold solution to \eqref{eq-shrinker} for some $3 \leq k\in \mathbb{N}$. Depending on $\alpha,k$, there exist constants $N, C <\infty  $ and $\varepsilon_0,\delta_0>0 $ with the following significance.  Given a solution $S(l,\theta)\in C^{\infty}([l_0,\infty)\times \mathbb{S}^1) $ to \eqref{Se},  suppose there exist a function $\theta_0 :[\tau_0,+\infty) \to \mathbb{S}^1$ with $\tau_0\geq \max\{ N, \ln  l_0+1\}$ and $\eps \in (0,\eps_0]$ that satisfy
\begin{equation}\label{eq:assumption_fitting}
e^{-\sigma \tau}\| w_{\theta_{0}(\tau)}  \|_{C^{4,\beta}_{\tau,1}}\leq   \varepsilon + e^{-\delta_0 \tau}  
\end{equation}
for $\tau \geq \tau_0$. Then, there exists a refinement $\theta_1 :[\tau_0+N,+\infty) \to \mathbb{S}^1$ which satisfies 
\begin{equation} \label{eq-4280}
 e^{-\sigma \tau}\| w_{\theta_{1}(\tau)}  \|_{C^{4,\beta}_{\tau,1}} \le   \tfrac{1}{2}  \varepsilon + e^{-\delta_0 \tau}  ,
\end{equation} 
for $\tau \geq \tau_0+N$, and
\be \label{eq-4290}\Vert h_{-\theta_0(\tau' )}-h_{-\theta_1(\tau)}\Vert _{C^{4,\beta}} \le C(\eps+ e^{-\delta_0\tau}),\ee
for $\tau \geq \tau_0+N$ and $\tau'\in [\tau-N,\tau]$. 
\end{lemma}
Before giving the proof of Lemma \ref{lem-key}, let us first prove Theorem \ref{thm-exponential1}. 
\begin{proof}[Proof of Theorem \ref{thm-exponential1}]
By the unique limit shape assumption (Definition \ref{def-41}), there exist large $\tau_0$ and a function $\theta_0:[\tau_0,\infty)\to \mathbb{S}^1$ that satisfy the assumption of Lemma \ref{lem-key} with  the choice $\eps=\eps_0$. By iterating Lemma \ref{lem-key}, we obtain a sequence of functions $\theta_n(\tau):[\tau_0+nN,\infty ) \mapsto \mathbb{S}^1$ such that, for $\tau \ge \tau_0+nN$ and $\tau'\in[
\tau-N,\tau]$,
\be \Vert h_{-\theta_{n-1}(\tau')}-h_{-\theta_n(\tau)}\Vert _{C^{4,\beta}}  \le C \left(\frac{\eps_0}{2^{n-1}} +e^{-\delta_0 \tau }\right)  \ee  and 
\be \label{eq-4056}e^{-\sigma \tau} \Vert w_{\theta _n(\tau)}\Vert _{C^{4,\beta}_{\tau,1}} \le \frac{\eps_0 }{2^n} +e^{-\delta_0\tau }.\ee 
Since
\be \Vert  h_{-\theta_{n-1}(\tau_0 +(n-1)N)}-h_{-\theta_k(\tau_0+nN)}\Vert _{C^{4,\beta}}  \le C \left(\frac{\eps_0}{2^{n-1}} +e^{-\delta_0 (\tau_0+nN) }\right)\ee  and the right hand side converges geometrically as $n\to \infty$, there is $\theta_{\infty}\in \mathbb{S}^1$ such that
\be \Vert h_{-\theta_\infty} - h_{-\theta_n(\tau )} \Vert_{C^{4,\beta}} \le C' e^{-\delta_1 nN},\ee for $\tau \in [\tau_0+nN , \tau_0+(n+1)N]$. Here $\delta_1= \min (\delta_0, \frac{\ln 2}{N})>0$ and $C'<\infty$ may depend on $\alpha, k, \tau_0, \eps_0,\delta_0$. Together with \eqref{eq-4056}, this implies that for $\tau \in[\tau_0+nN,\tau_0+(n+1)N]$
\be e^{-\sigma \tau} \Vert w_{\theta_{\infty}} \Vert _{C^{4,\beta}_{\tau,1}}\le C'' e^{-\delta_1 nN},\ee  for some possibly larger $C''<\infty$. This proves the theorem.
\end{proof}

\bigskip

\begin{proof} [Proof of Lemma \ref{lem-key}]
\textit{Step 1: $C^{4,\beta}$-norm in large regions.} By the assumption \eqref{eq:assumption_fitting}, there is some constant $C$ only depending on $\sigma=1-2\alpha$  such that $\|e^{-\sigma s}S(e^s,\cdot )-\sigma^{-1} h_{-\theta_{0}(\tau)}   \|_{C^{4,\beta}_{\tau,1}}\leq C (\varepsilon+e^{-\delta_0\tau})$. Hence,  
\begin{align}
\|h_{-\theta_0(\tau)}-h_{-\theta_0(\tau+\frac{1}{2} )}\|_{C^{4,\beta}} & \leq \|\sigma e^{-\sigma s}S(e^s,\cdot )-  h_{-\theta_{0}(\tau)}   \|_{C^{4,\beta}_{\tau,\frac{1}{2} }}+\|\sigma e^{-\sigma s}S(e^s,  )-  h_{-\theta_{0}(\tau+\frac{1}{2} )}   \|_{C^{4,\beta}_{\tau,\frac{1}{2} }}\\
&\leq   \|\sigma e^{-\sigma s}S(e^s,\cdot )-  h_{-\theta_{0}(\tau)}   \|_{C^{4,\beta}_{\tau,1}}+\|\sigma e^{-\sigma s}S(e^s,\cdot )-  h_{-\theta_{0}(\tau+\frac{1}{2} )}   \|_{C^{4,\beta}_{\tau+ \frac{1}{2},1 }}
\end{align}
yields
$\|h_{-\theta_0(\tau)}-h_{-\theta_0(\tau+\frac{1}{2} )}\|_{C^{4,\beta}}\leq C(\varepsilon +e^{ -\delta_0 \tau})$, and thus \be\label{eq-4330} \|h_{-\theta_0(\tau_1)}-h_{-\theta_0(\tau_2)}\|_{C^{4,\beta}}\leq C\varepsilon |\tau_1-\tau_2|+Ce^{-\delta_0\tau_1}\ee  holds whenever $\tau_0 \leq \tau_1 \leq \tau_2 -1$. Therefore, given $N \in \mathbb{N}$, 
\begin{align}
&\| e^{-\sigma s}w_{\theta_{0}(\tau)}(s,\cdot)  \|_{C^{4,\beta}_{\tau' ,1}} =  \| e^{-\sigma s}S(e^s,\cdot )- \tfrac{1}{\sigma}  h_{-\theta_{0}(\tau)}  \|_{C^{4,\beta}_{\tau',1}}\\
& \leq \|  e^{-\sigma s}S(e^s,\cdot )- \tfrac{1}{\sigma} h_{-\theta_{0}(\tau')}  \|_{C^{4,\beta}_{\tau',1}} +\tfrac{1}{\sigma}\|h_{-\theta_{0}(\tau')}-  h_{-\theta_{0}(\tau)}  \|_{C^{4,\beta}}  \leq C(N+1)\varepsilon+Ce^{-\delta_0(\tau-N)}
\end{align}
holds for $\tau' \in [\tau-N,\tau+N]$ and $\tau \geq \tau_0+N$. Namely, there is some $C$ only depending on $\alpha$ such that
\begin{equation}\label{eq:w0_rough.esti}
e^{-\sigma \tau' } \| w_{\theta_0(\tau)} \|_{C^{4,\beta}_{\tau',1}} \leq CN\varepsilon +C  e^{ -\delta_0 (\tau-N)}
\end{equation}
holds for $\tau' \in [\tau-N,\tau+N]$, $\tau \geq \tau_0+N$, and $N \in\mathbb{N}$.

\bigskip

\textit{Step 2: dynamics of projections.} We recall the eigenvalues $\lambda_\pm$ in Proposition \ref{prop:Error_est_interval}, and define $\delta_0,\varepsilon_0$ by
\begin{align}
&\delta_0=\tfrac{1}{10} \min\{\lambda_+,\lambda_-,\alpha \}, && \varepsilon_0=e^{-10N}.
\end{align}
Notice that $\delta_0$ only depends on $\alpha,k$. We  also define
\begin{align}
&x(s)=e^{-\sigma s}\|P_{+}[w_{\theta_0(\tau)}]\|_h, &&y(s)=e^{-\sigma s}\|P_{\sigma}[w_{\theta_0(\tau)}]\|_h, && z(s)=e^{-\sigma s}\|P_{-}[w_{\theta_0(\tau)}]\|_h+\eta e^{-s}.
\end{align}
Since $ \delta_0<1$, we have $e^{- s} \leq e^{- (\tau-N)}\leq e^{-\delta_0(\tau-N)}$ for $s\geq \tau-N$. Hence, \eqref{eq:w0_rough.esti} yields
\begin{equation}\label{eq:sum.proj_esti}
x+y+z\leq CN\varepsilon +  C e^{  -  \delta_0(\tau-N) },
\end{equation}
for $s \in [\tau-N,\tau+N]$.

\bigskip
Now, we choose large enough $N$ and then apply Proposition \ref{prop:Error_est_interval} to get $C_1=C_1(\alpha, k)>0$ satisfying
\begin{align}
 z'\leq  - (\min\{\lambda_-,1\})  z+ C_1  (N\varepsilon +e^{-\delta_0(\tau-  N)}+e^{-2 \alpha s})(x+y+z) ,
\end{align}
for $s\in [\tau- N,\tau+ N]$. We may assume 
\begin{equation}
10C_1N e^{-\delta_0 N} \leq \delta_0. 
\end{equation}
Then, remembering $s \geq \tau-  N \geq \tau_0\geq 10N$ and $\alpha \geq 10\delta_0 $, we have
\begin{equation}
z'\leq   - 9\delta_0  z+ \delta_0 e^{-9\delta_0 N }(x+y).
\end{equation}
In the same manner,  we can obtain
\begin{align}
& |y'| \leq  \delta_0 e^{- 9\delta_0 N} (x+y+z), && x' \geq  9\delta_0 x  -  \delta_0e^{-9\delta_0 N}(y+z). 
\end{align}
Thus, applying \cite[Lemma B.2]{CS2}, a quantitative variant of Merle-Zaag lemma, with \eqref{eq:sum.proj_esti} yields
\begin{equation}
x+z\leq C e^{-9\delta_0 N}(N\varepsilon+e^{-\delta_0(\tau- N) })+C(N\varepsilon+e^{-\delta_0(\tau- N) })e^{-\frac{9}{4} \delta_0 N}
\end{equation}
for $s\in [\tau- N/2,\tau+N/2]$. We may assume $N \geq 20$ so that we have 
\begin{equation}\label{eq:decay_interval}
 e^{-\sigma s}\|(P_{-}+P_+)[w_{\theta_0(\tau)}]\|_h \leq Ce^{-\delta_0 N}(\varepsilon   + e^{-\delta_0 \tau})    
\end{equation}
for $s\in [\tau- 10,\tau+10]$.

\bigskip

\textit{Step 3: rotation eliminating the neutral projection.}
Remembering $\mathbf{h}  := (  \sigma ^{-1}  h'   ,   h')$, we define
\begin{equation}\label{eq_thetas} 
\theta_1(\tau)= \theta_0(\tau)-e^{-\sigma \tau} \langle  [  w_{\theta_0(\tau)} ](\tau ,\cdot),  \mathbf{h}  \rangle_{h}\|\mathbf{h}  \|_h^{-2}.  
\end{equation}Then, \eqref{eq:int.Jacobi.field} implies
\begin{equation}\label{eq:angle.def}
P_\sigma [w_{\theta_0(\tau)}](\tau ,\cdot)=\left \langle [w_{\theta_0(\tau)}] (\tau ,\cdot), \tfrac{\mathbf{h}}{\|\mathbf{h}\|_h} \right \rangle_h \tfrac{\mathbf{h}}{\|\mathbf{h}\|_h}=-  (\theta_1-\theta_0)e^{\sigma \tau }\mathbf{h}  .
\end{equation}
Since  $\|\mathbf{h}\|_h \geq c$ holds for some $c>0$, \eqref{eq:w0_rough.esti}  yields
\begin{equation}\label{eq:angle_diff}
|\theta_1-\theta_0| \leq Ce^{-\sigma  \tau }\|[ w_{\theta_0(\tau)}](\tau ,\cdot)\|_h\leq  CN \varepsilon+Ce^{-\delta_0(\tau-N)}.
\end{equation}
Therefore, by using \eqref{eq:w0_rough.esti}, \eqref{eq:angle_diff}, and
\begin{align}
\|w_{\theta_1}\|_{C^{4,\beta}_{\tau',1}}=&\,\|S_{\theta_1}(e^s,\cdot)-\sigma^{-1}e^{\sigma s}h\|_{C^{4,\beta}_{\tau',1}}=\|S_{\theta_0}(e^s,\cdot)-\sigma^{-1}e^{\sigma s}h_{(\theta_0-\theta_1)}\|_{C^{4,\beta}_{\tau',1}}\\
\leq  &\,\|w_{\theta_0}\|_{C^{4,\beta}_{\tau',1}}+\tfrac{1}{\sigma} \| e^{\sigma s}h- e^{\sigma s}h_{(\theta_0-\theta_1)}\|_{C^{4,\beta}_{\tau',1}}\leq \|w_{\theta_0}\|_{C^{4,\beta}_{\tau',1}}+Ce^{\sigma\tau'}|\theta_1-\theta_0|,
\end{align}
we can obtain
\begin{equation}\label{eq:w1_rough.esti}
e^{-\sigma \tau' } \| w_{\theta_1(\tau)} \|_{C^{4,\beta}_{\tau',1}} \leq CN\varepsilon +C  e^{ -\delta_0 (\tau- N)}
\end{equation}
for $\tau' \in [\tau-N,\tau+N]$, $\tau \geq \tau_0+ N$, and $N \in\mathbb{N}$.

\bigskip

On the other hand, 
\begin{equation}
 (S_{\theta_1}(e^s, \theta')-S_{\theta_0}( e^s,\theta'))-\sigma^{-1}e^{\sigma s} (h_{\theta_1-\theta_0}(\theta')-h(\theta'))=\int_{\theta'}^{\theta'+\theta_1-\theta_0} \tfrac{\partial }{\partial \theta}w_{\theta_0}(s,\theta) d\theta
\end{equation}
implies
\begin{equation}
\|(S_{\theta_1}(e^s,\cdot  )-S_{\theta_0}( e^s, \cdot))-\sigma^{-1}e^{\sigma s} (h_{\theta_1-\theta_0}-h)\|_{C^1_{\tau',1}}\leq |\theta_1-\theta_0| \|w_{\theta_0} \|_{C^2_{\tau',1}}.
\end{equation}
In addition, the Taylor expansion of $h$ yields
\begin{equation}
\| \sigma^{-1}e^{\sigma s} (h_{\theta_1-\theta_0}-h)-(\theta_1-\theta_0)\sigma^{-1}e^{\sigma s}  h'\|_{C^1_{\tau',1}}\leq C |\theta_1-\theta_0|^2 e^{\sigma \tau' }.
\end{equation}
Hence, combining the above inequalities with \eqref{eq:w0_rough.esti} and \eqref{eq:angle_diff} leads us to
 \begin{equation}
    \|(S_{\theta_1}-S_{\theta_0})-(\theta_1-\theta_0)\sigma^{-1}e^{\sigma s} h' \|_{C^1_{\tau',1}}\leq C(N^2\varepsilon^2 +e^{-2\delta_0(\tau- N)} ) e^{\sigma \tau'}.
 \end{equation}
Therefore, recalling $w_{\theta_1}-w_{\theta_0}=S_{\theta_1}-S_{\theta_0}$, we have
 \begin{equation}\label{eq:rot.h.esti}
   e^{-\sigma s} \|[w_{\theta_1}]-[w_{\theta_0}]-(\theta_1-\theta_0)e^{\sigma  s}\mathbf{h}  \|_{h}\leq CN^2\varepsilon^2+Ce^{-2\delta_0(\tau- N)} 
 \end{equation}
 for $s\in [\tau- N,\tau+ N]$.  Notice that $2(\tau- N) \geq \tau_0+\tau- N \geq \tau+9N$ implies $e^{-2\delta_0(\tau- N)}\leq e^{-\delta_0  (\tau +N)}$, and we also  have $N^2\varepsilon\leq N^2e^{-10N}\leq e^{-N}$ for large $N$. Since \eqref{eq:int.Jacobi.field} implies $P_\sigma (e^{\sigma  s}\mathbf{h}) =e^{\sigma  s}\mathbf{h} $,   \eqref{eq:angle.def} and \eqref{eq:rot.h.esti} give us
  \begin{equation}\label{eq:neutral.center.esti}
e^{-\sigma  \tau}    \|P_\sigma  [w_{\theta_1}](\tau ,\cdot )\|_{h}\leq C e^{-\delta_0 N}(\varepsilon +   e^{-\delta_0  \tau} ).
 \end{equation}
 In addition, \eqref{eq:int.Jacobi.field} implies  $P_\pm \mathbf{h}=0$, and therefore \eqref{eq:rot.h.esti} yields
  \begin{equation}
   e^{-\sigma s} \|(P_-+P_+)( [w_{\theta_1}]-[w_{\theta_0}] ) \|_{h}\leq CN^2\varepsilon^2+Ce^{-2\delta_0(\tau- N)} .
 \end{equation}
 Hence, by using \eqref{eq:decay_interval} we can obtain 
\begin{equation}\label{eq:non.neutral.esti}
 e^{-\sigma s}\|(P_{-}+P_+)[w_{\theta_1}]\|_h  \leq   C_2   e^{-\delta_0 N}(\varepsilon +   e^{-\delta_0  \tau} ),    
\end{equation}  
for $s\in [\tau-10,\tau+10]$ and some constant $C_2 \geq 1$ only depending on $\alpha,k$.
 
 \bigskip
 
 \textit{Step 4: decay improvement.} 
We define $\hat y(s):=e^{-\sigma s}    \|P_\sigma  [w_{\theta_1}](s,\cdot )\|_{h}+ C_2   e^{-\delta_0 N}(\varepsilon +   e^{-\delta_0  \tau} )$. Then, \eqref{eq:w1_rough.esti} implies that we can apply Proposition \ref{prop:Error_est_interval} by choosing sufficiently large $N$. Hence, we have
\begin{equation}
|\hat y'| \leq C (N\varepsilon+e^{-\delta_0(\tau- N)})  e^{-\sigma s}\|[w_{\theta_1}]\|_h+C\eta e^{-(1+2\alpha)s}.
\end{equation}
for $s\in [\tau- N,\tau+ N]$. Notice that we have $e^{-(1+2\alpha)s}\leq e^{-\delta_0(\tau- N)}e^{-\delta_0(\tau+ N)}$. Thus, \eqref{eq:non.neutral.esti} implies
\begin{equation}
|\hat y'| \leq C (N\varepsilon+e^{-\delta_0(\tau- N)})\hat y'\leq C  \hat y'
\end{equation}
for $s\in [\tau-10,\tau+10]$. Therefore, for  $s\in [\tau-10,\tau+10]$ we have
\begin{equation}
\log \hat y(s)\leq \log \hat y(\tau )+10 \sup_{|s-\tau|\leq 10}\left|\tfrac{d}{ds}\log \hat y\right|\leq C+\log \hat y(\tau),
\end{equation}
and thus \eqref{eq:neutral.center.esti} implies $ e^{-\sigma s}\|P_\sigma [w_{\theta_1}]\|_h  \leq   C e^{-\delta_0 N}( \varepsilon +  e^{-\delta_0  \tau })$. So, adding to \eqref{eq:non.neutral.esti} yields
\begin{equation}
 e^{-\sigma s}\|[w_{\theta_1}]\|_h  \leq  Ce^{-\delta_0 N}( \varepsilon  +   e^{-\delta_0 \tau}),    
\end{equation}  
for $s\in [\tau-10,\tau+10]$ and $\tau \geq \tau_0+N$. Then, the interior regularity  in  \cite[Lemma B.1 (interior regularity)]{CCK_existence}, 
\begin{equation}
 e^{-\sigma \tau}\|w_{\theta_1(\tau)}\|_{C^{4,\beta}_{\tau,1}}  \leq  Ce^{-\delta_0 N}( \varepsilon  +   e^{-\delta_0 \tau})    
\end{equation}  
for $\tau \geq \tau_0+N$. Hence, we can choose sufficiently large $N$ to conclude \eqref{eq-4280}. Next, \eqref{eq-4290} is a direct consequence of \eqref{eq:angle_diff} and \eqref{eq-4330}.  
\end{proof}

\bigskip

\subsection{Unique blow-down and rate of convergence II}\label{sec-exponential-round}

In this subsection, we will prove the following fast convergence theorem for asymptotically round solutions.

\begin{theorem}\label{thm-exponential2}
Let $S:[l_0,+\infty)\times \mathbb{S}^1\to \mathbb{R}$ be a smooth solution to \eqref{Se} (on $\{l>l_0\}$) and assume $S(l,\theta)$  has the unique limit shape represented by the constant solution $h\equiv [2\alpha (1-2\alpha)]^\alpha$ to \eqref{eq-shrinker}. Then, there is some $ \delta>0$ satisfying 
\begin{equation}
\lim_{\tau \to +\infty}e^{\delta\tau} \|  (1-2\alpha)S(e^s,\theta ) e^{-(1-2\alpha)s} -h \|_{C^{10}_{\tau,1}}=0.
\end{equation} 
\end{theorem}
 
 To prove this theorem, we first explicitly compute the basis in \eqref{eq-basis} as follows:
\begin{equation}\label{eq-roundbasis} (1, \beta^\pm_0),\, (\cos( j\theta), \beta^\pm_j\cos (j\theta)), (\sin( j\theta), \beta^\pm_j\sin (j\theta)),
\end{equation}
where $j \in \mathbb{N}$ and  $\beta^\pm_j$ are solutions to $\beta^2+\beta- 2\alpha(1-2\alpha)( j^2-1)=0$, namely
\begin{equation}\label{eq:betapm}
\beta^\pm_j = -\tfrac12  \pm \tfrac12 \sqrt {1+8\alpha(1-2\alpha)(j^2-1)}.
\end{equation}
We recall the function $w(s,\theta)=S(e^s,\theta)-\sigma^{-1}e^{\sigma s}h$ with $h\equiv (2\alpha \sigma)^\alpha$, and define   
\begin{equation}\label{eq:coefficients} 
   \sigma^2 (he^{ \sigma s})^{-1} [w](\cdot,s)  = \sum_{j=0 }^\infty\sum_{*\in \{\pm\}} A^\ast_j(s)(\cos(j\theta), \beta^\ast_j\cos(j\theta))+B^\ast_j(s)(\sin(j\theta), \beta^\ast_j\sin(j\theta)).
  \end{equation}
Here, we ignore $B^\pm_0$ or simply define $B^{\pm}_0(s) \equiv 0$.  
\begin{remark}
Notice that if $\alpha=k^{-2}$ for some $ 3\leq k\in \mathbb{N}$ then we have
\begin{align}\label{eq:beta_k}
& \beta_k^+=\sigma, && |A^+_k|^2+|B^+_k|^2= c_0 e^{-2\sigma s} \|P_\sigma [w]\|_h^2,
\end{align}
where $c_0=\sigma^4 h^{-2} \|(\cos k\theta ,\sigma \cos  k\theta )\|_h^{-2}$. Here,
\begin{equation}
P_\sigma(v_1,v_2)=\mathrm{proj}_{\mathrm{span}\{ (\cos k\theta  ,\sigma  \cos k\theta ),(\sin k\theta  ,\sigma  \sin  k\theta )\}} (v_1,v_2).
\end{equation}
\end{remark}

\bigskip

Now, we provide the goal of this subsection. 
\begin{lemma}\label{lem:round_exp_goal}
Suppose that $\alpha=k^{-2}$ for some $3\leq k \in \mathbb{N}$, and there is some $s_0\gg 1$ such that the function 
\begin{equation}
\rho(s):=|A^+_k(s)|^2+|B^+_k(s)|^2
\end{equation}
satisfies 
\begin{align}\label{eq:rho_condition}
&\rho(s)>0 , &&\lim_{s\to +\infty}  \rho(s) =0, && \lim_{s\to +\infty} e^{\delta s}\rho(s) =+\infty, && e^{-2\sigma s}\|[w]-P_\sigma [w]\|_{H^{5}\times H^4}^2+\left|\tfrac{d}{ds} \rho \right|=o(\rho),
\end{align}
for all $s\geq s_0$ and for any $\delta>0$. Then, $\displaystyle \lim_{s\to +\infty}e^{-\sigma s}\|[w]\|_h=0$ fails.
\end{lemma}

\bigskip

This lemma proves Theorem \ref{thm-exponential2} by contradiction.
\begin{proof}[Proof of Theorem \ref{thm-exponential2}]
Since the limit shape is represented by the constant $h=[2\alpha \sigma]^\alpha $, the convergence assumption in Definition \ref{def-41} holds along arbitrary sequences $\lambda_i\to \infty$ and thus we have $\lim_{s\to \infty} e^{-\sigma s}\|[w]\|_h= 0$. Then, the regularity \cite[Lemma B.1 (interior regularity)]{CCK_existence} implies 
\begin{equation}\label{eq:C_infty_convergence}
\lim_{s\to +\infty} e^{-\sigma s}\|w\|_{C^m_{s,1}}=0,
\end{equation}
for each $m\geq 0$. Towards a contradiction, we assume that 
\begin{equation}\label{eq:assume_exp_decay}
\lim_{s\to+\infty} e^{\delta s} e^{-\sigma s}\|[w]\|_h =+\infty
\end{equation}
holds for every $\delta>0$. Then, by Theorem \ref{thm:MZ_ODE}, the neutral projection $P_\sigma[w]$ is non-zero unless $e^{-\sigma s}\|[w]\|_h$ decays exponentially fast. Therefore, observing $\beta_j^-\leq -\frac{1}{2}<\sigma$ in 
\eqref{eq:betapm}, there must be an integer $j$ such that $\beta^+_j=\sigma$. From this, we infer that $\alpha=k^{-2}<\frac14 $ for some $3\leq k\in \mathbb{N}$. Combining Theorem \ref{thm:MZ_ODE} with the assumption \eqref{eq:assume_exp_decay}, we obtain \eqref{eq:rho_condition}.  Hence, Lemma \ref{lem:round_exp_goal} completes the proof by contradiction.
\end{proof}

 \bigskip
 
To show Lemma \ref{lem:round_exp_goal}, we assume \eqref{eq:rho_condition} and  \eqref{eq:C_infty_convergence} in what follows. We define the functions 
\begin{align}
Q^\pm:= (|A^+_k|^2 -|B^+_k|^2) A^\pm_{2k}+  2 A^+_k B^+_k B^\pm_{2k}, 
\end{align} 
and 
\begin{align}
&  Q:= Q^++Q^-, &&A_j := A^+_j+A^-_j, && B_j := B^+_j+B^-_j,
\end{align}
for each $j\geq 0$. Also, we define the projections $\mathfrak{p}_j^c,\mathfrak{p}_j^s$ by
\begin{align}
\mathfrak{p}_0^c f= (2\pi)^{-1}\langle  f,  1 \rangle_{L^2(\mathbb{S}^1)}, && \mathfrak{p}_j^c f= \pi^{-1} \langle  f,  \cos j \theta \rangle_{L^2(\mathbb{S}^1)}, &&\mathfrak{p}_j^s f= \pi^{-1}\langle  f,  \sin j \theta \rangle_{L^2(\mathbb{S}^1)},
\end{align}
where $f\in L^2(\mathbb{S}^1)$ and $j \in \mathbb{N}$. Namely, assuming $\mathfrak{p}_0^s\equiv 0$, we have 
\begin{align}
f= \sum_{j=0}^\infty (\mathfrak{p}_j^c f) \cos j\theta + (\mathfrak{p}_j^s f) \sin j \theta.
\end{align}
Moreover, given a set $\mathcal{I}$ of indices $(j,*)\in \mathbb{Z}_{\geq 0}\times \{\pm\}$, where $\mathbb{Z}_{\geq 0}=\mathbb{N}\cup \{0\}$, we define the projection   ${P}^\perp_{\mathcal{I}}$ by
\begin{align}
   P^\perp_{\mathcal{I}}\mathbf{f}  =\mathbf{f}-\sum_{(j,*)\in \mathcal{I}} P_{=\beta_{j}^*}\mathbf{f}
\end{align}
with $\mathbf{f} \in C^\infty(\mathbb{S}^1)$.

\bigskip
   
\begin{lemma}[Error expansion upto third order] \label{lem_errorest} For each $m\geq 0$, the following holds
\begin{equation}
\left\|\sigma^2(  h e^{\sigma s} )^{-1} E(w)  + \mathcal{Q} (w;s,\theta)+\mathcal{C} (w;s,\theta)\right\|_{C^m(\mathbb{S}^1)} \lesssim   \left[e^{-\sigma s}\|w\|_{C^{m+2}_{s,1}}\right]^4+\eta e^{-4\alpha s},
\end{equation}
where
\begin{align}
\mathcal{Q}(w;s,\theta)&:= \sigma^3(h e^{\sigma s})^{-2} \left[2 w_s(  w +  w_{\theta\theta})+\tfrac{1-\alpha }{ \alpha\sigma }   w_s^2\right],\\ \label{eqCubic}
\mathcal{C}(w;s,\theta)&:=\tfrac{ 1-\alpha  }{3\alpha^2} \sigma^3(h e^{\sigma s})^{-3}  w_s^2 \left[   w_s + 3\alpha(  w + w_{\theta\theta})  \right].
\end{align} 
\end{lemma}

\begin{proof}
We denote $x=(h e^{\sigma s})^{-1} w_s$ and $y=(h e^{\sigma s})^{-1} (w_{\theta\theta}+w)$. Since $h^{\frac{1}{\alpha}}=2\alpha\sigma$,  
\begin{equation}
(  h e^{\sigma s} )^{-1} E_1=2\sigma  x y +2\alpha \sigma  \left[(1+x)^{\frac{1}{\alpha}}-1-\tfrac{1}{\alpha}x\right]( \sigma^{-1}+ y),
\end{equation}
by Lemma \ref{lem-w-eq}. Thus, Taylor expansion gives us 
\begin{align}
(1+x)^{\frac{1}{\alpha}}-1-\tfrac{1}{\alpha}x= \tfrac{1}{2\alpha}\left(\tfrac{1}{\alpha} -1\right) x+\tfrac{1}{6\alpha}\left(\tfrac{1}{\alpha} -1\right)\left(\tfrac{1}{\alpha} -2\right)x^3+O(x^4) 
= \tfrac{1-\alpha }{2\alpha^2}  x^2+\tfrac{(1-\alpha)\sigma}{6\alpha^3} x^3+O(x^4).
\end{align}
Hence,
\begin{equation}
( \sigma h e^{\sigma s} )^{-1}  E_1=2  x y +\tfrac{1-\alpha }{ \alpha\sigma }  x^2+\tfrac{ 1-\alpha }{3\alpha^2} x^3+\tfrac{ 1-\alpha}{ \alpha }  x^2y+O(|x^4|+|x^3y|).
\end{equation}
Also, Lemma \ref{lem-w-eq} says
\begin{equation}
(  h e^{\sigma s} )^{-1}  E_2= 2\alpha \sigma (1+ x)^{\frac{1}{\alpha}}(\sigma^{-1}+y)\left[\left(1+h^2e^{-4\alpha s}(1+x)^2\right)^{2-\frac{1}{2\alpha}}-1\right]= O(e^{-4\alpha s}).
\end{equation}
Therefore, we can obtain the desired result for $m=0$. For $m\geq 1$, we consider the derivatives of the Taylor expansions of $E_1$ and $E_2$, and obtain  $\partial_\theta^m E_1= O(  (e^{-\sigma s}\|w\|_{C^{m+2}_{s,1}})^4)$ and $\partial_\theta^m E_2= O(e^{-4\alpha s})$  so  the desired result follows. 
\end{proof}

\begin{proposition}\label{prop:err_rough_rho}
For each $m\geq 0$, the following hold
\begin{align}\label{eq-4100}
&e^{-\sigma s}\|w\|_{C^m_{s,1}} =O(\rho^{\frac{1}{2}}), && \left\|\sigma^2(  h e^{\sigma s} )^{-1} E(w)  + \mathcal{Q} (w;s,\theta)+\mathcal{C} (w;s,\theta)\right\|_{C^m(\mathbb{S}^1)}=O(\rho^2).
\end{align}
\end{proposition}

\begin{proof}
By the condition \eqref{eq:rho_condition}, we have $\frac{d}{ds}\log \rho =o(1)$. Hence, there is some constant $C$ only depending on $\alpha$ such that given $L\geq 1$
\begin{equation}
\sup_{|s-\tau|\leq L}e^{-2\sigma s}\|[w(s,\cdot)]\|_h^2\leq C \rho(\tau)
\end{equation}
holds for sufficiently large $\tau$. Thus, by the regularity \cite[Lemma B.1 (interior regularity)]{CCK_existence} $e^{-\sigma s}\|w\|_{C^m_{s,1}} =O(\rho^{\frac{1}{2}})$ for each $m\geq 0$. Therefore, combining the condition \eqref{eq:rho_condition} and Lemma \ref{lem_errorest} completes the proof. 
\end{proof}
 
\bigskip

\begin{lemma}[Expansion of quadratic error]\label{lem:Q.error.expansion} For each $m\geq 0$ and $s\gg 1$, the series
\begin{align}\label{eq:Q.error.expansion}
 \partial_\theta^m\mathcal{Q}= \lim_{k\to \infty}  \sum_{i,j\leq k}\sum_{*,*'}\sigma^{-1}  \beta^*_i   \left( 1-j^2  +\tfrac{  1-\alpha }{2\sigma \alpha  }\beta^{*'}_j   \right) \partial_\theta^m F_{(i,*),(j,*')}
\end{align}
absolutely converges in the $L^\infty(\mathbb{S}^1)$-norm, where
\begin{equation} \label{eq-Fvec}
F_{(i,*),(j,*')}(s,\theta) :=\bigg(\begin{array}{c}
A^*_iA^{*'}_j -B^*_i B^{*'}_j\\
B^*_iA^{*'}_j +A^*_i B^{*'}_j\end{array}\bigg ) \cdot
\bigg (\begin{array}{c}
\cos(i+j)\theta\\
\sin(i+j)\theta\end{array}\bigg)+
\bigg (\begin{array}{c}
A^*_iA^{*'}_j +B^*_i B^{*'}_j\\
B^*_iA^{*'}_j -A^*_i B^{*'}_j\end{array}\bigg ) \cdot
\bigg (\begin{array}{c}
\cos(i-j)\theta\\
\sin(i-j)\theta\end{array}\bigg ).
\end{equation}
In addition, given an index set $\mathcal{I}\subset \mathbb{Z}_{\geq 0}\times \{\pm\}$, there are some $s_0,C$ such that
\begin{equation}\label{eq:Q.error.estimate}
  \bigg|\partial_\theta^m \mathcal{Q}  -\sum_{(i,*),(j,*')\in \mathcal{I}} \sigma^{-1}\beta^*_i   \left( 1-j^2  +\tfrac{  1-\alpha }{2\sigma \alpha  }\beta^{*'}_j   \right) \partial_\theta^m F_{(i,*),(j,*')} \bigg| \leq C \rho^{\frac{1}{2}}e^{-\sigma  s}\| P^\perp_{\mathcal{I}} [w]\|_{H^{m+3}\times H^{m+2}} 
\end{equation}
holds for $s\geq s_0$.   
\end{lemma}
\medskip

\begin{proof} 
From \eqref{eq:betapm}, there is some $C$ depending on $\alpha$ satisfying
\begin{equation}\label{eq:beta_bound}
|\beta_j^*|\leq C(1+j),
\end{equation}
for all $j\geq 0$. Since $h$ is constant, we have
\begin{equation}
\langle (\cos j\theta, \beta_j^* \cos  j\theta),  [\partial_\theta^m  w]\rangle_h=(-1)^m\langle \partial_\theta^m (\cos j\theta, \beta_j^* \cos  j\theta),  [w]\rangle_h,
\end{equation} 
and thus combining \eqref{eq:coefficients} and \eqref{eq:beta_bound} yields
\begin{equation}
C^{-1}e^{-2\sigma s} \|[\partial_\theta^{m}w]\|_h^2\leq  \sum_{j=1}^\infty j^{2m+2} (|A_j^*| +|B_j^*|)^2 \leq Ce^{-2\sigma s} \|[\partial_\theta^{m}w]\|_h^2
\end{equation}
for every $m\geq 1$, where  $C$  depends on $\alpha,m$. Hence, by the Cauchy-Schwarz inequality, we have 
\begin{equation}\label{eq:coeff.CS}
\bigg[\sum_{j=1}^\infty j^{m}(|A_j^*|+|B_j^*|)\bigg]^2\leq \bigg[\sum_{j=1}^\infty j^{2m+2} (|A_j^*| +|B_j^*|)^2\bigg]\bigg[\sum_{j=1}^\infty j^{-2}\bigg]\leq Ce^{-2\sigma s} \|[\partial_\theta^{m}w]\|_h^2
\end{equation}
for each $m \geq 1$. Therefore, remembering \eqref{eq:beta_bound}, for each $m\geq 0$, we have representations of the following functions by absolutely convergent series:
\begin{align}\label{eq:Dw.expansion}
\begin{split}
   \sigma^2 (h e^{ \sigma  s})^{-1}\partial_\theta^m w_s &=\sum_{j=0}^\infty \sum_{* \in \{\pm\}} \beta^*_j \partial_\theta^m  (A_j^*\cos j\theta + B_j^* \sin j\theta),\\
\sigma^2 (h e^{ \sigma  s})^{-1} \partial_\theta^m  (w +w_{\theta\theta})&= \sum_{j=0}^\infty \sum_{ *\in \{\pm\}}  (1-j^2) \partial_\theta^m (A_j^*\cos j\theta + B_j^* \sin j\theta).
\end{split}
\end{align}
Here the equalities hold in the $L^\infty(\mathbb{S}^1)$-sense. Next, we observe
\begin{equation}
2(A^\ast_i \cos i\theta + B^*_i \sin i\theta)(A^{\ast'}_j \cos j\theta + B^{\ast '}_j \sin j\theta)=F_{(i,*),(j,*')}. 
\end{equation} Note that the product of two absolutely convergent series can be written as an absolutely convergent series (of doubled indices). Hence, from \eqref{eq:Dw.expansion} we infer that  for each $m\geq 0$ 
\begin{align}
\sigma^4 (h e^{ \sigma  s})^{-2} \partial_\theta^m [ w_s( w+  w_{\theta\theta})]=&    \lim_{k\to \infty}  \sum_{i,j\leq k}\sum_{*,*'} \tfrac12 \beta^*_i (1-j^2) \partial_\theta^m F_{i.j.*.*'},\\
\sigma^4 (h e^{ \sigma  s})^{-2}\partial_\theta^m [ w_s^2]=&    \lim_{k\to \infty}  \sum_{i,j\leq k}\sum_{*,*'} \tfrac12 \beta^*_i \beta^*_j\partial_\theta^m  F_{i.j.*.*'},
\end{align}
holds in the $L^\infty(\mathbb{S}^1)$-sense. This completes the proof of \eqref{eq:Q.error.expansion}. 

\bigskip

It remains to show \eqref{eq:Q.error.estimate}. By the Cauchy-Schwarz inequality as in \eqref{eq:coeff.CS}, we can obtain
\begin{equation}\label{eq:remainder.L_infty}
\mathcal{F}_{m,\mathcal{I}}:= \sum_{(j,*)\not \in \mathcal{I}}  (j+1)^{m}(|A_{j}^*|+|B_{j}^*|)   \leq  Ce^{-\sigma s}  \|  P^\perp _{\mathcal{I} }[w]\|_{H^{m+1}\times H^m}
\end{equation}
Next, observe a rough estimate
\bea\big|  \partial_\theta^m \bar F_{(i,*),(j,*')} \big|\le C (i+1)^{m+2} (j+1)^{m+2}(|A_i^*|+|B_i^*|)(|A_j^{*'}|+|B_j^{*'}|) ,\eea 
where $\bar F=\sigma^{-1}\beta^*_i   ( 1-j^2  +\frac{  1-\alpha }{2\sigma \alpha  }\beta^{*'}_j  )  F_{(i,*),(j,*')}$, and thus 
\begin{align}\label{eq4115}
    \sum _{((i,*),(j,*'))\notin \mathcal{I}\times \mathcal{I} }
\left| \partial_\theta^m \bar F_{(i,*),(j,*')}\right|&\leq \sum _{(i,*)\notin \mathcal{I}}\sum_{j,*'}
\left| \partial_\theta^m \bar F_{(i,*),(j,*')}\right|+\sum _{i,*}\sum_{(j,*')\notin \mathcal{I}}
\left| \partial_\theta^m \bar F_{(i,*),(j,*')}\right| \\
& \leq C \mathcal{F}_{m+2,\mathcal{I}}  \sum_{i,*}(i+1)^{m+2}(|A_i^{*}|+|B_i^{*}|) \leq C\rho^{\frac{1}{2}} \mathcal{F}_{m+2,\mathcal{I}}.
\end{align}
Now  \eqref{eq:Q.error.estimate} follows by combining  \eqref{eq:remainder.L_infty} and \eqref{eq4115}.
\end{proof}

The following lemma can be shown easily by using the proof of the previous lemma.
\begin{lemma}[Leading asymptotics of cubic error]\label{lem:C.error.expansion}  
 \begin{equation}
\mathcal{C}=-\tfrac{(1-\alpha)(2-\alpha)}{3\alpha^2\sigma}   (A_k^+\cos k\theta +B_k^+ \sin k\theta)^3+o(\rho^{\frac{3}{2}}).
\end{equation}
\end{lemma}

\begin{proof}
In this proof, we fix $\mathcal{I}=\{(k,+)\}$. We recall the definition and estimate of $\mathcal{F}_{m,\mathcal{I}}$ in \eqref{eq:remainder.L_infty}. By the condition \eqref{eq:rho_condition}, we obtain $\mathcal{F}_{m,\mathcal{I}}=o(\rho^{\frac{1}{2}})$, for $m=0,1,2$. Hence, \eqref{eq:Dw.expansion} yields
\begin{align}
   \sigma^2 (h e^{ \sigma  s})^{-1} w_s &=\beta_k^+    (A_k^+\cos k\theta + B_k^+ \sin k\theta)+o(\rho^{\frac{1}{2}}),\\
\sigma^2 (h e^{ \sigma  s})^{-1}   (w +w_{\theta\theta})&=    (1-k^2)  (A_k^+\cos k\theta + B_k^* \sin k\theta)+o(\rho^{\frac{1}{2}}),
\end{align} in $L^\infty(\mathbb{S}^1)$-sense.
Therefore, remembering $\beta_k^+=\sigma$ and $\alpha=k^{-2}$, we can obtain
\begin{equation}
\mathcal{C}= \tfrac{ 1-\alpha}{3\alpha^2\sigma} [\sigma+3\alpha (1-\alpha^{-1})]  (A_k^+\cos k\theta +B_k^+ \sin k\theta)^3+o(\rho^{\frac{3}{2}}).
\end{equation}
This completes the proof, because $\sigma:=1-2\sigma$.
\end{proof}

 \bigskip

\begin{lemma}[Evolution equations]\label{lem:coeff.evol}
For each $j \geq 0$, we have
\begin{align}
\tfrac{d}{ds}A_j^\pm &=(\beta_j^\pm-\sigma)A_j^\pm \mp     (\beta^+_j -\beta^-_j)^{-1} ( \mathfrak{p}_j^c\mathcal{Q}+\mathfrak{p}_j^c\mathcal{C})     +O(\rho^2),\label{eq:A_evol} \\
\tfrac{d}{ds}B_j^\pm  &=(\beta_j^\pm-\sigma)B_j^\pm \mp  (\beta^+_j -\beta^-_j)^{-1} ( \mathfrak{p}_j^s\mathcal{Q}+\mathfrak{p}_j^s\mathcal{C})    +O(\rho^2).\label{eq:B_evol}
\end{align}
In addition, for each $m\in \{ 0,1,2\}$ the following holds almost every $s$
\begin{align}\label{eq:minor.Proj_ineq}
\pm \lambda_\pm^{-1}    \tfrac{d}{ds}     e^{-\sigma  s} \Vert \partial_\theta^m  P_\pm  P^\perp_{\mathcal{J}} [w]   \Vert_h &\geq   e^{-\sigma  s}  \Vert \partial_\theta^m  P_\pm  P^\perp _{\mathcal{J}} [w]  \Vert_h + o(\rho),
\end{align}
where $\lambda^+$ is the least positive eigenvalue of $\mathcal{L}-\sigma $, $-\lambda^-$ is the greatest negative eigenvalue of $\mathcal{L}-\sigma$, and 
\begin{equation}
\mathcal{J}:=\{(k,+),(0,-),(0,+),(2k,-),(2k,+)\}.
\end{equation} 
\end{lemma}

\bigskip

\begin{remark}
Since $\beta_{2k}^-<\beta_0^-<\beta_0^+<\beta_k^+=\sigma  <\beta_{2k}^+$ and $P_+=P_{>\sigma}$, we know \begin{align*}
&P_+  P^\perp _{\mathcal{J}}=P_+-P_{=\beta_{2k}^+}, && P_- P^\perp _{\mathcal{J}}=P_--P_{=\beta_{2k}^-}-P_{=\beta_{0}^-}-P_{=\beta_{0}^+}.
\end{align*}
\end{remark}
\bigskip

\begin{proof}
We recall $\partial_s [w]=\mathcal{L}[w]+(0,E)$, which implies $(\partial_s -\mathcal{L}+\sigma ) e^{-\sigma s}[w]=(0, e^{-\sigma s} E)$. Remembering $( \mathcal{L}-\sigma ) (\cos j\theta, \beta_{j}^* \cos j\theta)=  ( \beta_j^*-\sigma ) (\cos j\theta, \beta_{j}^* \cos j\theta)$, we have
\begin{align}
\frac{d}{ds}A_{j}^*=\frac{\langle \partial_s(  \frac{\sigma^2 }{h e^{\sigma s}}  [w] ),(\cos j\theta, \beta_{j}^* \cos j\theta)\rangle_h }{\| (\cos j\theta, \beta_{j}^* \cos j\theta) \|_h^2}  =(\beta_j^*-\sigma)A_j^* +\frac{\langle (0,\frac{\sigma^2} {h e^{\sigma s}} E ),(\cos j\theta, \beta_{j}^* \cos j\theta)\rangle_h }{\| (\cos j\theta, \beta_{j}^* \cos j\theta) \|_h^2}.
\end{align} 
Proposition \ref{prop:err_rough_rho} implies
\begin{equation}
  -\sigma^2 (  h e^{\sigma s})^{-1} E =   \mathcal{Q} + \mathcal{C}+ O(\rho^2)=    \bigg[\sum_{j=0}^\infty (\mathfrak{p}_{j}^c \mathcal{Q}+\mathfrak{p}_{j}^c \mathcal{C})\cos j\theta+(\mathfrak{p}_{j}^s \mathcal{Q}+\mathfrak{p}_{j}^s \mathcal{C} )\sin j\theta \bigg]+O(\rho^2), 
\end{equation}
resulting in \begin{equation}
\frac{d}{ds}A_{j}^\pm =(\beta_j^\pm -\sigma)A_j^\pm   -\frac{ (\mathfrak{p}_{j}^c \mathcal{Q}+\mathfrak{p}_{j}^c \mathcal{C}) \langle (0,\cos j\theta),( \cos j\theta, \beta_{j}^\pm \cos j\theta)\rangle_h }{\| (\cos j\theta, \beta_{j}^\pm \cos j\theta) \|_h^2}+O(\rho^2).
\end{equation}
As we can write
\begin{align}\label{eq:eigenfunction_decomposition}
\begin{split}
(0,\cos j\theta)& = (\beta^+_j -\beta^-_j)^{-1} \left[(\cos j\theta,\beta^+_j\cos j\theta) -(\cos j\theta,\beta^-_j\cos j\theta)  \right],\\
(0,\sin j\theta) &= (\beta^+_j -\beta^-_j)^{-1} \left[(\sin j\theta,\beta^+_j\sin j\theta) -(\sin j\theta,\beta^-_j\sin j\theta)  \right],
\end{split}
\end{align}this yields \eqref{eq:A_evol}. In the same manner, we can obtain \eqref{eq:B_evol}.

\bigskip

For the last inequality \eqref{eq:minor.Proj_ineq}, we recall $\partial_\theta \mathcal{L}= \mathcal{L} \partial_\theta$ so that $\partial_s \partial_\theta^m [w]=\mathcal{L}\partial_\theta^m [w]+\partial_\theta^m E$. Remembering $\partial_\theta P_\pm P^\perp_{\mathcal{J}} =P_\pm P^\perp_{\mathcal{J}} \partial_\theta $, we can derive the following as in the proof of Proposition \ref{prop:ODE_general} 
\begin{equation}
\tfrac{d}{ds}e^{-2\sigma s} \|\partial_\theta^m  P_+ P^\perp_{\mathcal{J}} [w]\|_h^2 \geq 2 \lambda_+e^{-2\sigma s}  \|  \partial_\theta^m P_+ P^\perp_{\mathcal{J}} [w]\|_h^2+e^{-2\sigma s}\langle (0, \partial_\theta^m E),	\partial_\theta^m P_+ P^\perp_{\mathcal{J}} [w]\rangle_h.
\end{equation}
In view of the first estimate in Proposition \ref{prop:err_rough_rho}, the cubic error $\mathcal{C}$ (see \eqref{eqCubic}) has an estimate $\partial^m_\theta \mathcal{C}= O(\rho^{\frac32})$. Then the other estimate in Proposition \ref{prop:err_rough_rho} implies  
\begin{equation}
\sigma^2(  h e^{\delta s} )^{-1}\partial_\theta^m E=-\partial_\theta^m \mathcal{Q}+O(\rho^{\frac{3}{2}}). 
\end{equation}
Next,  by applying \eqref{eq:Q.error.estimate}  with 
$\mathcal{I}=\{(k,+)\}$ and using $\beta^+_k =\sigma$, we have 
\bea  \label{eq-4130} -\partial^m_\theta \mathcal{Q} =\tfrac{ 1-\alpha }{2\alpha}\partial_\theta^m F_{(k,+),(k,+)}+o(\rho) \eea 
for $m\leq 2$.
Note here we used the condition \eqref{eq:rho_condition} which implies  $\Vert P^\perp_\mathcal{I}[w]\Vert_{H^{m+3}\times H^{m+2}}=o(\rho ^{\frac 12})$ for $m\leq 2$.
Finally, \eqref{eq:eigenfunction_decomposition} implies 
\begin{equation} \label{eq-4131}
(0,\partial_\theta^m F_{(k,+),(k,+)}) \in \text{span}\{ (1,\beta_0^*), (\cos 2k\theta, \beta_{2k}^* \cos 2k\theta),(\sin 2k\theta, \beta_{2k}^* \sin 2k\theta): *=\pm \},
\end{equation}
namely $ \langle (0,\partial_\theta^m F_{(k,+),(k,+)}),  P_{+}  P^\perp_{\mathcal{J}} [w]\rangle_h=0$. Therefore, combining above, we infer \eqref{eq:minor.Proj_ineq} holds for $*=+$.

\bigskip

In the same manner, we can obtain
\begin{equation}
     \tfrac{d}{ds}     e^{-\sigma  s} \| \partial_\theta^m P_-  P^\perp_{\mathcal{J}} [w]   \|_h  \leq  -\lambda_- e^{-\sigma  s}  \| \partial_\theta^m P_-  P^\perp_{\mathcal{J}}   [w]  \|_h + o(\rho).
\end{equation}
Since $\lambda_->0$, this inequality is equivalent to \eqref{eq:minor.Proj_ineq} with $*=-$. 
\end{proof}

\bigskip

\begin{proposition}[Minor coefficient estimates]\label{prop:coeff.rough.bound}
 For each $(j,*)\neq (k,+)$, the following holds
\begin{equation}
|A_j^*|+|B_j^*|=O(\rho).
\end{equation}
\end{proposition}  

\begin{proof}
By using Proposition \ref{prop:err_rough_rho}, we derive from Lemma \ref{lem:coeff.evol} the following
\begin{align}
&\tfrac{d}{ds}A^\pm_j=(\beta_j^\pm-\sigma)A^\pm_j+O(\rho), && \tfrac{d}{ds}B^\pm_j=(\beta_j^\pm-\sigma)B^\pm_j+O(\rho).
\end{align}
Since $\beta_j^*\neq \sigma$ for $(j,*)\neq (k,+)$, the condition \eqref{eq:rho_condition} allows us to apply Lemma \ref{lem:ODE.large.O} with $\tilde{\rho}=\rho$ so that we can complete the proof.
\end{proof}

\bigskip

\begin{lemma}[Leading asymptotics of quadratic error]\label{lem:Q.fine.expansion}
The quadratic error $\mathcal{Q}$ satisfies
\begin{equation}
  \mathcal{Q} =  \sum_{(i,*),(j,*')\in \mathcal{J}}  \sigma^{-1}\beta^*_i   \big( 1-j^2  +\tfrac{  1-\alpha }{2\sigma \alpha  }\beta^{*'}_j   \big) F_{(i,*),(j,*')}+ o(\rho^{\frac{3}{2}}).
\end{equation}
where $\mathcal{J}:=\{(k,+),(0,+),(0,-),(2k,+),(2k,-)\}$.
\end{lemma}

\begin{proof}
Since we have the inequality \eqref{eq:minor.Proj_ineq} and the condition \eqref{eq:rho_condition}, applying  Lemma \ref{lem:ODE.small.O} with $\tilde{\rho}=\rho$  yields
\begin{equation}
e^{-2\sigma  s}\|     P^\perp_{\mathcal{J}} [w]   \|_{H^3\times H^2 }^2 = \sum_{* \in \{\pm\}}e^{-2\sigma s} \|   P_*  P^\perp _{\mathcal{J}} [w]   \|_{H^3\times H^2 }^2 \leq C \sum_{m=0}^{2}\sum_{* \in \{\pm\}}e^{-2\sigma s} \| \partial_\theta^m  P_*  P^\perp _{\mathcal{J}} [w]   \|_h^2=o(\rho^2).
\end{equation} 
Hence, \eqref{eq:Q.error.estimate} in Lemma \ref{lem:Q.error.expansion} completes the proof.
\end{proof}
  
  \bigskip

 \begin{lemma}[Projections of quadratic error]\label{lem:Q.projections}
The quadratic error $\mathcal{Q}$ satisfies
\begin{align}
 \mathfrak{p}_0^c\mathcal{Q}&=- \tfrac{ k^2-1  }{2}\rho +o(\rho^{\frac{3}{2}}),\\
   \mathfrak{p}_{2k}^c\mathcal{Q}&=- \tfrac{ k^2-1  }{2}(|A^+_k|^2-|B^+_k|^2) +o(\rho^{\frac{3}{2}}),\\
   \mathfrak{p}_{2k}^s\mathcal{Q}&=- ( k^2-1 ) A^+_k B^+_k+o(\rho^{\frac{3}{2}}),
\end{align}
and
\begin{align} 
 \mathfrak{p}_{k}^c   \mathcal{Q}& =  2 A^+_kA_0   - (4k^2-1)   (A^+_k A_{2k}  +B^+_k B_{2k})     +o(\rho^{\frac{3}{2}}),\\
  \mathfrak{p}_{k}^s   \mathcal{Q}& =  2 B^+_k A_0 - (4k^2-1)(A^+_k B_{2k} -B^+_kA_{2k})       +o(\rho^{\frac{3}{2}}).
\end{align}
 \end{lemma}
 
 \begin{proof}
By Lemma \ref{lem:Q.fine.expansion}, we only need to consider $F_{(i,*),(j,*')}$ with $(i,*),(j,*')\in \mathcal{J}$. In addition, Proposition \ref{prop:coeff.rough.bound} gives us $F_{(i,*),(j,*')}=O(\rho^2)$ when both of $(i,*)$ and $(j,*')$ are not equal to $(k,+)$. Hence, it is enough to consider $F_{(k,+),(j,*)}$ and $F_{(j,*),(k,+)}$ with $(j,*)\in \mathcal{J}$. We define 
\begin{equation}
 \gamma_{(i,*),(j,*')}:=\sigma^{-1}\beta_i^*\left(1-j^2+\tfrac{1-\alpha}{2\sigma \alpha}\beta_j^{*'}\right)+\sigma^{-1}\beta_j^{*'}\left(1-i^2+\tfrac{1-\alpha}{2\sigma \alpha}\beta_i^{*}\right).
\end{equation}
Then, \eqref{eq-Fvec} yields
\begin{align}
 \mathfrak{p}_0^c\mathcal{Q}&= \tfrac12 \gamma_{(k,+),(k,+)}(|A_k^+|^2 +|B_k^+|^2) +o(\rho^{\frac{3}{2}}),\\
 \mathfrak{p}_{2k}^c\mathcal{Q}&=\tfrac12 \gamma_{(k,+),(k,+)}(|A_k^+|^2 -|B_k^+|^2) +o(\rho^{\frac{3}{2}}),\\
 \mathfrak{p}_{2k}^s\mathcal{Q}&= \gamma_{(k,+),(k,+)} A_k^+B_k^{+}  +o(\rho^{\frac{3}{2}}).
\end{align}
 Since $\beta_k^+=\sigma$ and $\alpha=k^{-2}$, we have $\gamma_{(k,+),(k,+)}=1- k^2$. This completes the proof of the first three equations. To show the last two equations, we again derive from  Proposition \ref{prop:coeff.rough.bound} and Lemma \ref{lem:Q.fine.expansion} that
 \begin{align}
 \mathfrak{p}_{k}^c\mathcal{Q}&=  \sum_* 2 \gamma_{(k,+),(0,*)}  A_k^+A_0^*+ \sum_*\gamma_{(k,+),(2k,*)}  (A_k^+A_{2k}^*+B_k^+B_{2k}^*)  +o(\rho^{\frac{3}{2}}),\\
 \mathfrak{p}_{k}^s\mathcal{Q}& = \sum_* 2 \gamma_{(k,+),(0,*)}  B_k^+A_0^* +\sum_* \gamma_{(k,+),(2k,*)} (A_k^+B_{2k}^*-B_k^+A_{2k}^*)  +o(\rho^{\frac{3}{2}}).
\end{align}
Since we have $\sum_*A_j^*=A_j$ and $\sum_*B_j^*=B_j$ by definition, we can obtain the desired result by observing
 \begin{align}
\gamma_{(k,+),(0,*)}  &=1+\tfrac{1-\alpha}{2\sigma \alpha}\beta_0^*+\sigma^{-1}\beta_0^*\left(1-k^2+\tfrac{1-\alpha}{2\alpha}\right)=1,\\
\gamma_{(k,+),(2k,*)}  &=1-4k^2+\tfrac{1-\alpha}{2\sigma \alpha}\beta_{2k}^*+\sigma^{-1}\beta_{2k}^*\left(1-k^2+\tfrac{1-\alpha}{2\alpha}\right)=-(4k^2-1).
\end{align}
  \end{proof}

 \begin{lemma}[Projections of cubic error]\label{lem:C.projections}
The cubic error $\mathcal{C}$ satisfies
\begin{align}\label{eq-C2}
&\mathfrak{p}_{k}^c   \mathcal{C} =   - \tfrac{(k^2-1)(2k^2-1)}{4\sigma}\rho  A^+_k  +o(\rho^{\frac{3}{2}}) , && \mathfrak{p}_{k}^s   \mathcal{C} =   - \tfrac{(k^2-1)(2k^2-1)}{4\sigma}\rho   B^+_k  +o(\rho^{\frac{3}{2}}).
\end{align}
and $|\mathfrak{p}_{0}^c   \mathcal{C}|+|\mathfrak{p}_{2k}^c   \mathcal{C}|+|\mathfrak{p}_{2k}^s   \mathcal{C}|=o(\rho^{\frac{3}{2}})$.
 \end{lemma}

\begin{proof}
By using the identity $(\cos k\theta)^3 =\frac{1}{4}(3\cos k \theta+\cos 3k\theta)$, we obtain
\begin{align}\label{eq:tri.triple.identity}
\begin{split}
(A_k^+)^3 (\cos  k\theta)^3+3A_k^+|B^+_k|^2 \cos  k\theta (\sin k\theta)^2= &\, 3A_k^+|B_k^+|^2 \cos k\theta+A^+_k(|A^+_k|^2-3|B^+_k|^2)(\cos k\theta)^3\\
=&\,\tfrac{3}{4}\rho A^+_k\cos k\theta+ \tfrac{1}{4}A^+_k(|A^+_k|^2-3|B^+_k|^2)\cos 3k\theta.
\end{split}
\end{align} 
We have a similar identity for $(B^+_k)^3 (\sin k\theta)^3+ 3 B^+_k |A^+_k|^2 \sin k\theta (\cos k\theta)^2 $, and then  \eqref{eq-C2} follows from  Lemma \ref{lem:C.error.expansion}.  In addition, the previous two identities  imply that $(A_k^+ \cos k\theta + B_k^+ \sin k\theta)^3$ is orthogonal to the space spanned by $1,\cos 2k\theta, \sin 2k\theta$, which shows $|\mathfrak{p}_{0}^c   \mathcal{C}|+|\mathfrak{p}_{2k}^c   \mathcal{C}|+|\mathfrak{p}_{2k}^s   \mathcal{C}|=o(\rho^{\frac{3}{2}})$.
\end{proof}

 \bigskip

 \begin{lemma} [Dynamics among main coefficients] \label{lem:main.dynamics}
As $s\to \infty$, the following hold
\begin{align}
 \label{eq:evol.A_0} \frac{d}{ds} A^\pm_0& = (\beta^\pm_0-\sigma  ) A^\pm_0 \pm \frac{ k^2-1 }{2 (\beta^+_0-\beta^-_0)} \rho +o(\rho^{\frac{3}{2}}), \\
 \label{eq:evol.A_2k}\frac{d}{ds} A^\pm_{2k} &= (\beta^\pm_{2k}-\sigma  ) A^\pm_{2k} \pm \frac{k^2-1 }{2  (\beta^+_{2k}-\beta^-_{2k})} (|A^+_k|^2-|B^+_k|^2) +o(\rho^{\frac{3}{2}}), \\
 \label{eq:evol.B_2k}\frac{d}{ds} B^\pm_{2k} &= (\beta^\pm_{2k}-\sigma  ) B^\pm_{2k} \pm \frac{k^2-1 }{ \beta^+_{2k}-\beta^-_{2k} } A^+_kB^+_k + o(\rho^{\frac{3}{2}}),\\
 \label{eq:evol.A_k}  \frac{d}{ds} A^+_k& = \frac{    4k^2-1  }{\beta^+_{k}-\beta^-_{k}}(A^+_k A_{2k}   +  B^+_k B_{2k})  + \frac{ (k^2-1)(2k^2-1)\rho -8\sigma A_0  }{4\sigma (\beta^+_{k}-\beta^-_{k})}  A^+_k   + o(\rho^{\frac{3}{2}}),\\
\label{eq:evol.B_k} \frac{d}{ds} B^+_k& =  \frac{     4k^2-1 }{\beta^+_{k}-\beta^-_{k}}(A^+_k B_{2k}-B^+_kA_{2k})  + \frac{ (k^2-1)(2k^2-1) \rho  -8\sigma  A_0 }{4\sigma (\beta^+_{k}-\beta^-_{k})} B^+_k   + o(\rho^{\frac{3}{2}}).
\end{align}
\end{lemma}

\begin{proof}
Combine Lemma \ref{lem:coeff.evol}, Lemma \ref{lem:Q.projections}, and Lemma \ref{lem:C.projections}.
\end{proof}

 \bigskip

 \begin{lemma}The functions  $Q^\pm:= (|A^+_k|^2 -|B^+_k|^2) A^\pm_{2k}+  2 A^+_k B^+_k B^\pm_{2k}$ and $Q := Q^++Q^-$ satisfy
 \begin{align}
\label{eq:evol.rho}\frac{d}{ds}\rho&= \frac{  2(4k^2-1)  }{\beta^+_{k}-\beta^-_{k}}Q  + \frac{ (k^2-1)(2k^2-1) \rho -8\sigma A_0 }{2\sigma (\beta^+_{k}-\beta^-_{k})}\rho    + o(\rho^2),\\
\label{eq:evol.Q_pm}\frac{d}{ds} Q^\pm &=(\beta^\pm_{2k}-\sigma )Q^\pm \pm  \frac{  k^2-1}{2 (\beta^+_{2k}-\beta^-_{2k})}\rho^2 +o(\rho^{\frac{5}{2}}).
 \end{align}
 \end{lemma}

\begin{proof}
Since $\rho=|A_k^+|^2+|B_k^+|^2$, the equations \eqref{eq:evol.A_k} and \eqref{eq:evol.B_k} immediately imply \eqref{eq:evol.rho}.
 
To derive \eqref{eq:evol.Q_pm}, we first observe $|\frac{d}{ds}A_k^+|+|\frac{d}{ds}B_k^+|= O(\rho^{\frac{3}{2}})$ from \eqref{eq:evol.A_k}, \eqref{eq:evol.B_k}, and Proposition \ref{prop:coeff.rough.bound}. Therefore, by using Proposition \ref{prop:coeff.rough.bound} we obtain
\begin{equation}
\tfrac{d}{ds} Q^\pm = (|A^+_k|^2 -|B^+_k|^2) \tfrac{d}{ds} A^\pm_{2k}+  2 A^+_k B^+_k \tfrac{d}{ds} B^\pm_{2k} +O(\rho^3).
\end{equation}
Thus, combining with \eqref{eq:evol.A_2k} and \eqref{eq:evol.B_2k} completes the proof. 
\end{proof}

  \bigskip

  \begin{lemma}The functions $A^\pm_0$ and $Q^\pm$ can be represented only by using $\rho$ as follows:
  \begin{align}
\label{eq:A_0^pm.estimates}   A^\pm_0 &=\mp  \frac{ k^2-1 }{2 (\beta^\pm_0-\sigma  ) (\beta^+_0-\beta^-_0)} \rho+O(\rho^{\frac{3}{2}}),\\
\label{eq:Q^pm.estimates}   Q^\pm &=\mp \frac{k^2-1 }{2 (\beta_{2k}^\pm-\sigma) (\beta^+_{2k}-\beta^-_{2k})} \rho^2+O(\rho^{\frac{5}{2}}).
  \end{align}  
  \end{lemma}

 \begin{proof}
Notice that  Proposition \ref{prop:coeff.rough.bound} and \eqref{eq:evol.rho} imply $\frac{d}{ds}\rho = O(\rho^2)$. Hence, \eqref{eq:evol.A_0} gives us
 \begin{equation}
 (\beta^\pm_0-\sigma  )^{-1} \frac{d}{ds} \left[A^\pm_0 \pm   \frac{ k^2-1 }{2 (\beta^\pm_0-\sigma  ) (\beta^+_0-\beta^-_0)} \rho\right]  =  A^\pm_0 \pm \frac{ k^2-1 }{2(\beta^\pm_0-\sigma  )  (\beta^+_0-\beta^-_0)} \rho  +o(\rho^{\frac{3}{2}}).
 \end{equation}
Therefore, remembering the condition \eqref{eq:rho_condition}, we can apply Lemma \ref{lem:ODE.large.O} with $\tilde{\rho}=\rho^{\frac{3}{2}}$ to obtain \eqref{eq:A_0^pm.estimates}. Similarly, using $\frac{d}{ds}\rho = O(\rho^2)$ and \eqref{eq:evol.Q_pm} we can derive
   \begin{equation}
(\beta_{2k}^\pm-\sigma)^{-1}\frac{d}{ds} \left[ Q^\pm \pm \frac{k^2-1 }{2 (\beta_{2k}^\pm-\sigma) (\beta^+_{2k}-\beta^-_{2k})} \rho^2\right] = Q^\pm \pm \frac{k^2-1 }{2 (\beta_{2k}^\pm-\sigma) (\beta^+_{2k}-\beta^-_{2k})} \rho^2+ o(\rho^{\frac{5}{2}}).
\end{equation}  
 Hence, Lemma \ref{lem:ODE.large.O} with $\tilde{\rho}=\rho^{\frac{5}{2}}$ provides \eqref{eq:Q^pm.estimates}.
   \end{proof}
 
 \bigskip
 
 \begin{proof}[Proof of Lemma \ref{lem:round_exp_goal}] By \eqref{eq:A_0^pm.estimates}, the function $A_0:=A_0^++A_0^-$ satisfies
 \begin{align}
   A_0  =  \frac{ (k^2-1)\rho }{2  (\beta^+_0-\beta^-_0)} \left[  \frac{1}{ \beta^-_0-\sigma  }-\frac{1}{ \beta^+_0-\sigma  }\right] +O(\rho^{\frac{3}{2}})=  \frac{ (k^2-1)\rho }{2 (\beta^+_0-\sigma  )(\beta^-_0-\sigma  ) } +O(\rho^{\frac{3}{2}}).
  \end{align}
Since $\beta_j^\pm$ are the two solutions to $\beta^2+\beta-2\alpha \sigma (j^2-1)$, we have
\begin{equation}\label{eq:beta_roots}
(\beta^+_j-\sigma  ) (\beta^-_j-\sigma  ) = \beta^+_j\beta^-_j-(\beta^+_j+\beta^-_j)\sigma   + \sigma  ^2 =   -2\sigma \alpha  (j^2-1) + \sigma + \sigma (1-2\alpha)=2\sigma  (1- \alpha  j^2).  
\end{equation}
Hence,
\begin{align}
   A_0  =    \frac{  k^2-1  }{4\sigma }\rho +O(\rho^{\frac{3}{2}}).
\end{align}
In the same manner, by using \eqref{eq:Q^pm.estimates}, \eqref{eq:beta_roots}, and $\alpha=k^{-2}$, we can derive
  \begin{align}
   Q  =   \frac{  (k^2-1)\rho^2 }{2 (\beta^+_{2k}-\sigma  )(\beta^-_{2k}-\sigma  ) }   +O(\rho^{\frac{5}{2}})=-\frac{  k^2-1  }{12\sigma }\rho^2 +O(\rho^{\frac{5}{2}}).
  \end{align}
  Therefore, we can update the ODE \eqref{eq:evol.rho} as follows:
\begin{align}
   \frac{d}{ds}\rho   = \frac{k^2-1}{  \beta^+_{k}-\beta^-_{k}}\left[ -\frac{ 4k^2-1 }{6\sigma} + \frac{  (2k^2-1)-  2 }{2\sigma  }   + o(1) \right]\rho^2    =\frac{k^2-1}{  \beta^+_{k}-\beta^-_{k}}\left[\,\frac{    k^2-4 }{3\sigma }    + o(1)\right]\rho^2 . 
  \end{align}
 Since $k\geq 3$ and  $\beta^+_{k}>\beta^-_{k}$, there are some $c_0,s_1>0$ such that $\frac{d}{ds}\rho\geq c_0\rho^2$ holds for $s\geq s_1$. In particular, $\rho$ can not converge to zero and this contradicts to the condition \eqref{eq:rho_condition}.  
\end{proof}

\medskip

Combining main results from two sub-sections, we can now give the proof of the uniqueness of blow-downs and the fast convergence. 
\begin{proof}[Proof of Theorem \ref{thm-42}] 
The proof directly follows by combining Theorems \ref{thm-exponential1} and Theorem \ref{thm-exponential2}. 
\end{proof}


\subsection{Classification under unique blow-down assumption}  \label{sec-uniqueness} 


Let $S:[l_0,+\infty)\times \mathbb{S}^1\to \mathbb{R}$ be a smooth solution to \eqref{Se}  which has the unique limit shape represented by a  solution $h \in C^\infty(\mathbb{S}^1)$ to \eqref{eq-shrinker}.
 Then  Theorem \ref{thm-42} implies a convergence with a  decay estimate: 
\bea \label{eq-75}
{ \Vert w(s,\cdot) \Vert_{C^0(\mathbb{S}^1)} = \Vert S(e^s ,\cdot)- \sigma^{-1} e^{\sigma s} h  \Vert_{C^0(\mathbb{S}^1)} = o(e^{(\sigma-\delta)s})\quad \hbox{ as $ s\to \infty$ }   }
\eea 
  for some  $\delta>0$. Here,    $w$ solves \eqref{eq-w'}. 
We first    improve the decay estimate  by finding next leading coefficients in  the following proposition.  
The proof of Theorem \ref{thm-72} will be based on  the iteration of    improved decay estimates.    In what follows, we will use the spectral analysis and the notations introduced in Sections \ref{subsection-ODE} and \ref{sec:set.up.spect}. 
 
To begin with, we recall the definition  \eqref{eq-K} of the number $K$ that
\begin{equation}
    \beta_0^+\leq \cdots \leq \beta_{K-1}^+<\sigma \leq \beta_K^+.
\end{equation}

\begin{proposition} [Improved decay estimate]\label{prop-coeff-extract} 
Given a fixed $\eta \in \{0,1\}$ and a solution $h \in C^\infty(\mathbb{S}^1)$ to \eqref{eq-shrinker}, we let $S_1,S_2 $ be smooth solutions to \eqref{Se}  with the unique limit shape represented by    $h \in C^\infty(\mathbb{S}^1)$. 
 Suppose that for each $i=1,2$,  $w_i(s,\theta):= S_i(e^s ,\cdot)- \sigma^{-1} e^{\sigma s} h(\theta) $ satisfies  \eqref{eq-75} and 
 \be
 \Vert w_1(s,\cdot)-w_2(s,\cdot) \Vert _{C^0(\mathbb{S}^1)} =o(e^{ \gamma s})\quad\hbox{ as $ s\to \infty$ }
\ee 
for some $\beta ^+_m<{ {\gamma}} <\beta ^+_{m+1}$ with $m\in\{0,\ldots,K-1\}$. 
Let $N$ be the multiplicity of the eigenvalue $\beta^+_m$ in the spectrum \eqref{eq-basis}. 
Then,  there are constants $C_1$, $C_2$, $\ldots$, $C_{N}$ $ \in \mathbb{R}$ and $\e>0$ such that 
\bea 
  w_1(s,\cdot)-w_2(s,\cdot)= \sum_{i=1}^N C_i\varphi_{m-i+1}  e^{\beta^+_ms} + o\left(e^{(\beta^+_{m}-\e)s}\right)  \eea 
uniformly for  $ \theta\in \mathbb{S}^1$ as $s\to\infty $.

\begin{proof}
Although the theorem is written for two   cases $\eta=0$ and $\eta=1$, the same proof applies to both cases. So, we suppress the subscript $\eta$ and write one proof as we did in the other proofs in this section. 
 \smallskip

Let us first present a priori estimate of a   regularity improvement,      which we will   use repeatedly  in the proof. 
Suppose  that 
\be
{\Vert w_1(s,\cdot)-w_2(s,\cdot) \Vert _{L^2_h(\mathbb{S}^1)}} =o(e^{\tilde\sigma s})\quad\hbox{ as $ s\to \infty$ }
\ee 
 for some $\tilde \sigma<\sigma =1-2\alpha$.  Then we  have     
  \be \label{eq-reg-impr}
  \Vert w_1-w_2 \Vert_{C^{2,{{\beta}},\tilde\sigma   }_R} =o(1)\quad\hbox{ as $R\to \infty$ }
   \ee
   for   $0<\beta<1$ 
  by  \cite[Lemma B.1 (interior regularity)]{CCK_existence} which  follows from the elliptic regularity theory. Let us define 
\bea\label{eq-beta'} 
\tilde {\gamma}:= \inf \{\, \tilde \sigma   \,:\,  {\Vert w_1(s,\cdot)-w_2(s,\cdot) \Vert _{L^2_h(\mathbb{S}^1)}}= o(e^{\tilde \sigma  s}) \,\, \text{ as }s\to\infty \}. 
\eea 
Note   that we have $\tilde \gamma  \le \gamma <\sigma  $ by the assumption.   In light of   the definition of $\tilde\gamma$  and  the regularity improvement   as in \eqref{eq-reg-impr},   any $\delta>0$ satisfies
\bea\label{eq-770} 
\Vert w_1 -w_2 \Vert_{C^{2,\beta,\tilde \gamma +\delta }_R} = o(1)\quad\hbox{ as $R\to \infty$.} 
\eea 

\medskip
We may assume that  $\tilde \gamma \ge \beta^+_m$.  Otherwise, it  follows from          the estimate \eqref{eq-770}    that 
\bea 
\Vert w_1 -w_2 \Vert_{C^{2,\beta, \beta^+_m-\e }_R} = o(1)\quad\hbox{ as $R\to \infty$ } 
\eea  
for a  sufficiently small $\e>0$
 since        $\tilde \gamma < \beta^+_m$.  Thus the proposition holds true with $C_i=0$.
 \medskip
 
 Now   suppose   that $\tilde\gamma \ge \beta^+_m$.    Then by     \eqref{eq-75}    and  \cite[Lemma 3.3]{CCK_existence} (error estimate), we obtain that 
 \be 
\begin{aligned}
 \Vert E( w_1)-E(w_2) \Vert_{C_R^{0,\beta, \tilde\gamma-\delta}}
&\le  
C\left( e^{ -\delta R} + e^{(2\delta-4\alpha)R}\right)   \|w_1-w_2\|_{C_R^{2,\beta, \tilde\gamma+\delta}}  \end{aligned}
\ee 
for a sufficiently small $ \delta >0$ with some constant $C>0$ since $\tilde \gamma \le \gamma<\sigma$.  Thus   \eqref{eq-770} implies that 
\bea \label{eq-77}
\Vert E(w_1)-E(w_2) \Vert_{C^{0,\beta,\tilde \gamma -\delta }_R}=o(1) \quad\hbox{ as $R\to \infty$ }  
\eea 
for a sufficiently small $ \delta >0$. 
We define $v=w_1-w_2$ and denote
\bea 
x(s):= \Vert P_{>\beta_{m}^+}\,[ v ]\, \Vert_h ,\qquad y(s):= \Vert P_{=\beta_{m}^+}\, [ v ]\, \Vert_h ,\qquad z(s):= \Vert P_{<\beta_{m}^+}\,[ v ]\,  \Vert_h.
\eea
Then, $[v]=(v,\p_s v)$ solves 
\bea\label{eq-vectorwdifference}
 {\p_s}[v]=\mathcal{L}[v]+(0,E (w_1)-E (w_2) ),
\eea
and therefore Proposition  \ref{prop:ODE_general} yields
\be \label{eq--odesystem'} 
\left\{
\ba
x'-\beta^+_{m+1} x \,&\,\ge  -  \Vert (0,E (w_1)-E (w_2)) \Vert_h,  \\
|y'-\beta^+_{m} y |  &\,\le    \Vert (0,E (w_1)-E (w_2)) \Vert_h, 
 \\z' -\beta^+_{m-N} z\, &\,\le  \Vert (0,E (w_1)-E (w_2)) \Vert_h.
 \ea
 \right.
 \ee
 Here, we used the convention $\beta^+_{-1}:= \beta^-_0$ in the case   $m-N=-1$.

\medskip

Next, we will prove that $\tilde \gamma=\beta^+_m$.  Suppose to the contrary that 
 $\tilde \gamma>\beta^+_m$.  
  Let us choose a sufficiently small $\delta>0$ such that    
\bea\label{eq-choice-delta}
 \beta^+_m < \tilde \gamma-\delta < \tilde \gamma+\delta <\beta^+_{m+1} 
 \eea
and \eqref{eq-77} holds.  
Let 
 \be
 x_0(s):=  e^{-(\tilde \gamma+\delta)s}  x(s)\qquad \hbox{and}\qquad z_0(s):= e^{-( \tilde \gamma+\delta)s}\left(  y(s)+z(s)+e^{( \tilde \gamma-\delta)s}\right) . 
 \ee
 In light of \eqref{eq-770}, we have 
 \be
e^{-(\tilde \gamma+\delta)s} \Vert  [ v]\, \Vert_h= o(1)\quad \text{ as }s\to\infty ,
 \ee
and hence  
it follows   that 
\bea
 x_0+z_0 \to 0\quad \text{ as }s\to\infty
\eea 
by    the orthogonality of   $ \left\{  (\varphi_j, \beta^\pm_j \varphi_j) \right\}_{j=0}^\infty $ with respect to $\langle\cdot,\cdot\rangle_h$. 
Moreover,    we have 
\be\label{eq-ode-z-00}
\left\{
\ba x_0'  &\ge (\beta^+_{m+1}-\tilde\gamma-\delta)x_0- e^{-(\tilde \gamma+\delta)s}   \Vert (0,E (w_1)-E (w_2)) \Vert_h, \\  
z_0'  &\le \max(  \beta^+_{m}-\tilde\gamma-\delta,- 2\delta )z_0 + (N+1)e^{-(\tilde \gamma+\delta)s}   \Vert (0,E (w_1)-E (w_2)) \Vert_h .
 \ea
 \right.
 \ee 
Since Proposition \ref{prop-normequivalence} implies
 \be 
 \Vert (0,E (w_1)-E (w_2)) \Vert_h \leq C  \Vert E (w_1)-E (w_2)  \Vert_{L^2_h(\mathbb{S}^1)},
 \ee
the estimate \eqref{eq-77} implies  that 
  \be\label{eq-est-err-comp-h}
 \Vert (0,E (w_1)-E (w_2)) \Vert_h =  o\left(e^{(\tilde\gamma -\delta)s}\right)\quad \hbox{as $s\to\infty$},
 \ee
 and hence
\be\label{eq-ode-z-0}
\left\{
\ba x_0'  &\ge (\beta^+_{m+1}-\tilde\gamma-\delta)x_0-o(z_0), \\  
z_0'  &\le \max(  \beta^+_{m}-\tilde\gamma-\delta,- 2\delta )z_0-o(z_0)= -2\delta z_0 - o(z_0)
 \ea
 \right.
 \ee 
as $s\to\infty $. Here, we have chosen a sufficiently small    $\delta$ as in \eqref{eq-choice-delta}. 
Then utilizing Lemma \ref{lem-MZODE}, we obtain   that  
$
 x_0=o(z_0) $  as $s\to\infty $ since  $\beta^+_{m+1}-\tilde\gamma-\delta>0$. 
  This  yields   that 
  \be
   x_0+z_0 = o\left(e^{ -\frac 32\delta s}\right)\quad \hbox{as $s\to\infty $ }
   \ee 
because $z_0'  \le  -\frac 32\delta  z_0 
$ for large $s$ in light of \eqref{eq-ode-z-0}. 
Thus, it follows that 
$
\Vert  [v]  \Vert _h= o (e^{(\tilde\gamma- \frac12 \delta) s} ),$ 
 and hence by   Proposition \ref{prop-normequivalence}, we deduce  that 
 \be
 \Vert v (s,\cdot )  \Vert _{L^2_h(\mathbb{S})} = o\left(e^{(\tilde\gamma- \frac12 \delta) s}\right)\quad \hbox{as $s\to\infty $. }
 \ee This contradicts the definition of $\tilde\gamma$. 
 Therefore, we conclude that  $\tilde\gamma=\beta^+_m$.

\medskip

As seen  in \eqref{eq-77} and \eqref{eq-est-err-comp-h},   it holds that 
  $\Vert (0,E (w_1)-E (w_2))\Vert_h(s)= o\left(e^{(\beta^+_m -\delta)s}\right)$ 
  for a sufficiently small $\delta>0$ since $\tilde \gamma= \beta^+_m$. 
   By solving  the ODE for $y$ given in \eqref{eq--odesystem'}, we 
 obtain that 
 \bea\label{eq-yi} 
 y = c e^{\beta^+_ms}+ o\left(e^{(\beta^+_m-\delta)s}\right)\quad \hbox{as $s\to\infty $  }
 \eea 
with  some $c\ge0$.
Letting 
 \be
 x_1(s) = e^{-(\beta^+_m+\delta)s} x(s) \quad\hbox{  and } \quad z_1(s) = e^{-(\beta^+_m+\delta)s}\left(z(s) + e^{(\beta^+_m-\delta)s }\right) ,
 \ee
a similar argument as for \eqref{eq-ode-z-0} shows that  
 \be
 \left\{
 \ba
 x_1'  &\ge (\beta^+_{m+1}-\beta^+_m -\delta )x_1-o(z_1), \\  z_1'  &\le\max ( \beta^+_{m-N}-\beta^+_m-\delta,-2\delta )z_1-o(z_1)\leq -2\delta z_1 - o( z_1),
 \ea
 \right.
 \ee
 provided  that  $\delta>0$ is  sufficiently small.
 Then using 
 Lemma \ref{lem-MZODE} again  implies that $x_1=o(z_1)$, and hence  $x_1+z_1= o(e^{-\frac32 \delta s})$   since $z_1'  \le  -\frac 32\delta  z_1
$ for large $s$. 
 Thus  we obtain that
  \bea\label{eq-xz} 
  x+z = o\left(e^{(\beta^+_m- \frac12 \delta)s}\right) \quad \hbox{as $s\to\infty $}.
  \eea 
Therefore, combining \eqref{eq-yi} and \eqref{eq-xz} yields that 
\bea 
\left\Vert v(s,\cdot) -\sum_{i=1}^N C_{i} \varphi_{m-i+1} e^{\beta^+_ms} \right\Vert_h = o\left(e ^{(\beta^+_m-\frac12\delta)s}\right) 
\eea
as $s\to\infty $ 
 with  some  constants $C_i$ (for $i=1,\ldots,N$). 
 This proves the proposition in light of Proposition \ref{prop-normequivalence} 
  since $\Vert  f \Vert _{C^0(\mathbb{S}^1)} \le C\Vert  f \Vert _{H^1(\mathbb{S}^1)}   $  for any $f\in H^1(\mathbb{S}^1)$  by the Sobolev embedding. 
 \end{proof}

\end{proposition}

\medskip

Next, we are ready to prove the main classification  result in  Theorem \ref{thm-72}. 
The proof will   be a nonlinear analogue of the Gram-Schmidt process of finding correct `coordinates' $\mathbf{y}=(y_0,\ldots, y_{K-1})\in \mathbb{R}^K$ by an iteration of Proposition \ref{prop-coeff-extract}. 

\begin{proof}[Proof of Theorem \ref{thm-72}] 
Let  $\Sigma $ be a  given translator  with the unique limit shape represented by  a solution  $h\in C^\infty(\mathbb{S}^1)$  to \eqref{eq-shrinker}. Then, $\Sigma $  is the graph of a strictly convex entire function $u:\mathbb{R}^2\to \mathbb{R}$. Let $S$ be the level support function  of $\Sigma$ as defined in \eqref{eq-supportfunction}. 
 Note that $S $ is a smooth solution to \eqref{Se} with $\eta=1$ on $\{l> l_0=\inf u\}$.
By  Theorems \ref{thm-exponential1} and \ref{thm-exponential2}, 
$w(s,\theta)  :=   S(e^s ,\theta)- \sigma^{-1} e^{\sigma s} h(\theta)
$ satisfies the estimate \eqref{eq-75}:
\bea   \label{eq-75-1}
{ \Vert w(s,\cdot) \Vert_{C^0(\mathbb{S}^1)} = \Vert S(e^s ,\cdot)- \sigma^{-1} e^{\sigma s} h  \Vert_{C^0(\mathbb{S}^1)} = o(e^{(\sigma-\delta)s})\quad \hbox{ as $ s\to \infty$ }   }
\eea 
  for some  $\delta>0$.
  \medskip
  
For the solution  $h$, we let  $\Sigma_{\mathbf{0}}$ be the translator given by Theorem \ref{thm-existence} in Section \ref{subsec:existenceconstruction}, and let 
  $S_{\mathbf{0}} $ denote the level support function of  $\Sigma_{\mathbf{0}}$ given by \eqref{eq-supportfunction}.  Letting $w_{\mathbf{0}}(s,\theta):= S_{\mathbf{0}}(e^s ,\theta)- \sigma^{-1} e^{\sigma s} h(\theta) $, it follows from  Theorem \ref{thm-existence} $(i)$  that  
\bea  
w_{\mathbf{0}}(s,\cdot) = o(e^{(\sigma  -\e)s})  
\eea 
uniformly for $\theta\in\mathbb{S}^1$ as $s\to\infty $. 
Here and below in the proof,  $\e>0$  is a  sufficiently small  constant   which may vary from line to line.  
In light of  \eqref{eq-75-1}, we have  
\bea 
w - w_{\mathbf{0}} = o(e^{(\sigma  -\e)s}) 
\eea 
uniformly for $\theta\in\mathbb{S}^1$ as $s\to\infty $.
 Note that $\beta^+_{K-1}<\sigma  -\e < \beta^+_{K}$ for  a  sufficiently small $\e>0$. 
Applying   Proposition \ref{prop-coeff-extract}, there are   some $y_{K-N_1}, \ldots, y_{K-1} $ in $ \mathbb{R}$  such that 
\bea 	\label{eq-est-w-find-coeff}
w-w_{\mathbf{0}}= \sum_{i=1}^{N_1} y_{K-i} \varphi_{K-i}(\theta) e^{\beta^+_{K-i}s} + o\left(e^{(\beta^+_{K-N_1}-\e)s}\right) 
\eea 
  uniformly for $\theta\in\mathbb{S}^1$  as $s\to\infty $, where $N_1$ is the multiplicity of $\beta^+_{K-1}$. 
By using Theorem \ref{thm-existence},  we find $w_{(\mathbf{0}_{K-N_1}, y_{K-N_1}, \ldots, y_{K-1})}$, where $w_{\mathbf{y}}:=S_{\mathbf{y}}-\sigma^{-1}e^{\sigma s}h(\theta)$,  such that 
\bea 
  w_{(\mathbf{0}_{K-N_1}, y_{K-N_1}, \ldots, y_{K-1})}-w_{\mathbf{0}}= \sum_{i=1}^{N_1} y_{K-i} \varphi_{K-i}(\theta) e^{\beta^+_{K-i}s} + o\left(e^{(\beta^+_{K-N_1}-\e)s}\right)    \eea
  uniformly for $\theta\in\mathbb{S}^1$  as $s\to\infty $.  Combining with  \eqref{eq-est-w-find-coeff} yields that  
\bea 
w - w_{(\mathbf{0}_{K-N_1}, y_{K-N_1}, \ldots, y_{K-1})}= o\left(e^{(\beta^+_{K-N_1}-\e)s}\right).
\eea
As a next step,  by    applying  Proposition \ref{prop-coeff-extract}  (with $m=K-N_1-1$) we have  some $y_{K-N_1-N_2}, \ldots,  y_{K-N_1-1}$ in $ \mathbb{R}$ such that  
\bea 
w - w_{(\mathbf{0}_{K-N_1}, y_{K-N_1}, \ldots, y_{K-1})}= \sum_{i=N_1+1}^{N_1+N_2} y_{K-i} \varphi_{K-i}(\theta)e^{\beta^+_{K-i}s} +o\left(e^{(\beta^+_{K-N_1-N_2}-\e)s}\right),
\eea
where $N_2$ is the multiplicity of $\beta^+_{K-N_1-1}$.  Then  by a similar argument as above, we have
\bea 
	w - w_{(\mathbf{0}_{K-N_1-N_2}, y_{K-N_1-N_2}, \ldots, y_{K-1})}= o\left(e^{(\beta^+_{K-N_1-N_2}-\e)s}\right).
	 \eea

Repeating this argument finitely many times, we obtain $\mathbf{y}=(y_0,y_1,\ldots, y_{K-1})\in \mathbb{R}^K$ such that 
\bea 
w-w_{\mathbf{y}}=o(e^{(\beta^+_0-\e)s})=o(e^{(-2\alpha -\e)s}).
\eea
In terms of the support functions, this translates into 
\bea 
S-S_{\mathbf{y}} = o(l^{-2\alpha-\e}) 
\eea 
 uniformly for $\theta\in\mathbb{S}^1$  as $l\to\infty $ 
with  some $\e>0$.
Finally, we conclude that $\Sigma\equiv\Sigma_{\mathbf{y}}$ by \cite[Lemma 5.2]{CCK_existence}, an argument based on the comparison principle. This   completes the proof. 
\end{proof}

We remark  that, by the same argument, corresponding result of Theorem \ref{thm-72} can be obtained for blow-down translators: See \cite[Remark 5.3]{CCK_existence} for the construction of families of blow-down translators. 

\medskip

In the next proposition, using Proposition \ref{prop-coeff-extract} again, we obtain an asymptotics of blow-down translator. The main difference is that it shows the decay of a higher regularity norm of error as shown in \eqref{eq-expS}. This be used in the next Section \ref{sec-rigidityhomogeneous}. 

\begin{proposition} \label{prop-678}
 Let $\Sigma$ have unique limit shape and the graphical representation of $\Sigma$ satisfies the blow-down translator equation in the Alexandrov sense. Let  $S(l,\theta)$  be the representation of $\Sigma$  given by \eqref{eq-supportfunction}.    Then one of the following asymptotic behavior holds:   
\begin{enumerate}[(a)]
\item Either $\Sigma$ is exactly a translation of a  homogeneous profile  in $\mathbb{R}^3$, and hence it holds that  
\be \label{eq-567}
S(l,\theta) = \sigma^{-1} (l+L_3) ^\sigma h(\theta) + L_1\cos\theta + L_2 \sin\theta \ee 
for   some constants  $L_1$, $L_2$, and $L_3$ $\in \mathbb{R}$, and a solution $h $ to \eqref{eq-shrinker};

\item  Or the difference from the homogeneous solution is dominated by one of Jacobi fields of higher growth rate as follows:
\be \label{eq-expS}
 S(l,\theta) = \sigma^{-1} l^\sigma h(\theta) +  l^{ \beta^+_m}g(\theta) + o_2(l^{ \beta^+_m -\e }) 
  \ee 
   uniformly for  $ \theta\in \mathbb{S}^1$ as $l\to\infty $, 
 with  some  $m \in \{3,\ldots, K-1\}$  and     a nonzero  eigenfunction $g $ solving $L g = \lambda_m g$ in Lemma \ref{lemma-lambdavarphi}. Here, $f(l,\theta)=o_2(l^n)$ denotes that  
 \be
  \left\Vert |f(l,\cdot)|+|f_\theta(l,\cdot)|+|f_{\theta\theta}(l,\cdot)|+l | f_l(l,\cdot)|+ l^2 |f_{ll}(l,\cdot)|\right\Vert_{L^\infty(\mathbb{S}^1)}=o(l^n)\quad\hbox{  as $l \to \infty$.}
  \ee

\end{enumerate}

\begin{proof}
	Let $w (s,\theta) = S(e^s, \theta)- \sigma^{-1} e^{\sigma s} h(\theta)$.   Applying  Proposition \ref{prop-coeff-extract}   with $w_1   = w$, $w_2(s,\theta) \equiv 0$ and $m=K-1$, we obtain that 
\be\label{eq-w_1-bd-pre}
w(s,\theta)= \sum_{i=1}^N C_i \varphi_{K-i} e^{\beta^+_{K-1}s}+ o(e^{(\beta^+_{K-1} -\e )s} ) 
\ee
uniformly for  $ \theta\in \mathbb{S}^1$ as $s\to\infty $  with  some $N\in \mathbb{N}$ and $C_i \in \mathbb{R}$ ($i=1,\cdots, N$).  Here and below, $\e>0$ is some small   constant which may  vary from lime to line.  
Then  we have two cases: (i) some  $C_i$ is nonzero;  (ii)  all $C_i$ are zero. 
 \smallskip
 
(i) Suppose that some  $C_i$ is nonzero and we fix $m=K-1$ and  $g=\sum_{i=1}^N C_i \varphi_{K-i}$. Then, we will show that    
   the assertion \eqref{eq-expS}  holds by improving the error $o(l^{\beta^+_m -\e })$ in \eqref{eq-w_1-bd-pre} to $o_2(l^{\beta^+_m-\e})$. Remembering \eqref{eq-w_1-bd-pre}, we apply \cite[Lemma B.1 (interior regularity)]{CCK_existence} to   $w_1   = w $ and $w_2  = 0$ so that given $\beta\in (0,1)$ we have
\be\label{eq-162}
\Vert w \Vert_{  C^{2,\beta } _{s,1}} = O(e^{\beta^+_ms })\quad \hbox{ as  $s\to\infty $. }
\ee  
We remind that blow-down translators have $E(0)=0$  since $\eta=0$ in Lemma \ref{lem-w-eq}. Applying \cite[Lemma 3.3]{CCK_existence} (error estimate) using  \eqref{eq-162}, we obtain 
  \be \label{eq-E67}
 \Vert
  E(w)\Vert_{C^{0,\beta}_{s,1}}  =O(e^{\gamma s }) \quad \hbox{ as  $s\to\infty $}
  \ee  
with $\gamma :=   2\beta^+_{m}-\sigma <\beta^+_{m}$. 
 Since    $\beta_m^{+} $ solves   $ \beta^2+ \beta + \lambda_m =0 $,   one can see that   $e^{\beta^+_{m}s}g=\sum_{i=1}^N C_i \varphi_{K-i}  e^{\beta^+_{m}s}$ belongs to   the kernel of the operator $\hat Lw = w_{ss}+w_s+{h}^{\frac1\alpha}( w_{\theta\theta}+w) $  and hence   Lemma \ref{lem-w-eq} implies
\be 
\hat L\left(w-e^{\beta^+_ms}g \right)=E(w).
\ee
 Thus utilizing  \eqref{eq-w_1-bd-pre}, \eqref{eq-E67} and the Schauder estimate implies that
\be
\left\|w-e^{\beta^+_ms}g \right\| _{C^{2,\beta}_{s,1}}=  o\left(e^{(\beta^+_m-\e)s}\right) \quad\hbox{as  $s\to\infty $ }
\ee  
for some small $\e>0$, 
 which yields the assertion \eqref{eq-expS}. Here, we have used      $\gamma  <\beta^+_m$.

\medskip

(ii) In the case when   all $C_i$ are zero, we   apply Proposition \ref{prop-coeff-extract}  to $w_1   = w $ and $w_2  = 0$ with $m=K-N-1$ and repeat the argument until we find a   nonzero  dominating Jacobi field with the growth rate $\beta_m^+$  for  some  $m \in \{3,\ldots, K-1\}$.  Then we have      
\be  \label{eq-678} 
 S(l,\theta)= \sigma^{-1} l^\sigma h(\theta) + l^ {\beta^+_m} \sum_{i=1}^{N_m} C_i \varphi_{m-i+1}(\theta)+ o(l^{\beta^+_m-\e}) 
  \ee  
  uniformly for $\theta\in \mathbb{S}^1$ as $l\to\infty$, with  some $m \in \{3,\ldots, K-1\}$ and   some  nonzero $C_i$, where $N_m$ is 
  the multiplicity of the eigenvalue $\beta^+_m$ in the spectrum \eqref{eq-basis}. Thus    letting $g=\sum_{i=1}^{N_m} C_i \varphi_{m-i+1}$,  \eqref{eq-678} implies   the assertion \eqref{eq-expS}    with    $o(l^{\beta^+_m -\e })$  replacing $o_2(l^{\beta^+_m-\e})$.  
By repeating the same argument   as for the case (i),  we can improve this to $o_2(l^{\beta^+_m-\e})$ and obtain \eqref{eq-expS}.

\medskip

Lastly, suppose 
that  we successively fail to identify such a nonzero Jacobi field with the growth rate $\beta_m^+$  for     all $m \in \{3,\ldots, K-1\}$. Then  we deduce that   
 \be\label{eq-reach}
  S(l,\theta) = \sigma^{-1} l ^\sigma h(\theta) + o(l^{\beta^+_3-\e }) 
 \ee 
   uniformly for  $ \theta\in \mathbb{S}^1$ as $l\to\infty $.  
 We will show that a solution to \eqref{Se} (with $\eta=0$) satisfying \eqref{eq-reach} has to be a translation of the homogeneous profile.  In fact,    applying  Proposition \ref{prop-coeff-extract}   with  $m=2$ implies that  
\be
S(l,\theta) = \sigma^{-1} l^{\sigma} h(\theta) +L_1 \cos\theta +L_2 \sin\theta + o(l^{-\e}) 
\ee
   uniformly for  $ \theta\in \mathbb{S}^1$ as $l\to\infty $ with  some constants $L_1$ and $L_2$ since 
   $\beta_1^+=\beta_2^+=0$.  
Let  
\be
\hat S (l,\theta) =\sigma^{-1} l^{\sigma} h(\theta) +L_1 \cos\theta +L_2 \sin\theta.
\ee   Then 
$\hat S$ solves the same equation   \eqref{Se} (with $\eta=0$). Letting  $\hat w(s,\theta)= \hat S (e^s,\theta)-\sigma^{-1} e^{\sigma s} h(\theta) $, we     
  apply Proposition \ref{prop-coeff-extract} to $w_1=w$ and $w_2=\hat w$
  with $m=0$ to get
\be
S(l,\theta)=  \hat S(l,\theta) + L_3 l^{\sigma -1} h(\theta) + o (l^{\sigma -1 -\e}) 
\ee
    uniformly for  $ \theta\in \mathbb{S}^1$ as $l\to\infty $   with  some constant $L_3 $.
   Here, we note $\beta^+_0 =-2\alpha=\sigma-1 $. Thus we 
    deduce  that 
 \be 
  S(l,\theta)= \sigma^{-1} (l+L_3)^\sigma h(\theta) +L_1 \cos\theta + L_2 \sin \theta+ o(l^{\sigma-1-\e})\ee 
      uniformly for  $ \theta\in \mathbb{S}^1$ as $l\to\infty $.  
  Therefore,  \eqref{eq-567} follows  by \cite[Lemma 5.2]{CCK_existence}  and the solution falls into the case (a).   
  \end{proof}

\end{proposition}

 \subsection{Rigidity of entire solutions to homogeneous equation} \label{sec-rigidityhomogeneous}

In the subsection, we prove  a rigidity  result of entire solutions to the homogeneous equation  \eqref{eq-asymptranslatorgraph} under the assumption that  solutions have unique limit shapes (see Definition \ref{def-41}). This result will be used to prove the unique blow-down  of translators in Theorem \ref{thm:tangent_flow-shrinker} based on a compactness argument. 
 
Before stating the result, let us recall  from Remark \ref{rmk-homo-J} that  
\be
J^*=(2\alpha)^{-4\alpha} (h^2 +h'^2 + 4\alpha ^2 h^{2-\frac{1}{\alpha}})
\ee  is a constant  for     a  $ {\sigma}^{-1}$-homogeneous  solution  to \eqref{eq-asymptranslatorgraph}.  
By a classification result in Corollary \ref{cor-anderews2003}, we   have       the $k$-fold symmetric homogeneous solution  $\bar u_k$       to \eqref{eq-asymptranslatorgraph}      for   $3\geq k\in \mathbb{N}$,  and  the 
rotationally symmetric solution  $\bar u_\infty$.   For simplicity, we denote  $J^*_k=J_{\bar u_k}$ and $J^*_\infty=J^*_{\bar u_\infty}$.

\begin{theorem}\label{thm-rigiditywithcondition}
Let $\Sigma=\p\{x_3 > \bar u (x)\}\subset \mathbb{R}^3$ be a convex surface   having unique limit shape, where $\bar u \in C^0(\mathbb{R}^2) \cap C^{\infty}(\mathbb{R}^2\setminus \{0\} ) $  is an entire Alexandrov solution to   \eqref{eq-asymptranslatorgraph}.
  Suppose that the blow-downs of $\Sigma$ are $k$-fold homogeneous solutions \rm($k\neq \infty$)\em, and either $J^*_{\bar u}(x)\leq J_k^*$ or $J^*_{\bar u}(x)\geq J_k^*$ holds in $\mathbb{R}^2\setminus B_R(0)$ for some $R>1$.    
 Then, $\Sigma$ is a translation of $k$-fold homogeneous solution in $\mathbb{R}^3$.

\begin{proof}
Suppose to the contrary that $\Sigma$ is not a translation of $k$-fold homogeneous solution. Let us denote  by $S(l,\theta)$ the level set support function of $\Sigma$. We recall Lemma \ref{prop-J^*-supp} that
 \begin{equation}\label{eq-J^*-def}
 J^*  =   \left[  S_l    (-S_{ll})^{-1} \right]^{4\alpha }  [    S_l^2+S_{\theta l}^2   - S_{ll} (S_{\theta\theta}+S)].
\end{equation}
Using \eqref{eq-expS} in Proposition \ref{prop-678},   it follows that    
  \bea \label{eq-expSderiv}
S_l&= l^{\sigma -1} h +   \beta l ^{\beta -1 }g +o(l^{\beta -1 }),\\
S_{ll}&=(\sigma-1) l^{\sigma -2} h +  \beta (\beta-1) l^{\beta-2 }g  +o(l^{\beta-2 }),\\
S_{\theta l}&= l^{\sigma-1} h' +  \beta l^{\beta -1 }g'  +o (l^{\beta-1 }),\\
S_{\theta\theta}+S&=\sigma ^{-1}l^\sigma  (h''+h) +l^{\beta  }(g''+g) +o(l^{\beta  }),
\eea 
uniformly for $\theta\in\mathbb{S}^1$ as $l\to\infty$, where $\beta = \beta^+_m$ for some  $m \in \{3,\ldots, K-1\}$ and $g=g(\theta)$ is a nonzero eigenfunction  of  the operator $L$ in \eqref{eq-linearopL} such that $L g=\lambda_m g$ with $\lambda_m=-\beta_m^+(\beta_m^++1)$. Remembering \eqref{eq-shrinker} that $h''+h=(1-\sigma) \sigma h^{1-\frac{1}{\alpha}}$, we reformulate $Lg=\lambda_m g$ to
\begin{equation}
g''+g= -\beta(\beta+1) h^{-\frac1\alpha} g=  \frac{\beta(\beta+1)}{\sigma (\sigma-1)} (h''+h) h^{-1} g. \end{equation}
This yields
\begin{equation}\label{eq-est-Sthth}
S_{\theta\theta}+S  =  {\sigma^{-1} (\sigma-1)^{-1}}l^\sigma (h''+h) h^{-1}  \left[(\sigma -1)  h  +\beta ( \beta+1) l^{  -(\sigma-\beta) }g\right] +o(l^{\beta })
\end{equation}
uniformly for $\theta\in\mathbb{S}^1$ as $l\to\infty$. Using \eqref{eq-expSderiv},   we observe that 
\begin{align}
   l^{-(\sigma-1)}   {S_l}  &= h+  {\beta  l^{-(\sigma-\beta)} g} +o(l^{-(\sigma-\beta)}), \\
    ( 2\alpha)^{-1}l^{-(\sigma-2)}   (-{S_{ll}} )  &= h+\frac{\beta (\beta-1)}{\sigma-1}     l^{-(\sigma-\beta)}  {g} +o(l^{-(\sigma-\beta)}) 
\end{align}
uniformly for $\theta\in\mathbb{S}^1$ as $l\to\infty$. Thus,
\begin{equation}\label{eq-4550}
   (2\alpha)^{4\alpha} l^{2(\sigma-1)} \left[  S_l    (-S_{ll})^{-1}\right]^{4\alpha }  =1- 2 \beta  (\sigma-\beta)     l^{-(\sigma-\beta)} h^{-1} { g}  +o(l^{-({\sigma-\beta})})
\end{equation} 
uniformly for $\theta\in\mathbb{S}^1$ as $l\to\infty$. 
In light of  \eqref{eq-expSderiv} and \eqref{eq-est-Sthth}, we also have 
\bea 
l^{-2(\sigma-1)}S_l^2&=h^2+2\beta l^{-(\sigma-\beta)}hg+o(l^{-(\sigma-\beta)}),\\
l^{-2(\sigma-1)}S_{\theta l}^2&=h'^2+2\beta l^{-(\sigma-\beta)}h'g'+o(l^{-(\sigma-\beta)}),\\
l^{-2(\sigma-1)}S_{ll}(S_{\theta \theta}+S) 
&=-  \sigma ^{-1}(h''+h) \left[2\alpha  h   -2 \beta ^2 l^{ -(\sigma-\beta)}g\right]+o(l^{-(\sigma-\beta)}) 
\eea 
uniformly for $\theta\in\mathbb{S}^1$ as $l\to\infty$. 
 By Remark \ref{rmk-homo-J} we have 
\be\label{eq-Jkrecall}  
J^*_k =  (2\alpha)^{-4\alpha} \left(h^2 +h'^2 + 4\alpha ^2 h^{2-\frac{1}{\alpha}}\right) =(2\alpha)^{-4\alpha} \left[h^2 +h'^2 + 2\sigma^{-1}\alpha    (h''+h)h\right],
 \ee  
and thus we can obtain
\bea\label{eq-455}
&l^{-2(\sigma  -1)}[S_l^2+S_{\theta l}^2- S_{ll}(S_{\theta \theta}+S) ]\\
&= (2\alpha)^{4\alpha}  J^*_k +2\beta l^{-(\sigma-\beta)}  [hg+h'g'- \sigma^{-1} \beta  (h''+h)g] +o(l^{-(\sigma-\beta)} )
\eea
uniformly for $\theta\in\mathbb{S}^1$ as $l\to\infty$. 
We multiply \eqref{eq-4550}  and \eqref{eq-455}, and then combine with \eqref{eq-J^*-def} to get 
\begin{equation}
  J^*-   J^*_k =  -2\beta(\sigma -\beta) l^{-(\sigma-\beta)}  h^{-1}g   J^*_k+ (2\alpha)^{-4\alpha} (2\beta) l^{-(\sigma-\beta)}    [hg+h'g'- \sigma^{-1}{\beta } (h''+h)g] +o(l^{-(\sigma-\beta)}), 
\end{equation}
and hence remembering \eqref{eq-Jkrecall} and $\sigma=1-2\alpha$ we have
\bea
  & (2\alpha)^{4\alpha}   l^{(\sigma-\beta)}h  (J^*-J^*_k) \\&=2\beta \left[ -   (\sigma-\beta) \left\{h^2+h'^2+  2\sigma^{-1}\alpha(h''+h)h\right\} g+h^2g+h'hg'- \sigma^{-1}{\beta }{  }(h''+h)hg \right]+o(1)\\
  &=2\beta  \left[ - (\sigma-\beta) h'^2 g  - (2\alpha+\beta)   h''h g + h'hg'\right]+o(1) 
\eea
 uniformly for $\theta\in\mathbb{S}^1$ as $l\to\infty$.

To prove the desired result,  it suffices to show that   the function $D=D(\theta)$, given by 
 \be 
 D:=2\beta  \left[ - (\sigma-\beta) h'^2 g  - (2\alpha+\beta)   h''h g + h'hg'\right],
 \ee
 attains both positive and negative values on $   [0,2\pi]$, because the assumption $J^*_{\bar u}\leq J^*_k$ and the equation above imply $D \leq 0$.  Since $h$ is $k$-fold symmetric,  we notice that $h' $ has $2k$ zeros, say $\theta_i=\theta_0+\frac{\pi }{k }i   $   for   $i=1,\ldots, 2k $ with  some $\theta_0$.
  On those angles $\theta_i$, we have 
 \begin{align}\label{eq-D-thei}
 D = -2\beta (2\alpha+\beta)  h''hg. 
 \end{align}
Here  we note that $\beta (2\alpha +\beta)\neq 0$ since $0<\beta = \beta^+_m<\sigma$.  
We    observe that   $h'$ solves $L h'=-\sigma(\sigma+1) h'$,  and that $h'$ and $h''$ can not be both zero at any angle as otherwise $h$ has to be a constant.   
Similarly, $g$ and $g'$ can not be both zero at any angle.  
 By \cite[Theorem 2.4]{ChoiSun}, which is a consequence of the Sturm-Liouville theory, $g$ should have strictly less number of zeros than $h'$ since $\lambda_m=-\beta(\beta+1)> -\sigma(\sigma+1)$. 
  Thus, there should be two consecutive zeros of $h'$, say $\theta_{k}$ and $\theta_{k+1}$, such that $g$ does not change its sign on $  (\theta_k,\theta_{k+1})$.  
  Firstly, suppose $g(\theta_{k})=0$ or $g(\theta_{k+1})=0$. Then,    at such angle,    it holds that   $h'=0$, $g=0$,    $D=0$ and   
  \begin{align} 
  D' =2\beta (\sigma -\beta)   h'' h g'\neq 0 ,
  \end{align}
where we note that
   $\beta(\sigma-\beta)  \neq 0$. 
   Thus $D$ changes its sign near the angle.
    Next, suppose that  $g(\theta_{k})\neq 0$ and $g(\theta_{k+1})\neq0$.   Then, $g $ is either positive or negative on $ [ \theta_k,\theta_{k+1}]$.  
   However, we have $h''(\theta_k)h''(\theta_{k+1})<0$, because $h$ attains its local maximum at one of the angles $\theta_k,\theta_{k+1}$ and attains its local minimum at the other angle. Thus, \eqref{eq-D-thei} yields
   $D(\theta_k)D(\theta_{k+1})<0$ so that we conclude that $  J^*_{\bar u}-   J^*_k $ changes its sign  for large $l$, which completes the proof by contradiction.	
\end{proof}
\end{theorem}

\section{Proofs of Theorems \ref{thm:tangent_flow-shrinker} and \ref{thm:classification-2}}\label{sec-uniq-bd1}

In this section, we prove the uniqueness of the blow-down limit and the  classification of translators. Working with the Legendre transformation, recall that the dual of the height function of a translator satisfies 
\bea  \label{eq-dualtranslator-1}
      \det D^2 u^* = (1+|x|^2)^{\frac{1}{2\alpha}-2} .
      \eea
Let us first recall  some results from   \cite{DSavin}.   Daskalopoulos and Savin  \cite{DSavin}    considered  the equation with the homogeneous right hand side:
\begin{equation}\label{homo.eq}
\det D^2 \bar v(x)=|x|^{\frac{1}{\alpha}-4},
\end{equation}
and introduced the following function $J_v$ which will play a role of entropy.
\begin{definition}\label{def:J_v}
Given a strictly convex function $v$, we define a function $J_v$ by
\begin{equation}
 J_v(y) = \Delta  v(y) \left\{y_iy_j v_{ij}(y)\right\}^{4\alpha-1}. 
\end{equation}
\end{definition}

 In fact, the entropy $J^*$ in Definition \ref{def-J*} and $J$ in Definition \ref{def:J_v} represent the same quantity. 

\begin{lemma}\label{lem-JJ*} For a striclty convex function $u$ and its   Legendre transform  $u^*$, there holds
\be
J^*_u (x) = J_{u^*}(\nabla u(x)).
\ee
\end{lemma} 
\begin{proof}
We recall the well-known identity between the Hessians of  $u$ and $u^*$:  
\be
D^2 u^* (\nabla u(x)) = [D^2 u (x)]^{-1}=\{u^{ij}(x)\},
\ee
where   the inverse matrix of $D^2 u(x)$ is denoted by $\{u^{ij}(x)\}$ as in Definition \ref{def-J*}. Then it is straightforward to calculate that 
\be
J_{u^*}(\nabla u(x))= \left\{u^{11}(x)+u^{22}(x)\right\}\left\{u_i(x)u_j(x) u^{ij}(x)\right\}^{4\alpha -1}= J_u (x) .
\ee
This finishes the proof.
\end{proof}
 \smallskip

We notice that $J$ is invariant under the homogeneous rescaling, namely for every $\lambda>0$
\begin{align}
 J_{v_\lambda}(x)=J_v(\lambda  x) \text{\;\;\;holds, \;where \;} v_\lambda(x):=\lambda^{-\frac{1}{2\alpha}}v(\lambda x).
\end{align}
  In light of Theorem \ref{thm-class-homo-dual}, given $\alpha\in (0,1/4)$, 
  we let $m$ be the largest integer satisfying $m^2 <   1/\alpha$. 
  Then, there are exactly $(m-1)$-homogeneous solutions to  \eqref{homo.eq}  up to rotations, which are the rotationally symmetric solution and  $k$-fold symmetric solutions for  $3\leq k\leq m$ (if $m\geq3$). 

For the sake of brevity,  we denote 
\be
J_k:=J_{\bar v_k} \qquad \hbox{for any   integer $3\le k < \alpha^{-1/2}$,  }
\ee 
 where $\bar v_k$ is a $k$-fold symmetric homogeneous solution to \eqref{homo.eq}   for    $3\le k < \alpha^{-1/2}$. Moreover, we denote  $J_\infty:=J_{\bar v_\infty}$, where $\bar v_\infty$ is the rotationally symmetric solution to \eqref{homo.eq}. 
 Then,  $J_{ k}$ is a constant function by Theorem \ref{thm_Entorpy order}. 

\smallskip
 From   Section \ref{sec-rigidityhomogeneous}, we  recall   that  $J^*_k:=J^*_{\bar u_k}$  and $J^*_\infty:=J^*_{\bar u_\infty}$,  where $\bar u_k$ is a $k$-fold symmetric homogeneous solution to \eqref{eq-asymptranslatorgraph} for any   integer $3\le k < \alpha^{-1/2}$ and   $\bar u_\infty$ is  the rotationally symmetric solution to \eqref{eq-asymptranslatorgraph}. Since   $\bar u_k$ and $\bar v_k$ are dual to each other (upto a rotation), it follows by Lemma \ref{lem-JJ*}  
 that  
 \be
J^*_k=J_k\qquad \hbox{and}\qquad J^*_\infty = J_\infty. 
\ee

 Next, we give some results of \cite{DSavin} regarding the entropy $J$, which will be used in this section.

\begin{proposition}\label{prop_homo.int.max}
Suppose that $\bar v\in C^\infty(\mathbb{R}^2\setminus\{0\})$ is an Alexandrov solution to \eqref{homo.eq} in $\mathbb{R}^2$, and $J_{\bar v}$ attains its local maximum or minimum at $x_0\neq 0$. Then, we have $J_{\bar v}(x_0)=J_{\infty}$ unless $J_{\bar v}$ is a constant.
\end{proposition}

\begin{proof}
The   proof of  Proposition 5.6 in \cite{DSavin} already implies  this proposition, which is slightly stronger than the original statement in \cite{DSavin}. We briefly explain the proof for the reader's convenience.
 Let $M=\log J_{\bar v}$ and $\mathcal{O}=\{x\in \mathbb{R}^2\setminus\{0\} : J_{\bar v}(x)\not=J_{\infty}\}$.
From  the proof  of   \cite[Proposition 5.6]{DSavin}, there exists a smooth function  $C(x)$ defined in $\mathcal{O}$ such that
\be
|\bar v^{ij}M_{ij}|\leq C(x)\left(|\nabla M|+|\nabla M|^2\right)
\ee
holds in $\mathcal{O}$. Therefore, the strong maximum principle completes the proof.
\end{proof}

\begin{proposition}\label{prop_homo.constnat}
 Suppose that $\bar v\in C^\infty(\mathbb{R}^2\setminus\{0\})$ is   an Alexandrov solution to \eqref{homo.eq} in $\mathbb{R}^2$,  $\bar v (0)=0$, $\nabla \bar v (0)=0$ and   $J_{\bar v}$ is a constant on $\mathbb{R}^2\setminus \{0\}$. Then, $\bar v$ is a homogeneous function. 
\end{proposition}

\begin{proof} 
See the second part of the proof of Proposition 5.7 in \cite{DSavin}, from line 27 pp. 663 to line 9 pp. 664. 
\end{proof}

 In the following, we  study subsequential limits of the blow-down.

 \begin{proposition}\label{prop_subconvergence}
  Given $\alpha \in (0,1/4)$, let  $  v\in C^\infty(\mathbb{R}^2\setminus\{0\})$ be   an Alexandrov solution to \eqref{eq-dualtranslator-1} or \eqref{homo.eq} in $\mathbb{R}^2$. 
 For every divergent sequence $\{\lambda_i\}_{i\in\mathbb{N}}\subset \mathbb{R}^+$,    rescaled functions $v_i(x):=\lambda_i^{-1/{2\alpha}}v(\lambda_i x)$ have a subsequence converging to an Alexandrov solution $\bar v$ to the equation \eqref{homo.eq} in the $C^\infty_{\text{loc}}(\mathbb{R}^2\setminus\{0\})$-topology.
 \end{proposition}
 
 \begin{proof}
By the growth rate estimates in \cite[Theorem 1.2 and Remark 1.3]{CCK_regularity}, there exist some  constants $C>1 $ and $R>0$ such that
 \begin{align}
 C^{-1}|x|^{\frac{1}{2\alpha}}\leq v(x) \leq  C|x|^{\frac{1}{2\alpha}}\qquad \hbox{for $|x|\geq R$}.
 \end{align}
Then  $v_i$ also satisfies 
  \be
   C^{-1}|x|^{\frac{1}{2\alpha}}\leq v_i(x) \leq  C|x|^{\frac{1}{2\alpha}} \qquad \hbox{for  $|x|\geq \lambda_i^{-1} R$}.
   \ee 
 Since   $v_i$  solves
 $
\det D^2 v_i= \left(\lambda_i^{-2} +|x|^2\right)^{\frac{1}{2\alpha}-2}$ or $\det D^2 v_i=  |x|^{\frac{1}{2\alpha}-4}$,
 the interior regularity theory of Monge--Amp\`ere equations yields the desired result; we refer to \cite{MA-figalli-book}.
 \end{proof}

 \begin{lemma}\label{lem_lim.J}
 Given $\alpha \in (0,1/4)$, let $m$ be the largest integer such that $m^2< 1/\alpha$, and  
 let  $  v\in C^\infty(\mathbb{R}^2\setminus\{0\})$ be    an Alexandrov solution to \eqref{eq-dualtranslator-1} or \eqref{homo.eq} in $\mathbb{R}^2$. 
  Then, there exist $k_1,k_2\in \{3,4,\cdots,m\}\cup \{\infty\}$ such that
\begin{align}\label{eq_limsup.J}
\limsup_{|x|\to +\infty}J_{ v}(x)= J_{k_1}, \qquad\quad \liminf_{|x|\to +\infty}J_{ v}(x) = J_{k_2}.
\end{align}
If a sequence $\{x_i\}_{i\in \mathbb{N}}$ satisfies that $|x_i|\to +\infty$ and 
\be
 J_{ v}(x_i)\,\,\to\,\, J_{k_1}:=\limsup J_v \quad  \hbox{ as $i\to +\infty$,}
\ee
 then    rescaled functions $v_i(x):=|x_i|^{-1/{2\alpha}}   v(|x_i| x)$ have a subsequence converging to a $k_1$-fold symmetric homogeneous solution to \eqref{homo.eq} in the $C^\infty_{\text{loc}}(\mathbb{R}^2\setminus\{0\})$-topology. 
 \smallskip
 
 Moreover, if $k_1=\infty$, then $v$ has the unique blow-down $\,\bar v_{\infty}(x)=\displaystyle\lim_{\lambda\to +\infty}\lambda^{-1/{2\alpha}}v(\lambda x)$, which is the rotationally symmetric homogeneous solution to \eqref{homo.eq}.
\end{lemma}  
 
 \begin{proof} 
 Let $\{x_i\} $ be a sequence such that $|x_i |\to+ \infty$ and $\displaystyle J_v(x_i) \to \limsup  J_v$ as $i\to+ \infty$. By Proposition \ref{prop_subconvergence}, the sequence of  $ v_{i}(x):= |x_i|^{- {1}/{2\alpha}}  v(|x_i| x)$ sub-converges to a solution $\bar v$ to \eqref{homo.eq} in the $C^\infty_{\text{loc}}(\mathbb{R}^2\setminus\{0\})$-topology.   We may further assume that $x_i/|x_i|$ sub-converges to  $x_0 \in \mathbb{S}^1:= \{|x|=1\}$. Then, $J_{\bar v}$ attains its global maximum $J_{\bar v}(x_0)=\limsup J_v<+\infty$  at  $x_0 \in \mathbb{S}^1$.   
 Here,   we note that $J_{v_i}(x)= J_v(|x_i| x)$  due to the invariance of $J$ under the homogeneous rescaling.

\medskip
 
 First, let us assume that  $\limsup J_v =J_{\bar v}(x_0)\neq J_{\infty}$.   Then,    Proposition \ref{prop_homo.int.max} implies that  $J_{\bar v}$ is a constant since $J_{\bar v}$ has its global maximum $J_{\bar v}(x_0)$.  
Thus,    $\bar v$ is a homogeneous solution to \eqref{homo.eq} by Proposition \ref{prop_homo.constnat},   and hence there exists some $k_1 \in \{3,4,\cdots,m\}$ such that $\limsup J_v=J_{k_1}$   by  Theorem \ref{thm-class-homo-dual}.  Furthermore, we deduce that     $\{v_i\}$ sub-converges to a $k_1$-fold symmetric homogeneous solution to \eqref{homo.eq}  in the $C^\infty_{\text{loc}}(\mathbb{R}^2\setminus\{0\})$-topology. In the same manner,   there exists some $k_2 \in \{3,4,\cdots,m\}\cup \{\infty\}$ such that $\liminf J_v=J_{k_2}$. 

\medskip
Next,   suppose that   $\limsup J_v = J_{\infty}$. Then, by Theorem \ref{thm_Entorpy order}, we have  that 
\be
\limsup J_v=J_{\infty}\leq J_{k_2}=\liminf J_v,
\ee
which implies that  
\be
\lim_{|x|\to +\infty} J_v=J_{\infty}.
\ee
Using  Proposition \ref{prop_subconvergence}, this yields that for a given divergent sequence $\{\lambda_i\} $,  a sequence of  $\lambda_i^{-\frac{1}{2\alpha}}v(\lambda_i x)$ sub-converges to a solution $\bar v$ to  \eqref{homo.eq} and  $J_{\bar v}(x)= J_\infty$ for  all $x\in \mathbb{R}^2\setminus \{0\}$. By Proposition \ref{prop_homo.constnat}, $\bar v$ is a rotationally symmetric homogeneous solution.  Therefore, $v$ has the unique blow-down $\bar v$ since  the subsequential limit $\bar v$ is unique.
 \end{proof}
 
 To prove    Theorem \ref{thm:tangent_flow-shrinker}, we   first show that     the unique blow-down can be obtained  by the existence of the entropy limit at infinity.

 \begin{proposition}\label{prop_support.convergence}
 Given $\alpha \in (0,1/4)$, let $m$ be the largest integer such that $m^2< 1/\alpha$.  Suppose that   $  v\in C^\infty(\mathbb{R}^2\setminus\{0\})$ is   an Alexandrov solution to \eqref{eq-dualtranslator-1} or \eqref{homo.eq} in $\mathbb{R}^2$, and $\lim_{|x|\to +\infty} J_v(x)=J_k$ for some $k\in \{3,4,\cdots,m\}\cup \{\infty\}$. Let $S(l,\theta)$ denote the level support function of  the Legendre transform $v^*$ of  $v$. 
 Then, the rescaled functions $\lambda^{-\sigma}S(\lambda l, \theta)$ with $\sigma:=1-2\alpha$ converge  to $\sigma^{-1} l^{\sigma}h(\theta)$ in the $C^{\infty}_{\text{loc}}((l_0,\infty) \times \mathbb{S}^1)\cap C^{0}_{\text{loc}}([l_0,\infty)\times \mathbb{S}^1) $-topology as $\lambda\to +\infty$. Here,     $h\in C^{\infty}(\mathbb{S}^1)$ is a solution to \eqref{eq-shrinker}, and $l_0:=\min_{\mathbb R^2} v^*$.
 \end{proposition}
 
 \begin{proof}
To begin with, let us  fix a $k$-fold symmetric homogeneous solution $\bar v_k$ to \eqref{homo.eq}. In light of Remark \ref{remark-shrinker}, there exists  a solution  $h\in C^{\infty}(\mathbb{S}^1)$  to \eqref{eq-shrinker} such that $\sigma^{-1}l^{\sigma}h(\theta)$ is the support function of the level curve $\{{(\bar v_k)^*=l}\}$. Here, $(\bar v_k)^*$ denotes    the Legendre transform $\bar v_k$, which  solves \eqref{eq-asymptranslatorgraph}.   

\medskip
 
Next, by the growth rate estimates in \cite[Theorem 1.2 and Remark 1.3, and Corollary 1.4]{CCK_regularity},   $v$ and $v^*$ satisfy
 \begin{align}\label{eq_growth.rate.dual}
& C^{-1}|x|^{\frac{1}{\sigma}}\leq v^*(x) \leq  C|x|^{\frac{1}{\sigma}}, &&  C^{-1}|x|^{\frac{1}{2\alpha}}\leq v(x) \leq  C|x|^{\frac{1}{2\alpha}},
 \end{align}
in $\{|x| \geq R\}$ for some large  constants $C $ and $R$. Note that  $v^*$ is an entire function. Given $\lambda>0$, we define
 \begin{align}
& v^*_\lambda(x):=\lambda^{-\frac{1}{\sigma}} v^*(\lambda x), && v_\lambda(x):=\lambda^{-\frac{1}{2\alpha}} v(\lambda x).
 \end{align}
We observe that $v^*_\lambda$ and $v_\lambda$ satisfy \eqref{eq_growth.rate.dual} in $\{|x|\geq \lambda^{-1}R \}$. We denote by $(v_\lambda)^*$ the dual function of $v_\lambda$, and we should notice that $v^*_\lambda \neq (v_\lambda)^*$. In fact, we claim that 
 \begin{equation}\label{eq-ddual}
(v^*_\lambda)^* =((v_{\lambda^{  2\alpha/\sigma}})^*)^*=v_{\lambda^{  2\alpha/\sigma}}.
\end{equation}
To prove the claim, we manipulate 
\begin{equation}
(v_\lambda)^*(\nabla v_\lambda(x))=  \nabla v_\lambda (x)\cdot x-v_\lambda (x)=\lambda^{-\frac{1}{2\alpha}}[ \nabla v(\lambda x)\cdot \lambda x- v(\lambda x)]=\lambda^{-\frac{1}{2\alpha}} v^*(\nabla v(\lambda x))
\end{equation}
 by using the definition of the Legendre transform.  
Since  $\nabla v_\lambda(x)=\lambda^{1-\frac{1}{2\alpha}} \nabla v (\lambda x)=\lambda^{-\frac{\sigma}{2\alpha}}\nabla v(\lambda x)$, replacing $\nabla v_\lambda(x)$ by $p$ yields that
\begin{align}
(v_\lambda)^*(p)=\lambda^{-\frac{1}{2\alpha}}v^*(\lambda^{\frac{\sigma}{2\alpha} }p)=v^*_{\lambda^{{\sigma}/(2\alpha)}}(p).
\end{align}
This implies \eqref{eq-ddual}.

\medskip
Now, let  $\{\lambda_i\} $ be   any    divergent sequence.  By using Lemma \ref{lem_lim.J} and the assumption $\lim J_v=J_k$, we can find a subsequence $\{\lambda_{i_m}\}$ and a rotation $A\in SO(2)$ such that the sequence $\{v_{\lambda_{i_m}^{ 2\alpha/\sigma}}\}$ converges in $C^{\infty}_{\text{loc}}(\mathbb{R}^2\setminus\{0\})$ to $\bar v_k \circ A$; see Theorem \ref{thm-class-homo-dual}. 
  Then, in light of \eqref{eq-ddual},   one can show that  the dual sequence $\{v^*_{\lambda_{i_m}}\}$ converges at least in the $C^0_{\text{loc}}(\mathbb{R}^2)$-topology to $(\bar v_k\circ A)^*$,    the Legendre  dual of $\bar v_k\circ A$.   Note that $(\bar v_k\circ A)^*$ is a homogeneous solution to  \eqref{eq-asymptranslatorgraph}.  
Moreover, since   \eqref{eq_growth.rate.dual} implies 
 \be
  C^{-1}|x|^{\frac{1}{\sigma}}\leq  v^*_\lambda \leq C|x|^{\frac{1}{\sigma}}
  \ee   in $\{|x|\geq \lambda^{-1}R\}$ for the   constants  $C $ and $R$ (independent of $\lambda$), the   interior regularity theory for Monge--Amp\`ere equations guarantees the $C^{\infty}_{\text{loc}}(\mathbb{R}^2\setminus\{0\})$-convergence of $\{v^*_{\lambda_{i_m}}\}$ to  $(\bar v_k\circ A)^*$.

\smallskip

On the other hand, we see that  the level support function of $v^*_{\lambda_{i_m}}$ is given by $\lambda_{i_m}^{-1}S(\lambda_{i_m}^{\frac{1}{\sigma}} l,\theta)$. In addition, there exists some $\theta_A\in \mathbb{S}^1$ such that $\sigma^{-1}l^{\sigma}h(\theta+\theta_A)$ is the level support function of the limit $(\bar v_k\circ A)^*$ by Corollary \ref{cor-anderews2003}. 
 Thus, any divergent sequence $\{\lambda_i\}$ has a subsequence $\{\lambda_{i_m}\}$ such that $\lambda_{i_m}^{-1}S(\lambda_{i_m}^{\frac{1}{\sigma}} l,\theta)$ converges to $\sigma^{-1}l^{\sigma}h(\theta+\theta_A)$ for some $\theta_A\in \mathbb{S}^1$ in the $C^{\infty}_{\text{loc}}((l_0,\infty)\times \mathbb{S}^1)\cap C^{0}_{\text{loc}}([l_0,\infty)\times \mathbb{S}^1)$-topology.
This shows that  $S(l,\theta)$ has unique limit shape  as in Definition \ref{def-41}. 
By  Theorem \ref{thm-42}, $S(l,\theta)$ has the unique blow-down $\sigma^{-1}l^{\sigma}h(\theta+\theta_A)$ for some $\theta_A$ independent of blow-down sequences. We replace $h(\theta+\theta_A)$ by $h(\theta)$ so that we complete the proof.
 \end{proof}

\medskip

 Now, we are ready to prove Theorem \ref{thm:tangent_flow-shrinker}. Here is a brief idea. By Lemma \ref{lem_lim.J} and Proposition \ref{prop_support.convergence}, the theorem fails only if $\limsup J_v\neq \liminf J_v$. Suppose $\limsup J_v = J_{k_1}$ with $k_1\neq \infty$. By a point picking argument, we obtain a solution $\hat v$   to the blow-down equation \eqref{homo.eq}, which is an asymptotically $k_1$-fold solution   with a  non-constant function $J_{\hat v} \leq J_{k_1}$. This contradicts   Theorem \ref{thm-rigiditywithcondition}, where we showed such  a solution  $\hat v$ should be $k_1$-fold symmetric homogeneous.  
 
\begin{proof}[Proof of Theorem \ref{thm:tangent_flow-shrinker}]
By Theorem \ref{thm-smooth}, $\Sigma$ is the graph of a smooth entire function $u\in C^{\infty}(\mathbb{R}^2)$ satisfying
\begin{equation}
 C^{-1}|x|^{\frac{1}{\sigma}}\leq u(x) \leq  C|x|^{\frac{1}{\sigma}},
\end{equation}
on $\{|x| \geq R\}$ for some large constants  $C$ and $R$. Hence, its Legendre transformation $v:=u^*$ is also a smooth entire function. We claim that $\lim_{|x|\to +\infty}J_v(x)$ exists. If the claim is true, the result follows from Proposition \ref{prop_support.convergence}. 

\medskip

Since $v$ satisfies \eqref{eq-dualtranslator-1}, there exists $k_1\in \{3,4,\cdots\}\cup \{\infty\}$ satisfying 
\be
\limsup_{|x|\to +\infty}J_{ v}(x)= J_{k_1}.
\ee 
Towards a contradiction,  suppose that $\liminf J_v \neq \limsup J_v$. Then,  we have $k_1\neq \infty$ by Lemma \ref{lem_lim.J}, and hence there exists $\delta>0$ such that
\begin{equation}\label{eq-est-J-1}
J_{k_1}-3\delta \geq J_{k_1+1} \ge \liminf_{|x|\to +\infty}J_v(x),
\end{equation} 
where $J_{k_1+1}=J_{\infty}$ if the $(k_1+1)$-fold shrinker does not exist for chosen $\alpha$. Here, the second inequality  follows from Lemma \ref{lem_lim.J}. 

 Let   $\{x_i\} $ be a   sequence  such that $|x_i|\to +\infty$ and $J_v(x_i)\to J_{k_1}$ as $i\to +\infty$. Then, by Lemma \ref{lem_lim.J}, the rescaled functions $v_i(x):=|x_i|^{-\frac{1}{2\alpha}}v(|x_i|x)$    have a   subsequence, still denoted by  $\{v_{i}\}$, converging to a $k_1$-fold symmetric homogeneous solution to \eqref{homo.eq}. Hence, there exists  $N> 1$ such that
\begin{equation}\label{eq-est-J-2}
\min_{|x|=1}J_{v_{i}}(x)=\min_{|x|=|x_{i}|}J_{v}(x)\geq J_{k_1}-\delta\qquad\hbox{for all $i\geq N$}
\end{equation}
 due to the invariance of $J$ under the homogeneous rescaling. 
Using  \eqref{eq-est-J-1} and \eqref{eq-est-J-2},    we find some $N_1\geq N$ such that for each $i\geq N$, there is a point $y_{i}$   satisfying  
 \be\label{eq-yixi}
|y_{i}|<|x_{i}|, \qquad  J_{v}(y_{i})=J_{k_1}-2\delta, \qquad J_{v}(x) \geq J_{k_1}-2\delta \quad\hbox{on  $\{ |y_{i}|\leq |x| \leq |x_{i}|\}$}. 
 \ee
 Such a sequence of point $y_i$ necessarily diverges, namely $|y_i|\to\infty$ as $i\to \infty$, as otherwise we find $\liminf J_v \ge J_{k_1}-2\delta$, a contradiction to \eqref{eq-est-J-1}. Moreover, it holds that  $\lim_{i\to \infty} \tfrac{|y_i|}{|x_i|}=0$ due to the  $C^{\infty}_{loc}(\mathbb{R}^2 \setminus \{0\})$-sub-convergence of arbitrary subsequences of $v_i$ to a $k_1$-fold symmetric homogeneous function whose $J$ value is constant $J_{k_1}$.

\medskip

Now, let us consider the rescaling of $v$ by $|y_i|$, $\hat v_i (x):= |y_i|^{-\frac1 {2\alpha}} v(|y_i|x)$. By Proposition \ref{prop_subconvergence}, $\{\hat v_{i}\}$ sub-converges to an Alexandrov\ solution $\hat  v\in C^\infty(\mathbb{R}^2\setminus\{0\})\cap C^0(\mathbb{R}^2)$ to \eqref{homo.eq}. In addition, $\hat  v$ satisfies
\begin{align}\label{eq_J.non.constant}
&\inf_{|x|\geq 1 }J_{\hat  v}(x)=\min_{|x|=1}J_{\hat  v}(x)=J_{k_1}-2\delta, && \sup_{|x|\neq 0}J_{\hat  v}(x)\leq J_{k_1}.
\end{align}
Here, the first assertion is a consequence of $|y_i|<|x_i|$ and $|x_i|/|y_i| \to \infty$ as $i \to \infty$. 
Since $J_{k_1}-2\delta>J_{k_1+1}$, Lemma \ref{lem_lim.J} implies 
\be
\lim_{|x|\to +\infty}  J_{\hat  v}(x)=J_{k_1}.
\ee
Hence, by Proposition \ref{prop_support.convergence}, the Legendre dual   $\hat  u:=\hat  v^*$ of $\hat  v$ has the unique blow-down, which is a $k_1$-fold symmetric homogeneous solution to \eqref{eq-asymptranslatorgraph}.
 Moreover, in view of Lemma \ref{lem-JJ*}, the second assertion in \eqref{eq_J.non.constant} shows  that 
  \be
  \sup_{|x|\neq 0} J^*_{\hat  u}(x) =  \sup_{|x|\neq 0}J_{\hat  v}(x)\leq J_{k_1} = J^*_{k_1} .
 \ee
 By Theorem \ref{thm-rigiditywithcondition},  $\hat  u$ is a homogeneous function, and thus $\hat  v=\hat  u^*$ is also a homogeneous function. This contradicts \eqref{eq_J.non.constant}.  	Therefore, we have proved that $\lim_{|x|\to+\infty} J_v(x)$ exists, which finishes the proof.
\end{proof}
 
Lastly, we prove the classification of translators. 
\begin{proof}[Proof of Theorem \ref{thm:classification-2}]
By Theorem \ref{thm:tangent_flow-shrinker}, every given translator is asymptotic to a homogeneous solution to \eqref{eq-asymptranslatorgraph}.  Note that  the conclusion of Theorem \ref{thm:tangent_flow-shrinker} implies the unique limit shape assumption (see Definition \ref{def-41}). Then, by Theorem \ref{thm-72}, the translator has to be one of solutions constructed in Theorem \ref{thm-existence} in Section \ref{subsec:existenceconstruction}. 
\end{proof}

\bigskip

\section{Proof of Theorem \ref{thm:moduli.space.topology}} \label{sec-uniq-bd2}
Given $\alpha\in (0,\frac{1}{4})$, we denote by $\mathcal{S}_k$ the space of $k$-fold symmetric solutions to the shrinker equation \eqref{eq-shrinker}. Note that $\mathcal{S}_{\infty}$ is a point, and $\mathcal{S}_k$ with $k\neq \infty$ is diffeomorphic to $\mathbb{S}^1$. 

Given $h\in \mathcal{S}_k$, by  Theorem \ref{thm-existence} in Section \ref{subsec:existenceconstruction}, there exists a $K$-parameter family of translators $\{\Sigma_{\mathbf{y}}^h\}_{\mathbf{y}\in \mathbb{R}^K}$ whose level support functions are asymptotic to $\sigma^{-1}l^\sigma h(\theta)$ as $l\to +\infty$, where $K$ is the dimension of slowly decaying Jacobi fields. 
In light of the classification result of Theorem \ref{thm-72}, we can define a map $\psi:\mathcal{X}_k\to \mathbb{R}^K\times \mathcal{S}_k$   by $\psi(\Sigma_{\mathbf{y}}^h)=(\mathbf{y},h)$. Note that we consider not only the parameter $\mathbf{y}$ but also the asymptotic $h$ in this subsection because the asymptotics of translators in $\mathcal{X}_k$ are not necessarily the same.  Moreover, let us assume that we corotationally define the eigenfunctions  $\{(\varphi_j^h,\beta_j^\pm \varphi_j^h)\}_{j=0}^\infty$ of $\mathcal{L}$ with respect to $h$, i.e., rotation $\mathcal{R}_{\theta}$ becomes an equivariant action: $(\varphi_j^{ \mathcal{R}_\theta h},\beta_j^\pm \varphi_j^{ \mathcal{R}_\theta h}) =(\mathcal{R}_\theta \varphi_j^{  h},\beta_j^\pm \mathcal{R}_\theta \varphi_j^{  h}) $. Similarly, $\Sigma^h_{\mathbf{y}}$ are corotationally defined so that $\mathcal{R}_{\theta} \Sigma^h_{\mathbf{y}}= \Sigma ^{\mathcal{R}_{\theta} h}_{\mathbf{y}}$.  

To discuss the topology of moduli spaces $\mathcal{X}^k$, we impose the topology induced by the locally uniform convergence on surfaces. See \cite[pp 4]{CCK_existence} for a discussion on the necessity of this choice of topology. We say that a sequence of translators $\Sigma_i$  converges to a complete convex surface $\Sigma$ in $C^0_{\text{loc}}$ (or $C^0_{\text{loc}}(\mathbb{R}^3)$ for clarification) if given a compact set $K\subset \mathbb{R}^3$, the sequence $\Sigma_i\cap K$ converges to $\Sigma \cap K$ in the Hausdorff distance. Note that the $C^0_{\text{loc}}$-convergence improves to the $C^{\infty}_{\text{loc}}$-convergence, and the limit surface $\Sigma $ is also a translator. Indeed, the Legendre duals of the height function of $\Sigma_i$ converge to the Legendre dual of the height function of $\Sigma $ in $C^{0}_{\text{loc}}$, and  hence by a standard regularity theory (see \cite[Corollary 2.12]{MA-figalli-book}),   the limit is an Alexandrov solution to the dual equation of the translator equation.  Then,  the result on growth rate   of the Legendre dual in \cite[Theorem 1.2 and Remark 1.3]{CCK_regularity} and the standard interior regularity theory imply the smooth convergence of the duals with locally uniform convexity bounds.

Now, we provide a key lemma in this subsection.
\begin{lemma}\label{lem-para.bound}
If a sequence $\{\Sigma_i\}_{i=1}^\infty\subset \mathcal{X}_k$ converges in $C^0_{\text{loc}}$ to $\Sigma \in \mathcal{X}_k$, then $\{\psi(\Sigma_i)\}_{i=1}^\infty$ is bounded in $\mathbb{R}^K\times \mathcal{S}_k$.
\end{lemma}
By assuming this lemma, we can prove the last main theorem.
\begin{proof}[Proof of Theorem \ref{thm:moduli.space.topology}]  We have already established that $\psi$ is a bijection by the classification theorem. Hence, it suffices to prove the continuity of $\psi$ and $\psi^{-1}$. This reduces down into two statements: 
\begin{enumerate}
\item if $\{(\mathbf{a}_i,h_i)\}_{i\in \mathbb{N}}\subset \mathbb{R}^K\times \mathcal{S}_k$ converges to $(\mathbf{a},h)$, then  $\{\Sigma_{\mathbf{a}_i}^{h_i}\}_{i\in \mathbb{N}}\subset \mathcal{X}_k$ converges to $\Sigma_{\mathbf{a}}^h\in \mathcal{X}_k$ in $C^0_{\text{loc}}$,
\item if $\{\Sigma_i\}_{i\in \mathbb{N}}\subset \mathcal{X}_k$ converges to $\Sigma \in \mathcal{X}_k$ in $C^0_{\text{loc}}$, then $\{\psi(\Sigma_i)\}_{i\in \mathbb{N}}$ converges to $\psi(\Sigma)$.
\end{enumerate}

 The most part of (1) has been already established in the continuity statement of Theorem \ref{thm-existence} in Section \ref{subsec:existenceconstruction}, part (iv). There is an additional complication due to the possible rotations brought by $h_i$, but the corotational construction rules out this issue. Note that, for entire graphical surfaces, the locally uniform convergence is the same as the locally uniform convergence of height functions. Let $\Sigma^{h_i}_{\mathbf{a}_i}$ be the graph of $u_{(\mathbf{a}_i,h_i)}(x)$ on $\mathbb{R}^2$. On every compact ball $\bar B_r$, observe 
\[ \Vert u_{(\mathbf{a}_i,h_i)}-u_{(\mathbf{a},h)}\Vert_{C^0(\bar B_r)} \le \Vert u_{(\mathbf{a}_i,h_i)}-u_{(\mathbf{a},h_i)}\Vert_{C^0(\bar B_r)}+\Vert u_{(\mathbf{a},h_i)}-u_{(\mathbf{a},h)}\Vert_{C^0(\bar B_r)} . \]
In view of  the corotational construction,  
$\Vert u_{(\mathbf{a}_i,h_i)}-u_{(\mathbf{a},h_i)}\Vert_{C^0(\bar B_r)} $ is the same as $\Vert u_{(\mathbf{a}_i,h)}-u_{(\mathbf{a},h)}\Vert_{C^0(\bar B_r)}$ 
 which 
 converges to $0$   by Theorem \ref{thm-existence}, part (iv). Similarly,  $u_{(\mathbf{a}_i,h_i)}$ is a rotation of $u_{(\mathbf{a}_i,h)}$ in $x_1x_2$-plane by an angle, say $\theta_i$, and $\theta_i \to 0$ as $i\to\infty$. Thus, $\Vert u_{(\mathbf{a},h_i)}-u_{(\mathbf{a},h)}\Vert_{C^0(\bar B_r)} $ also converges to $ 0$. This proves (1). 

\medskip 

 Let us show the part (2). For the sake of contradiction, we assume $\{\Sigma_i\}$ has a subsequence $\{\Sigma_{i_m}\}$ such that $|\Psi(\Sigma_{i_m})-\Psi(\Sigma)| \geq \varepsilon$ for some $\varepsilon>0$.  Remembering Lemma \ref{lem-para.bound}, however, the Bolzano-Weierstrass theorem applied to $\Psi(\Sigma_{i_m})$, and then the part (1) imply that  there is  a   further  subsequence, still denoted $\Sigma_{i_m}$, converging to some $\overline{\Sigma}\in \mathcal{X}_k$ and satisfying   $ \Psi(\Sigma_{i_m})\to \Psi(\overline\Sigma)$. Therefore, we have $|\Psi(\overline{\Sigma})-\Psi(\Sigma)| \geq \varepsilon$ which contradicts   the assumption that $\Sigma_i \to \Sigma$ in $C^0_{\text{loc}}$.
\end{proof}

\medskip

In the remaining of this section, we prove Lemma \ref{lem-para.bound}.  Recall the transformation $\mathbf{b}:\mathbb{R}^K\to \mathbb{R}^K$ from Definition \ref{def-b(a)}
 \begin{equation}
     \mathbf{b}(a_0,\cdots,a_{K-1})=(\textrm{sgn}(a_0)|a_0|^{1/({\sigma-\beta^+_0})},\ldots , \textrm{sgn}(a_{K-1})|a_{K-1}|^{1/({\sigma-\beta^+_{K-1}})}).
 \end{equation}
For $\Sigma _{\mathbf{a}}^h \in \mathcal{X}_k$ and  $\lambda>0$, we define the rescaled surface
 \begin{equation}\label{eq-rescale.down}
     \Sigma_{\mathbf{a},\lambda}^h:= \{(x_1,x_2,x_3)\in \mathbb{R}^3\,:\, (\lambda^\sigma x_1, \lambda ^\sigma x_2, \lambda x_3)\in \Sigma _{\mathbf{a}}^h\}. 
 \end{equation}
  Recalling the difference function $w_{\mathbf{a}}^h(s,\theta):=S_{\mathbf{a}}^h(e^s,\theta)-\sigma^{-1}e^{\sigma s}h(\theta)$, we observe that $\Sigma _{\mathbf{a},\lambda}^h$ has the support function $S_{\mathbf{a},\lambda}^h$ and the difference $w_{\mathbf{a},\lambda}^h$ satisfying
 \begin{equation}
S_{\mathbf{a},\lambda}^h(l,\theta) := \lambda ^{-\sigma } S_{\mathbf{a}}^h(\lambda l,\theta)\quad \text{ and } \quad w_{\mathbf{a},\lambda}^h(s,\theta) := \lambda^{-\sigma } w_{\mathbf{a}}^h(s + \ln \lambda , \theta).     
 \end{equation}
 
\bigskip

\begin{lemma}\label{lem:nontrivial.blow.down}
Given a sequence $\{(\mathbf{a}_i,h_i)\}_{i\in \mathbb{N}}\subset \mathbb{R}^K\times \mathcal{S}_k$ with $|\mathbf{a}_i|\to +\infty$, the rescaled surfaces $\Sigma_{\mathbf{a}_i,|\mathbf{b}(\mathbf{a}_i)|}^{h_i}$ subconverge  in $C^0_{\text{loc}}$ to a blow-down translator $\Sigma_{\infty}$, where  the   support function $S_{\infty}$ of $\Sigma_{\infty}$  satisfies $\sigma l^{-\sigma} S_{\infty}(l,\cdot)\to  h_{\infty}$  as $l\to +\infty$  for some $ h_{\infty} \in \mathcal{S}_k$, but $S_{\infty} \not \equiv \sigma^{-1}l^\sigma h_{\infty}$.
\end{lemma}

\begin{proof}
By  Corollary \ref{cor-wexpression} in Section \ref{subsec:existenceconstruction}, for every $\mathbf{a}=(a_0,\ldots, a_{k-1})\in \mathbb{R}^K$ and $h\in \mathcal{S}_k$, we can write $w_{\mathbf{a}}^h$ as 
\begin{equation}
    \label{eq-wexpression} w_{\mathbf{a}}^h= g_{\mathbf{0}}^h+\sum_{j=0}^{K-1} a_j ( \varphi_j^h e^{\beta^+_j s} + { g_{\pi_j(\mathbf{a})}^{h}})+f_\mathbf{a}^h
\end{equation}
with  the families of functions $\{g_{\mathbf{a}}^{h}\}$ and  $\{f_\mathbf{a}^h\}$. Here, given $\beta \in (0,1)$,  there are   $\rho >0$  and $\varepsilon>0$ (only depending on $\alpha$ and $\beta $) satisfying
\begin{align}\label{eq-gf}
\| g_{\mathbf{0}} \|_{C^{2,\beta,\sigma -\eps}_{R_{\rho,\mathbf{0}}}}\ \le    e^{\eps  R_{\rho,\mathbf{0}}}, \quad  \|g_{\pi_j(\mathbf{a})}^{h} \|_{C^{2,\beta,\beta^+_{j}-\eps}_{R_{\rho,\mathbf{a}}}}\ \le    e^{\eps  R_{\rho,\mathbf{a}}}  ,\quad  \Vert f_\mathbf{a}^h \Vert _{C^{2,\beta,-\frac{1}{2} }_{R_{\rho,\mathbf{a}}}} \le e^{(\sigma +\frac{1}{2} )R_{\rho,\mathbf{a}}},  
\end{align}
where $R_{\rho,\mathbf{a}}:=\rho+\ln (|\mathbf{b}(\mathbf{a})|+1)$. We rescale \eqref{eq-wexpression} by $|\mathbf{b}|:=|\mathbf{b}(\mathbf{a})|$ so that
\begin{equation}\label{eq-wab} 
w_{\mathbf{a},|\mathbf{b}|}^h ={ g_{\mathbf{0},\vert \mathbf{b}\vert}^h} +\sum_{j=0}^{K-1} a_j |\mathbf{b}|^{\beta^+_j - \sigma } (\varphi_j^h e^{\beta^+_j s} +  g_{\pi_j(\mathbf{a}),|\mathbf{b}|}^{h} ) + f_{\mathbf{a},|\mathbf{b}|}^h,
\end{equation}
     where  $g_{\mathbf{a},|\mathbf{b}|}^{h}(s,\theta)= |\mathbf{b}|^{-\sigma} g_{\mathbf{0}}^{h}(s+\ln \vert \mathbf{b}\vert,\theta )$,  $g_{\pi_j(\mathbf{a}),|\mathbf{b}|}^{h}(s,\theta)= |\mathbf{b}|^{-\beta^+_j} g_{\pi_j(\mathbf{a})}^{h}(s+\ln \vert \mathbf{b}\vert,\theta )$ and $f_{\mathbf{a},|\mathbf{b}|}^{h}(s,\theta)= |\mathbf{b}|^{-\sigma} f_{\mathbf{a}}^{h}(s+\ln \vert \mathbf{b}\vert,\theta )$. Hence, by \eqref{eq-gf}, it holds that 
 \begin{align}
     \| f_{\mathbf{a},|\mathbf{b}|}^{h}\|_{C_{s,1}^{2,\beta}}=|\mathbf{b}|^{-\sigma}\| f_{\mathbf{a}}^{h}\|_{C_{s+\ln |\mathbf{b}|,1}^{2,\beta}}\leq |\mathbf{b}|^{-\sigma}e^{(\sigma +\frac{1}{2} )R_{\mathbf{a}}}e^{-\frac{1}{2}(s+\ln |\mathbf{b}|)}=C_0e^{-\frac{1}{2}s}
 \end{align}
  for $s \geq 1+\rho+\ln(1+|\mathbf{b}|^{-1})$, where $\ln C_0=(\sigma+\frac{1}{2})(\rho+\ln (1+|\mathbf{b}|^{-1}))$. Similarly, we   have 
\begin{equation}
\Vert g_{ \pi_j(\mathbf{a})}^{h} \Vert _{C^{2,\beta, \beta^+_j-\eps}_{\rho+\ln(1+|\mathbf{b}|^{-1})}}\le  C \quad \text{ and } \quad \Vert f_{\mathbf{a},|\mathbf{b}|}^h \Vert _{C^{2,\beta,-1/2}_{\rho+\ln(1+|\mathbf{b}|^{-1})}}  \le C .
\end{equation} 

Now, given a sequence $\{(\mathbf{a}_i,h_i)\}_{i=1}^\infty$ with $|\mathbf{a}_i|\to +\infty$, by passing to a subsequence, we may assume $(\mathbf{b}_i/|\mathbf{b}_i|,h_i)\to ( \mathbf{b}_{\infty},  h_{\infty})$, where $\mathbf{b}_i:=\mathbf{b}(\mathbf{a}_i)$. Then, we plug in $(\mathbf{a}_i,\mathbf{b}_i)$ into $(\mathbf{a},\mathbf{b})$ at \eqref{eq-wab}, and we observe
\begin{align}
&    \lim_{i\to +\infty} \varphi_j^{h_i}=\varphi_j^{h_\infty}, && \lim_{i\to +\infty} a_{i,j} |\mathbf{b}_i|^{\beta^+_j - \sigma }=    b_{\infty,j} | b_{\infty,j}|^{ -1-(\beta^+_j-\sigma) }
\end{align}
where $\mathbf{a}_i:=(a_{i, 0},\ldots, a_{i, k-1})$ and $ \mathbf{b}_{\infty}=(b_{\infty,0} ,\ldots, b_{\infty,K-1})$. Also,  $g_{\mathbf{0},\vert \mathbf{b}_i\vert }^{h_i}$ smoothly converges to $0$ on every compact set of $(s,\theta)$ as $i\to +\infty$, because $\Sigma_{\mathbf{0}}^{h_i}$ is asymptotic to the graph of the homogeneous profile. Therefore,  $w_{\mathbf{a}_i,|\mathbf{b}_i|}^{h_i}$ subconverge to   $w_\infty  = \sum_{j=0}^{J}  b_{\infty,j} | b_{\infty,j}|^{-1-(\beta^+_j-\sigma) }  \varphi_i e^{\beta^+_i s} + O(e^{(\beta^+_J -\eps)s })  $ as $s\to \infty$, where $J$ is the largest index such that $b_{\infty,J}\neq0$.    This shows that  $\sigma l^{-\sigma}S_\infty \to h_\infty$ but $\sigma l^{-\sigma}S_\infty \not \equiv h_\infty$. 

By the subsequential convergence of $w_{\mathbf{a}_i,|\mathbf{b}_i|}^{h_i}$, the rescaled surfaces $\Sigma_{\mathbf{a}_i,|\mathbf{b}_i |}^{h_i}$ subconverge    in $C^0_{\text{loc}}(\{l>e^\rho\})$.  Moreover, by the convexity, we can also conclude that $\Sigma_{\mathbf{a}_i,|\mathbf{b}_i |}^{h_i}$ subconverge to   a convex hypersurface $\Sigma_{\infty}$ in $C^0_{\text{loc}}(\mathbb{R}^3)$, which is a  blow-down translator.
\end{proof}
 
 If we further assume the original translators $\Sigma_{\mathbf{a}_i}^{h_i}$ are also convergent, then the rescaled limit $\Sigma_\infty$ has homogeneous growth rate near the origin.  

\begin{proposition}\label{prop:blow.donw.origin.asymp}
Suppose that $\{\Sigma_{\mathbf{a}_i}^{h_i}\}_{i=1}^\infty \subset\mathcal{X}_k$ is convergent, and $\{\Sigma_{\mathbf{a}_i,\lambda_i}^{h_i}\}_{i=1}^\infty$ with $\lambda_i\to +\infty$ is also convergent in $C^0_{\text{loc}}(\{x_3>0\})$. Then, $\{\Sigma_{\mathbf{a}_i,\lambda_i}^{h_i}\}_{i=1}^\infty$ converges in $C^0_{\text{loc}}(\mathbb{R}^3)$ to a blow-down translator $\Sigma_\infty=\partial \{x_3>u_{\infty}(x_1,x_2)\}$ with $u_\infty \in C^1(\mathbb{R}^2)\cap C^\infty(\mathbb{R}^2\setminus \{0\})$ satisfying $u_\infty(0)=0$ and $\nabla u_\infty(0)=0$. Moreover,  $J^*_{u_\infty}(x)\to J^*_j$ for some $j$ as $x\to 0$.     
\end{proposition}

\begin{proof} Since $\Sigma^{h_i}_{\mathbf{a}_i}$ are convergent, we may assume that their height functions $u_i$ satisfy $u_i(0)=0$ and $\nabla u_i(0)=0$.  Let us denote by $S^i_l$ the section  of $u_i$ at height $l$.  
Recalling  \cite[Definition 3.6]{CCK_regularity} (or\cite[pp 644]{DSavin}),  we write   $S^i_l\sim A$ for unimodular matrix $A$ when the eccentricity of the section  $S^i_l$ is proportional to $|A|$. We also recall   the result  \cite{CCK_regularity}  on the growth rate of $u_i$ and its Legendre dual $u_i^*$, and    \cite[{Corollary 1.4}]{CCK_regularity}     on the behavior of eccentricities of level sets with the constants $l_0,C_0,C_1$.  
Since $\Sigma_{\mathbf{a}_i}^{h_i}$ are convergent, sections $S^i_{l_0}$ are also convergent. Hence, there is some $M\geq C_0$ (depending on the limiting translator) such that   if $S^i_{l_0} \sim A^i_{l_0}$ for a matrix $A^i_{l_0}$, then $|A^i_{l_0}|\leq M$ for sufficiently large $i$. Therefore, there are some large $i_0$, $C>1$ and   $R >1$ (independent of $i$) such that the dual function $u^*_i$ of $u_i$ satisfies
\begin{equation}\label{eq-BcC} 
0< C^{-1} \le |x|^{-\frac1{2\alpha}}u^*_i(x)\le C
\end{equation}
for $i\geq i_0$ and  $|x| \geq R$. Thus, $\Sigma_{\mathbf{a}_i,\lambda_i}^{h_i}\to \Sigma_{\infty}$ in $C^0_{\text{loc}}$ and $C^{-1}|x|^{\frac{1}{2\alpha}} \leq u^*_{\infty}(x)\leq C|x|^{\frac{1}{2\alpha}}$ holds in $\mathbb{R}^2$. This implies $u_\infty \in C^1(\mathbb{R}^2)$, $u_\infty(0)=0$, and $\nabla u_\infty(0)=0$. Finally, by \cite[Corollary 2.12]{MA-figalli-book}, $u^*_\infty$ is an Alexanrov solution to \eqref{homo.eq}, and therefore $\Sigma_\infty$ is a blow-down translator with profile $u_\infty \in C^\infty(\mathbb{R}^2\setminus \{0\})$. Finally,  the convergence of $J^*_{u_{\infty}}(x)$  was already shown in \cite[Proposition 5.7]{DSavin}.
\end{proof}

\bigskip

\begin{proposition}\label{prop-276}Let $v\in C^1(\mathbb{R}^2)\cap C^{\infty}(\mathbb{R}^2\setminus\{0\})$ be a solution to \eqref{homo.eq} in $\mathbb{R}^2\setminus\{0\}$. Suppose $v(0)=0$ and $\nabla v(0)=0$ hold, and the entropy $J_v(x)$ has the limits
\begin{equation}
\lim_{x\to 0} J_v(x) = J_j\quad\text{ and }\quad \lim_{x\to \infty } J_v(x) = J_k ,
\end{equation}
where $J_j$ and $J_k$ are entropies for $j$-fold and $k$-fold homogeneous solutions, respectively.  If $j \ge  k$, then $j=k$ and $v$ has to a $k$-fold homogeneous solution.      
\end{proposition}

\begin{proof}
If $j=k=\infty$, then $v$ is the rotationally symmetric homogeneous solution by Proposition \ref{prop_homo.int.max} and \ref{prop_homo.constnat}. Thus, we assume $k<\infty$. By Theorem \ref{thm_Entorpy order} and $j\ge k$, we have $J_j\leq J_k$ and $J_\infty<J_k$. Hence, Proposition \ref{prop_homo.int.max} implies $J_v(x) \leq J_k$ for every $x\in \mathbb{R}^2$. Therefore, Theorem \ref{thm-rigiditywithcondition} completes the proof.
\end{proof}

\bigskip

\begin{proof}[Proof of Lemma  \ref{lem-para.bound}]
Suppose that  $\{\Sigma_{\mathbf{a}_i}^{h_i}\}_{i=1}^\infty\subset \mathcal{X}_k$ converges to $\Sigma \in \mathcal{X}_k$, but $|\mathbf{a}_i|\to +\infty$. Then, by Lemma \ref{lem:nontrivial.blow.down}, $\Sigma_{\mathbf{a}_i,|\mathbf{b}_i|}^{h_i}$ subconverge to a blow-down translator $\Sigma_\infty=\partial\{x_3>u_\infty(x_1,x_2)\}$ whose profile is not a homogeneous function. Thus, by Proposition \ref{prop:blow.donw.origin.asymp} and \ref{prop-276}, we have
\begin{equation}
\lim_{x\to 0} J_{u_\infty}(x) = J_j\quad\text{ and }\quad \lim_{x\to \infty } J_{u_\infty}(x) = J_k
\end{equation}
for some $j<k$, namely $J_j>J_k$. For simplicity, we may assume  $\Sigma_{\mathbf{a}_i,|\mathbf{b}_i|}^{h_i} \to \Sigma_\infty$, and denote the profile of $\Sigma_{\mathbf{a}_i,|\mathbf{b}_i|}^{h_i}$ by $u_i$.

Remembering Theorem \ref{thm_Entorpy order}, we choose a small $ \delta = \frac{1}{4} \min (  J_m -J_\infty ,\min_{\ell \in \{3,\ldots,m-1\} }  (J_{\ell }-J_{\ell +1}))>0$. Since $\Sigma = \partial \{x_3 > u(x)\}\in \mathcal{X}_k$, i.e. the blow-down is $k$-fold symmetric, there exists $R>0$ such that 
\begin{equation}\label{eq-entropy.asympt}
    J_{u^*}(x)\le J_k +\delta
\end{equation}
holds for all $|x|  \ge R $. Since  $\lim_{x\to 0^+}J_{u^*_\infty}(x)=J_j>J_k$,  $\sup_{x\in \mathbb{R}^2 \setminus B_R}J_{u^*_i}(x)$ is achieved at a point $x_i\in \mathbb{R}^2\setminus B_R$ for sufficiently large $i$ and $\liminf _{i\to \infty}J_{u_i^*}(x_i)\ge J_j$. Next, we claim $|x_i|\to +\infty$. If not, $x_i$ subconverge to a certain $\bar x \in \mathbb{R}^2\setminus B_R$ and $J_{u^*}(\bar x)\geq J_j$. This contradicts \eqref{eq-entropy.asympt}. Now, by \eqref{eq-BcC} and by passing to a subsequence, we may assume $z_i= x_i/|x_i|\to z\in \partial B_1$ and $v_i(x):=|x_i|^{-\frac{1}{2\alpha}}u^*_i(|x_i|x)$ converges to an Alexandrov solution $ v$ to \eqref{homo.eq} with $ v(0)=0$ and $\nabla   v(0)=0$. Since  $J_{ v_j}(x)$ on $\{x\, :\, |x|\ge R/|x_i|\}$ attains its maximum at $z_j\in  \partial B_1$ and $v_j$ converges to $v$ in $C^{\infty}_{\mathrm{loc}}(\mathbb{R}^2\setminus \{0\})$, we conclude that $J_v(x)$  on $  \mathbb{R}^2\setminus \{0\}$ attains its maximum at $z$. Moreover, $J_v(z)=\lim _{j\to \infty}J_{ v_j}(z_j)\ge J_j >J_k \ge J_{\infty} $. Therefore,  by Proposition \ref{prop_homo.int.max} and \ref{prop_homo.constnat}, $  v$ is a $j'$-fold homogeneous solution for some $j'\le j$. From this result, we conclude that for large $i$, there is $y_i$ such that 
\begin{align}
|y_{i}|<|x_{i}|, \qquad  J_{u^*_i}(y_{i})=J_{j'}-2\delta, \qquad  J_{u^*_i}(x) \geq J_{j'}-2\delta \quad\hbox{on $\{ x\,:\, |y_{i}|\leq |x| \leq |x_{i}|\}$}. 
\end{align}

We observe that $|y_i|\to \infty$ for the same reason as $|x_i|\to \infty$, and $|x_i|/|y_i|\to \infty$ as $i\to \infty$ due to the convergence of $|x_i|^{-\frac{1}{2\alpha}}u^*_i(|x_i|x)$ to the $j'$-fold solution $ v$. Using \eqref{eq-BcC} again, $|y_i|^{-\frac{1}{2\alpha}}u^*_i(|y_i|x)$ subsequentially converges to an Alexandrov solution $\check v$ satisfying $J_{\check v}(x)\le J_{j'}= \lim_{i\to \infty} (\sup_{x\in \mathbb{R}^2 \setminus B_R}J_{u^*_i}(x))$ on $x\in \mathbb{R}^2\setminus \{0\}$. Moreover, we have $\lim_{|x|\to \infty}J_{\check  v}(x)= J_{j'}$ since $J_{j'}-2\delta \le J_{\hat v}(x)\le J_{j'}$ on $|x|\ge 1$ (here we used $|x_i|/|y_i|\to \infty$) and hence the blow-down of $\hat v$ has to be $j'$-fold symmetric by the choice of $\delta$. 

Finally, from Theorem \ref{thm-rigiditywithcondition}, we conclude $J_{\check v}\equiv J_{j'}$.  This contradicts $\min_{|x|=1} J_{\check v}(x)= J_{j'}-2\delta$.
\end{proof}

\appendix 
\section{ODE lemmas} \label{sec-ODE}


\begin{lemma} [Merle-Zaag ODE lemma A.1 in \cite{merle1998optimal}] \label{lem-MZODE}
	Suppose $x(s)$, $y(s)$, and $z(s)$ be non-negative absolutely continuous functions such that  $x+y+z >0$ and $\lim_{s\to +\infty} x + y + z = 0$. Suppose for each small $\e>0$ there is $s_0$ such that for each $s\ge s_0$  
	\begin{align}\label{eq:mz.ode.system}
&	x' - x \geq  -\eps(y + z), &&	|y'| \leq \eps (x + y + z), && z' + z \leq  \eps(x + y),
	\end{align}
	holds almost every $s\geq s_0$. Then either one of $x+z = o(y)$ or $x+y = o(z)$ holds as $s\to +\infty$.
\end{lemma}

The following elementary ODE lemmas are excerpted from \cite{CS2}.

\begin{lemma}[Lemma B.5 in \cite{CS2}]\label{lem:ODE.small.O} Let $\tilde{\rho},f:[s_0,+\infty)\to \mathbb{R}$ be absolutely continuous functions satisfying
\begin{align}\label{eq:fr.ODE.condition}
& \lim_{s\to +\infty} f(s)=0, && \lim_{s\to +\infty} \tilde{\rho}(s) e^{\frac{|\lambda|}{2}s}=+\infty, &&   |\tilde{\rho}'| \leq \tfrac{1}{2}|\lambda|\tilde{\rho},
\end{align}
and $\lambda^{-1}f'\geq f +o(\tilde{\rho})$ for some $\lambda\neq 0$. Then, $f_+=\max\{0,f\} =o(\tilde{\rho})$ holds.
\end{lemma}

 \medskip

\begin{lemma}[Lemma B.6 in \cite{CS2}]\label{lem:ODE.large.O} Let $\tilde{\rho},f:[s_0,+\infty)\to \mathbb{R}$ be absolutely continuous functions satisfying  \eqref{eq:fr.ODE.condition} and
$\lambda^{-1}f'=f+O(\tilde{\rho})$ for some $\lambda\neq 0$. Then, we have $f=O(\tilde{\rho})$. 
\end{lemma}

 \medskip

\section{Order of entropies}

In this section, we study entropies    of homogeneous solutions to 
\be \label{eq-wmain}
 \det D^2 v = |x|^{\frac 1\alpha -4}\quad\text{ on }\mathbb{R}^2. 
 \ee 
For $\alpha<  \frac{1}{2}, $ we say that an entire solution to the equation
\eqref{eq-wmain}  is homogeneous if $v$ is of the form $r^{\frac{1}{2\alpha}}g(\theta)$ in polar coordinates. By a direct computation, $g(\cdot)>0$ satisfies the ODE 
\be \label{eq-odeg} 
\tfrac g{2\alpha}\left(g'' + \tfrac{g}{2\alpha} \right) - \left(\tfrac1{2\alpha} -1\right) (g')^2 = \tfrac{2\alpha}{1-2\alpha}.
\ee 
If $g $ is not a constant, then $g $ has to be a periodic solution   of  period $ 2\pi/k$ for some integer  $k\ge1$. If there is no larger such $k$, we call it a $k$-fold symmetric solution.  Daskalopoulos and Savin   \cite{DSavin}  defined a function $J_v$ associated with a strictly convex function $v$, which behaves like an entropy in many aspects  \be \label{eq-J}
 J_v(x) = \Delta v (r^2 v_{rr})^{4\alpha-1}.
 \ee
 First, the entropies $J_v$ are constants for homogeneous solutions. This can be shown by the same result for dual equation (Remark \ref{rmk-homo-J}), but we prove here as we will need the derivation of $J_v$ in its proof.

\begin{lemma}[Remark 5.3 of  \cite{DSavin}] 
\label{lem-Jconstant} If $v $ is a  homogeneous    solution to \eqref{eq-wmain}, then $J_v $ is a constant function.
	\begin{proof} 
Let $v=r^{\frac{1}{2\alpha}}g(\theta)$ for some   $g$   on $\mathbb{S}^1$. By a direct computation, we have 
		\bea \label{eq-J_w-g}J_v (x)= \left(  g'' + \tfrac{g}{4\alpha^2} \right) \left( \tfrac{1-2\alpha}{4\alpha^2} g\right)^{4\alpha-1}.\eea 
	Since $J_v$ is trivially constant for the rotationally symmetric case, we may assume that  $v$ is not rotationally symmetric and hence $g$ is non-constant. For such   $g$, let us choose $[a,b]$, a maximal increasing interval. By regarding $g=g(\theta)$ on $\theta\in[a,b]$ as a new variable, we define   $f:[g(a),g(b)] \to (0,\infty)$  by  the relation 
\be
g'=\sqrt{2f(g)}.
\ee
  Then $f$ is a solution to the ODE 
  \be\label{eq-odeh}
  \tfrac t{2\alpha}\left(f'(t) + \tfrac{t}{2\alpha} \right)-2 \left(\tfrac1{2\alpha} -1\right) f(t) = \tfrac{2\alpha}{1-2\alpha}.
  \ee 
By solving this equation, we obtain  
\bea\label{eq-hc} 
2f(t)= 2f_c(t)=c\, t^{2(1-2\alpha)}- \tfrac{1}{4\alpha^2}t^2 - (\tfrac{2\alpha}{1-2\alpha})^2,
  \eea  
  where   a  constant $c>0$ is    chosen so that   $\{t\,:\, f_c(t)> 0 \}=(g(a),g(b))$.  In terms of $f_c(\cdot)$, we can rewrite \eqref{eq-J_w-g} as 
\be\label{eq-hh}
J_v=  \left( f_c'(t)+\tfrac{t}{4\alpha^2}\right)\left( \tfrac{1-2\alpha}{4\alpha^2}t \right)^{4\alpha-1}  = c \, (1-2\alpha)\left( \tfrac{1-2\alpha}{4\alpha^2} \right)^{4\alpha-1}.
\ee 
In particular, this shows that   $J_v(\cdot)$ is constant.
	\end{proof}
\end{lemma} 
Next, if $v= |x|^{\frac{1}{2\alpha}}g(\theta)$ is a homogeneous solution to \eqref{eq-wmain}, then the Legendre dual $v^*$ is also a  homogeneous function (of degree $\frac{1}{1-2\alpha}$),  which solves $\det  D^2 v^* = |Dv^*|^{4-\frac{1}{\alpha}}$. From the classification of such solutions $v^*$ in Corollary \ref{cor-anderews2003}, we infer the classification of homogeneous solutions to \eqref{eq-wmain}.
\begin{theorem}[Classification of homogeneous solutions]\label{thm-class-homo-dual}
Given $\alpha \in (0,1/4)$, let $m$ be the largest integer such that $m^2< 1/\alpha$. Then, modulo rotations, there are exactly $(m-1)$ homogeneous solutions to \eqref{eq-wmain}. These are the rotationally symmetric solution and $k$-fold symmetric solutions for  $3\leq k\leq m$  (which exist only if $\alpha<\frac19$).

\end{theorem}

In what follows, we will use the subscript $_\infty$ to denote   quantities that correspond to the rotationally symmetric solution. 

 \begin{definition}[Entropy of homogeneous solution]\label{def-entrok} For a given     $\alpha \in (0,\frac19)$  and    any   integer $3\le k< \alpha^{-1/2}$, we define the constant $J_k$ to be the entropy $J_{v_k}$ of a $k$-fold symmetric homogeneous solution $v_k$ to \eqref{eq-wmain}. Similarly, we define $J_\infty$ to be the entropy of the rotationally symmetric homogeneous solution to \eqref{eq-wmain}. 
\end{definition}

We will show  that $J_k$ are strictly decreasing with respect to $k$ in Theorem \ref{thm_Entorpy order}.  
\begin{lemma}\label{lem-J0} There holds that $J_\infty<J_k$ for any integer $3\le k < \alpha^{-1/2}$. 

\begin{proof}
Observe that there exists a unique $c_{\infty}>0$ which makes $f_{c_\infty}$ in \eqref{eq-hc} touches the zero function from below, i.e., $f_{c_\infty}(t)\le 0$ for all $t>0$ and there is a unique $t_{\infty}>0$ such that $f_{c_\infty}(t_\infty)=0$.  
Using  \eqref{eq-odeh} and the  touching conditions: $f_{c_\infty}(t_\infty)=f'_{c_\infty}(t_\infty)=0$, we directly check that 
\be
t_\infty=2\alpha\sqrt{\tfrac {2\alpha}{1-2\alpha}}. 
\ee 
Since $g(\cdot)\equiv t_\infty  $ is a  constant solution to   \eqref{eq-odeg},  we have  that 
\be
J_\infty= { c_\infty} (1-2\alpha)\left( \tfrac{1-2\alpha}{4\alpha^2} \right)^{4\alpha-1}.
\ee 
For any other $k$-fold symmetric solutions, the corresponding $c$ in \eqref{eq-hc} should be strictly bigger than $c_{\infty}$ due to the choice of $c$. Thus we deduce that   $J_\infty <J_k$ by \eqref{eq-hh}. 
\end{proof}

\end{lemma}

\begin{theorem}\label{thm_Entorpy order}  For $\alpha \in (0,\frac 19)$, let $m\ge3$ be the largest integer such that $m^2< 1/\alpha$. Then the entropies of homogeneous solutions in Definition \ref{def-entrok} are strictly decreasing with respect to the fold number as follows: 
\be
J_{{  \infty} } < J_{ m} < J_{ {m-1}} < \ldots <J_{ 3}.
\ee

\begin{proof} 
  For each $3\le k\le m$,  let $ v_{k}=r^{1/2\alpha}g_{k}(\theta)$ be 
 a  homogeneous solution to  \eqref{eq-wmain}. 
First of all, the strict minimality of $J_{{  \infty} }$ is  shown in Lemma \ref{lem-J0}. 
As seen in the proof of Lemma \ref{lem-Jconstant},  $k$-fold symmetric homogeneous solutions $g_k$  for $3\le k <\alpha^{-1/2}$ are represented by $f_{c_k}(\cdot)$ for some $c_k>0$. From \eqref{eq-hh}, the proof reduces to show that $c_k$ is strictly decreasing with respect to $k$.   
\smallskip

  For each $3\le k\le m$,    let $v_k^*$ be    the Legendre dual of $v_k$. Then,   $v_k^*$  is a  $k$-fold symmetric homogeneous solution  to $   \det D^2 v^* = |Dv^*|^{4-\frac{1}{\alpha}}$, 
and    the level curves of $v^*_{k} $ are $k$-fold symmetric shrinking solitons  of the  $\frac{\alpha}{1-\alpha}$-CSF with $\frac{\alpha}{1-\alpha} <\frac13$;   see Remark \ref{remark-shrinker}. That is,   the level support function $S_k $ of  $v_k^*$  is given by $S_k(l,\theta)=\sigma^{-1} l^{\sigma} h_k(\theta)$ with some $k$-fold symmetric solution $h_k$ to \eqref{eq-shrinker}. 
Now, let us define the ratio between inner and outer radii  of these shrinkers by   
\be
 r_k :=  \frac{ \max h_k(\cdot) }{\min h_k(\cdot)}=  \frac{\max \left\{|x|\,:\,  v^*_{k}(x) =1\right\} }{\min \left\{|x|\,:\,  v^*_{k}(x) =1\right\}} .
\ee 
Then,  an important result of    \cite[Corollary 5.6]{andrews2003classification} shows that $r_k$ is strictly decreasing with respect to $k$.  
  Observe that  
 \bea
 v_k(e^{i\theta})&=\sup_{y\in\mathbb{R}^2} \left\{\langle e^{i\theta},y \rangle-  v_k^*(y) \right\}=\sup_{l\geq0  } \left\{ S_k(l,\theta)-  l \right\}=\sup_{l\geq0  } \left\{\sigma^{-1} l^{\sigma} h_k(\theta)-  l \right\}=\tfrac{2\alpha}{1-2\alpha}h_k^{\frac{1}{2\alpha}}(\theta).
 \eea
 In terms of $ v_{k}=|x|^{\frac1{2\alpha}}g_k(\theta)$, one can see  that 
\be
  \gamma_k :=  \frac{ \max g_k(\cdot) }{\min g_k(\cdot)}= \frac{\max \{   v_{k} (x)\,:\, |x|=1 \}}{\min \{  v_ {k} (x)\,:\, |x|=1 \}} =  r_k^{\frac{1}{2\alpha}}
\ee
 is strictly decreasing with respect to $k$.

 Recall that $\max g_k$ and $\min g_k$ are two zeros of $f_{c_k}(\cdot)=0$. 
 Since the two zeros are determined by the intersections of 
 \be
  y=c_k t^{2(1-2\alpha)}\quad\text{ and } \quad   y=\tfrac{1}{4\alpha^2}t^2 + (\tfrac{2\alpha}{1-2\alpha})^2,
  \ee 
 we conclude that $\gamma_k$  is strictly increasing with respect to $c_k$.  
This implies that $c_k$ strictly decreases with respect to $k$. This finishes the proof.
\end{proof}
\end{theorem}

\subsection*{Acknowledgments}
BC thanks KIAS for the visiting support and also thanks University of Toronto where the research was initiated. BC and KC were supported by  the National Research Foundation(NRF) grant funded by the Korea government(MSIT) (RS-2023-00219980). BC has been supported by NRF of Korea grant No. 2022R1C1C1013511, POSTECH Basic Science Research Institute grant No. 2021R1A6A1A10042944, and POSCO Science Fellowship. KC has been partially supported by KIAS Individual Grant MG078901 and POSCO Science Fellowship. SK has been partially supported by NRF of Korea grant No. 2021R1C1C1008184. 

We are grateful to Liming Sun for fruitful discussion. We would like to thank Simon Brendle, Panagiota Daskalopoulos, Pei-Ken Hung, Ki-ahm Lee, Christos Mantoulidis, Connor Mooney, Dongjun Noh for their interest and insightful comments.

 \bibliography{GCF-ref.bib}
\bibliographystyle{alpha}

\end{document}